\documentclass[11pt]{article}

\usepackage{a4wide}
\usepackage{amssymb}
\usepackage{amsmath}
\usepackage{theorem}
\usepackage{epsfig}
\usepackage{epstopdf}
\usepackage{verbatim}
\usepackage{graphicx}
\usepackage{subfigure}
\usepackage{color}
\usepackage{longtable}
\usepackage{pdflscape}
\usepackage{indentfirst}
\usepackage{rotating}
\usepackage{listings}
\usepackage{colortbl}
\newtheorem{theorem}{Theorem}[section]

\newtheorem{remark}[theorem]{Remark}

\newtheorem{conjecture}[theorem]{Conjecture}

\newcommand{\R}{\mathbb{R}}
\newcommand{\C}{\mathbb{C}}

\newcommand{\T}{\mathbb{T}}

\newcommand{\al}{\alpha}

\newcommand{\da}{\partial_{\alpha}}

\newcommand{\si}{\sigma}
\newcommand{\vp}{\varphi}
\newcommand{\om}{\omega}

\newcommand{\Dcal}{\mathcal{D}}

\allowdisplaybreaks[2]

\title{Structural stability for the splash singularities\\
 of the water waves problem}
\author{Angel Castro, Diego C\'ordoba, Charles Fefferman, \\ Francisco Gancedo, Javier G\'omez-Serrano}

\begin{document}

\maketitle
\begin{abstract}In this paper we show a structural stability result for water waves. The main motivation for this result is that we would like to exhibit a water wave
whose  interface starts as a graph and ends in a splash. Numerical simulations lead to an approximate solution with the desired behaviour. The stability result will conclude that near the approximate solution to water waves there is an exact solution.\end{abstract}

\section{Introduction}

The water waves problem models the motion of an incompressible fluid with constant density $\rho$  in a domain $\Omega(t)$ with a free boundary $\partial \Omega(t)$,
which satisfies the Euler equation with the presence of gravity and whose flow in potential.    The system, in $\R^2$, can be written, after some computations, as an equation for the free boundary, \begin{equation}\label{Parametriza}
\partial\Omega(t)=\{z(\alpha,t)=(z_1(\alpha,t),z_2(\alpha,t)):\alpha\in\R\},
\end{equation} and an equation for the amplitude of the vorticity, $\omega(\alpha,t)$, in the following way
\begin{equation}\label{em}
z_t(\al,t)=BR(z,\omega)(\al,t)+c(\al,t)z_{\al}(\al,t),
\end{equation}
\begin{align}
\begin{split}\label{cEuler}
\omega_t(\al,t)&=-2BR_t(z,\omega)(\al,t)\cdot
z_{\al}(\al,t)-\Big(\frac{\omega^2}{4|\da z|^2}\Big)_{\al}(\al,t) +(c\omega)_{\al}(\al,t)\\
&\quad+2c(\al,t) BR_{\al}(z,\omega)(\al,t)\cdot z_{\al}(\al,t)-2
(z_2)_\al(\al,t),
\end{split}
\end{align}
where $BR(z,\omega)$ is the classical Birkhoff-Rott integral
\begin{equation}\label{BR}
BR(z,\omega)(\al,t)=\frac{1}{2\pi}PV\int_{\R}\frac{(z(\al,t)-z(\beta,t))^{\bot}}{|z(\al,t)-z(\beta,t)|^2}\omega(\beta,t)d\beta.
\end{equation}
The function $c(\al,t)$ is arbitrary since the boundary is convected by the normal component of the velocity of the fluid. Also, we notice that, in order to get an explicit equation for $\partial_t\omega$, we need to invert the operator $$I + T = I + 2 \langle BR(z,\cdot), z_{\al} \rangle$$
and we have taken the acceleration due to gravity and the density $\rho$ equal to one.

Once one has solved this system for $(z,\omega)$ the velocity of the fluid and the pressure in the domain $\Omega(t)$ can be recovered by using Biot-Savart and Bernoulli laws. For details see \cite{Castro-Cordoba-Fefferman-Gancedo-GomezSerrano:finite-time-singularities-free-boundary-euler}.

In the last two decades these equations have been intensively studied. For an extensive survey about analytical results on water waves see the monograph \cite{Lannes:water-waves-book}.

\begin{figure}[h!]\centering
\includegraphics[scale=0.4]{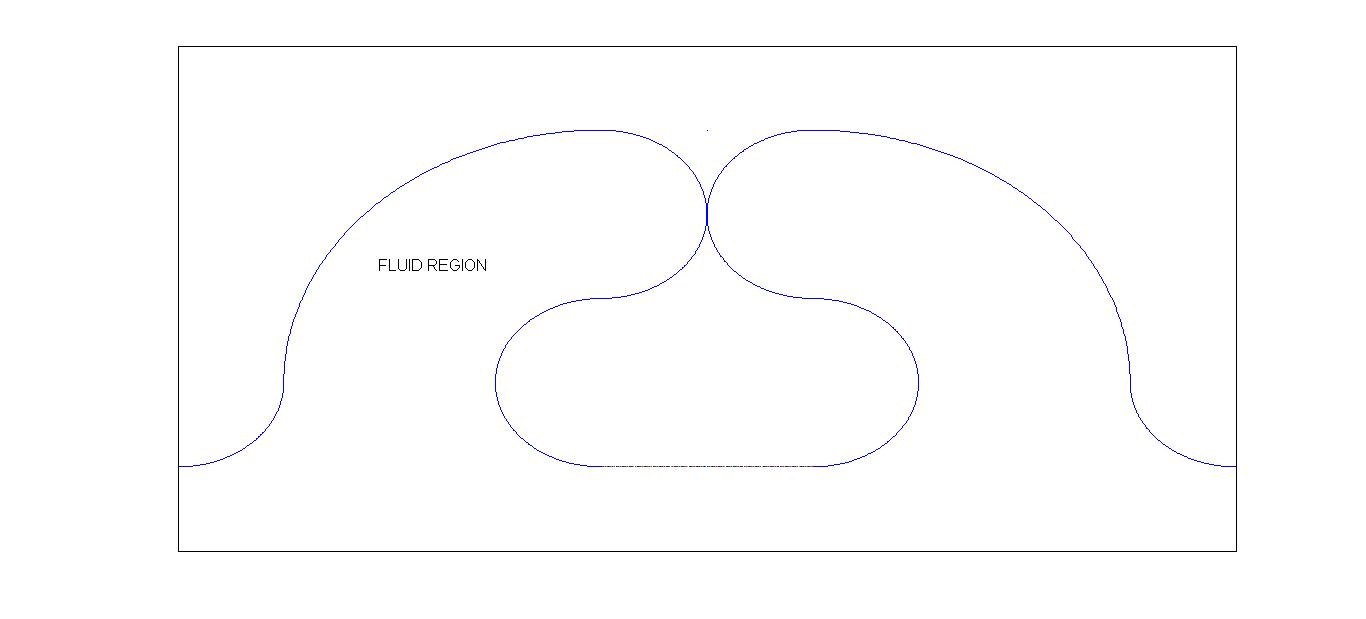}
\caption{Splash singularity. A smooth interface that  collapses in a point.}
\label{PictureSplash}
\end{figure}

In this paper we are concerned with the problem of the existence of water waves which start as a graph and become a splash curve in finite time. Roughly speaking, a splash curve is a smooth curve that  collapses with itself in a single  point such as the curve of fig. \ref{PictureSplash}. A rigorous definition can be found in \cite{Castro-Cordoba-Fefferman-Gancedo-GomezSerrano:finite-time-singularities-free-boundary-euler} where the existence of splash singularities has been shown. Coutand and Shkoller \cite{Coutand-Shkoller:finite-time-splash} have proven the existence of splash singularities in presence of vorticity. Fefferman, Ionescu and Lie \cite{Fefferman-Ionescu-Lie:absence-splash-singularities} have proven  the non existence of splash singularities for internal waves, i.e. for an interface between two incompressible fluids.

We are interested in the following statement:
\begin{conjecture}
There exist initial data $z_0(\al), \om_0(\al)$ of solutions of the water wave equations such that at time $0$ the curve  $z_0(\al)$ can be parameterized as a graph, the interface then turns over at a finite time $T_1 > 0$, and finally produces a splash at a finite time $T_2 > T_1$.
\end{conjecture}

\begin{figure}[h!]\centering
\includegraphics[scale=0.5]{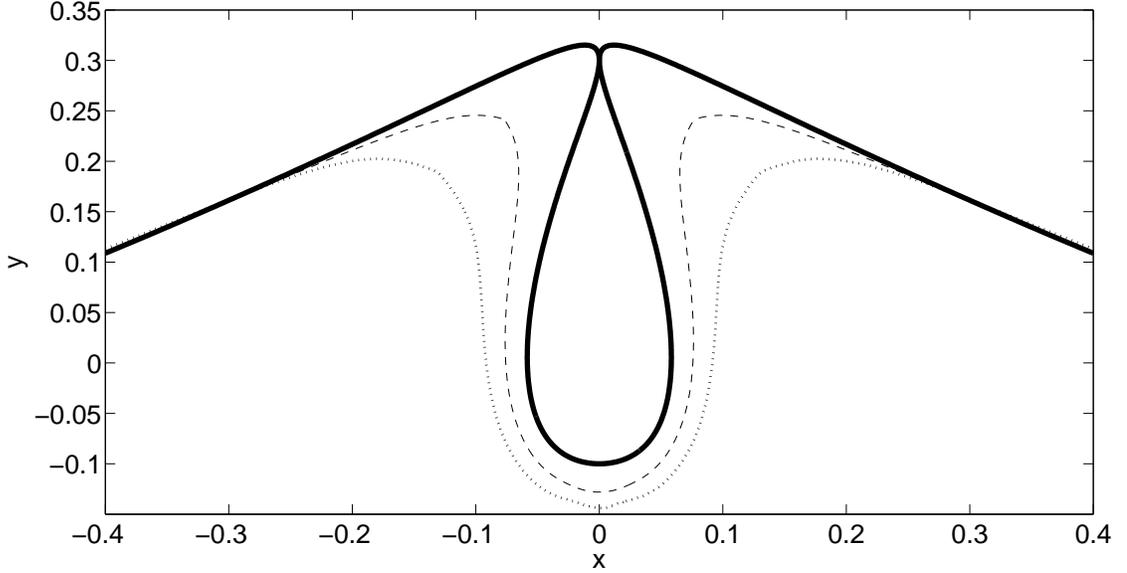}
\caption{Evolution from a graph to a splash.}
\label{splash}
\end{figure}

We should remark that this conjecture is a combination of the scenarios in theorems \cite[Theorem I.1] {Castro-Cordoba-Fefferman-Gancedo-GomezSerrano:finite-time-singularities-free-boundary-euler} and \cite[Theorem 7.1]{Castro-Cordoba-Fefferman-Gancedo-LopezFernandez:rayleigh-taylor-breakdown} and is supported by numerical evidence that we can see in Fig. \ref{splash}. This numerical simulation was carried out using the method of Beale, Hou and Lowengrub \cite{Beale-Hou-Lowengrub:convergence-boundary-integral}.

The proof of this conjecture could follow along these lines. First of all, we will move backwards in time, 0 being  the time of the splash, $T_2 - T_1$ the time of the turning and $T_2$ the time in which the solution can be parameterized as a graph. Also we write the water waves equation in a new  domain given by the projection of $\Omega(t)$ by the conformal map
 $$
 P(w)=\Big(\tan\Big(\frac{w}{2}\Big)\Big)^{1/2},\quad w\in\C,
 $$
whose intention is to keep apart the self-intersecting point by taking the branch of the square root above passing through this crucial point. The equation in this new domain can be written as follows:

\begin{align}\label{zeq}
\tilde{z}_t(\al,t) & = Q^2(\al,t)BR(\tilde{z},\tilde{\omega})(\al,t) + \tilde{c}(\al,t)\tilde{z}_{\al}(\al,t),\end{align}
\begin{align}\label{eqomega}
\tilde{\omega}_{t}(\al,t)  =& -2 BR_t(\tilde{z},\tilde{\omega})(\al,t) \cdot \tilde{z}_{\al}(\al,t) - (Q^{2})_{\al}(\al,t)|BR(\tilde{z},\tilde{\omega})|^{2} (\al,t) - \Big(\frac{Q^2(\al,t)\tilde{\omega}(\al,t)^2}{4|\tilde{z}_{\al}(\al,t)|^{2}}\Big)_{\al}\nonumber \\
& + 2\tilde{c}(\al,t) BR_\al(\tilde{z},\tilde{\omega}) \cdot \tilde{z}_{\al}(\al,t) + \left(\tilde{c}(\al,t)\tilde{\omega}(\al,t)\right)_{\al}
- 2 \left(P^{-1}_2(\tilde{z}(\al,t))\right)_{\al} \nonumber  \end{align}
where
$$\tilde{z}(\al,t)=P(z(\al,t)),\quad Q^2(\al,t) = \left|\frac{dP}{dw}(P^{-1}(\tilde{z}(\al,t)))\right|^{2} \text{ and } \al \in \mathbb{T}.$$

(From now on we will omit the superscript tilde in the notation).

We start computing a numerical approximation of a solution  to the water waves equation \ref{zeq} that starts as a splash, turns over and finally is a graph. Such a candidate is depicted in Fig. \ref{splash}. With this aproximation we can construct explicit functions $(x,\gamma)$ that solve the system
\begin{equation}
%\label{CharlieFlat}
\left\{
\begin{array}{rl}
x_t  =& Q^2(x)BR(x,\gamma) +  b x_{\al} + f\\
 \gamma_t = &-2BR_{t}(x,\gamma) \cdot x_{\al}  - (Q^2(x))_{\al}|BR(x,\gamma)|^{2}-\left(\frac{Q^2(x)\gamma^2}{4|x_{\al}|^{2}}\right)_{\al}\\
 &+ 2bBR_{\al}(x,\gamma) \cdot x_{\al} + (b\gamma)_{\al} -2(P^{-1}_{2}(x))_{\al} + g
\end{array}
\right.
\end{equation}

where $f$ and $g$ are errors that we hope are small. By using the computer we are able to give rigorous bounds for these errors. The  question we want to answer is if there exists an exact solution $(z,\omega)$ of the water waves equation close to these functions $(x,\gamma)$. That means we need to prove the following theorem:

\begin{theorem}
\label{stabilitytheorem}
Let
$$ D(\al,t) \equiv z(\al,t) - x(\al,t), \quad d(\al,t) \equiv \omega(\al,t) - \gamma(\al,t), \quad \mathcal{D}(\al,t) \equiv \varphi(\al,t) - \psi(\al,t)$$
where $(x,\gamma,\psi)$ are the solutions of
\begin{equation}
\label{CharlieFlat}
\left\{
\begin{array}{rl}
x_t & = Q^2(x)BR(x,\gamma) +  b x_{\al} + f\\
b & = \underbrace{\frac{\al + \pi}{2\pi}\int_{-\pi}^{\pi}(Q^2
BR(x,\gamma))_{\al}\frac{x_\al}{|x_{\al}|^{2}} d\al - \int_{-\pi}^{\al}(Q^2 BR(x,\gamma))_{\beta}\frac{x_{\al}}{|x_{\al}|^{2}}d\beta}_{b_s} \\
& + \underbrace{\frac{\al + \pi}{2\pi}\int_{-\pi}^{\pi}f_{\al}\frac{x_\al}{|x_{\al}|^{2}} d\al - \int_{-\pi}^{\al}f_{\beta}\frac{x_{\beta}}{|x_{\beta}|^{2}}d\beta}_{b_e}\\
\gamma_t & + 2BR_{t}(x,\gamma) \cdot x_{\al} = - (Q^2(x))_{\al}|BR(x,\gamma)|^{2} + 2bBR_{\al}(x,\gamma) \cdot x_{\al} + (b\gamma)_{\al}  \\
& \qquad \qquad \qquad \quad \; \, -\left(\frac{Q^2(x)\gamma^2}{4|x_{\al}|^{2}}\right)_{\al} - 2(P^{-1}_{2}(x))_{\al} + g \\
\psi(\al,t) & = \frac{Q^2_x(\al,t)\gamma(\al,t)}{2|x_{\al}(\al,t)|} - b_s(\al,t)|x_{\al}(\al,t)|,
\end{array}
\right.
\end{equation}
where $(z,\omega)$ are the solutions of \eqref{CharlieFlat} with $f \equiv g \equiv 0$, $\varphi$ is the function
$$\varphi=\frac{Q^2_z(\alpha,t)\omega(\alpha,t)}{2|z_\alpha(\alpha,t)|}-b(\alpha,t)|z_\alpha(\alpha,t)|,$$

 and $E$ is the following norm for the difference
$$E(t) \equiv \left(\|D\|^{2}_{H^{3}} + \int_{-\pi}^{\pi}\frac{Q^2\sigma_{z}}{|z_{\al}|^{2}}|\partial^{4}_{\al}D|^{2} + \|d\|^{2}_{H^{2}} + \|\mathcal{D}\|^{2}_{H^{3+\frac{1}{2}}}\right).$$
Then we have that
$$\left|\frac{d}{dt}E(t)\right|\leq \mathcal{C}(t)(E(t)+E^{k}(t))+c\delta(t)$$
where $$\mathcal{C}(t)= \mathcal{C}(\mathcal{E}(t),\|x\|_{H^{5+\frac12}}(t),\|\gamma\|_{H^{3+\frac12}}(t),
\|\zeta\|_{H^{4+\frac12}}(t),\|F(x)\|_{L^\infty}(t))$$ and 
$$\delta(t)=(\|f\|_{H^{5+\frac12}}(t)+\|g\|_{H^{3+\frac12}}(t))^k+(\|f\|_{H^{5+\frac12}}(t)+\|g\|_{H^{3+\frac12}}(t))^2, \text{ $k$ big enough}$$
depends on the norms of $f$ and $g$, and $\mathcal{E}(t)$ is given by
\begin{align*}
\mathcal{E}(t)=&\|z\|^2_{H^3}(t)+\int_{\T}\frac{Q^2\sigma_z}{|z_\al|^2}|\da^4 z|^2d\al+\|F(z)\|^2_{L^\infty}(t)\\
&+\|\om\|^2_{H^{2}}(t)+
\|\varphi\|^2_{H^{3+\frac12}}(t)+\frac{|z_\al|^2}{m(Q^2\sigma_z)(t)}+\sum_{l=0}^4\frac{1}{m(q^l)(t)}
\end{align*}
where the $L^\infty$ norm of the function
$$
F(z)\equiv \frac{|\beta|}{|z(\al,t)-z(\al-\beta,t)|},\quad \al,\beta\in\T
$$
measures the arc-chord condition,
\begin{align}
\begin{split}\label{R-T}
 \sigma_{z}  \equiv& \left(BR_{t}(z,\omega) + \frac{\varphi}{|z_{\al}|}BR_{\al}(z,\omega)\right) \cdot z_{\al}^{\perp} + \frac{\omega}{2|z_{\al}|^{2}}\left(z_{\al t} + \frac{\varphi}{|z_{\al}|}z_{\al \al}\right) \cdot z_{\al}^{\perp} \\
& + Q\left|BR(z,\omega) + \frac{\omega}{2|z_{\al}|^{2}}z_{\al}\right|^{2}(\nabla Q)(z) \cdot z_{\al}^{\perp}
 + (\nabla P_{2}^{-1})(z) \cdot z_{\al}^{\perp}
\end{split}
\end{align}
is the Rayleigh-Taylor function,
$$
m(Q^2\sigma_z)(t)\equiv\min_{\al\in\T}Q^2(\al,t)\sigma_z(\al,t),
$$
and finally
$$
m(q^l)(t)\equiv\min_{\al\in\T}|z(\al,t)-q^l|
$$
for $l=0,...,4$, with
\begin{equation}\label{points}
q^0=\left(0,0\right),\quad
q^1=\left(\frac{1}{\sqrt{2}},\frac{1}{\sqrt{2}}\right),\quad
q^2=\left(\frac{-1}{\sqrt{2}},\frac{1}{\sqrt{2}}\right),\quad
q^3=\left(\frac{-1}{\sqrt{2}}, \frac{-1}{\sqrt{2}}\right),\quad
q^4=\left(\frac{1}{\sqrt{2}}, \frac{-1}{\sqrt{2}}\right),
\end{equation}
which are the singular points of the transformation $P$.
\end{theorem}

\begin{remark}
We can absorb the terms in $\mathcal{E}(t)$ by $E(t)$ raised to an appropriate power and terms in $(x,\gamma)$ by performing the splitting $\|z\| = \|z-x\| + \|x\|$ (or the analogous one for a different variable) for any norm or any quantity that appears in $\mathcal{E}(t)$.
\end{remark}
Theorem \ref{stabilitytheorem}  was announced in \cite{Castro-Cordoba-Fefferman-Gancedo-GomezSerrano:splash-water-waves}.

%We can construct a solution $(z,\omega,\varphi)$ that satisfies \eqref{CharlieFlat} with $f = g = 0$ and very similar same initial conditions
%\begin{align*}
%z_0(\al) \approx x_0(\al), \quad \omega_0(\al) \approx \gamma_0(\al), \quad \varphi_0(\al) \approx \psi_0(\al).
%\end{align*}

If we knew $\mathcal{C}(t), f(t), g(t), k$ or bounds on them, a priori, then we could provide bounds on $\mathcal{E}(t)$ at any time $T$. We point out here that $E(t)$ controls the norm $\|\da z^{1}(\al) - \da x^{1}(\al)\|_{L^{\infty}}$. Let $T_g$ be a time in which the approximate solution is a graph, i.e. $\da x^{1}(\al,T_g) > 0 \quad \forall \alpha$. Now, if $E(T_g) < \da x^{1}(\al,T_g)$ then
\begin{align*}
\da z^{1}(\al,T_g) > -\|\da z^{1}(\al) - \da x^{1}(\al)\|_{L^{\infty}} + \da x^{1}(\al,T_g) > 0,
\end{align*}
and this shows that $z$ is a graph. In other words, the possible set of solutions of the water waves equation is a ball centered at $(x,\gamma,\zeta)$ with the topology given by $E$. All of the elements of this ball are graphs, therefore the solution is necessarily a graph. Thus, the problem is reduced to study and find bounds for $\mathcal{C}(t), f(t), g(t), k$.

The recent developments of computer architecture have boosted their use in mathematics, giving birth to a full set of new results only achievable by this enormous power. However, it has the drawback that floating-point operations can not be performed exactly, resulting in numerical errors. In order to overcome this difficulty  and be able to prove rigorous results, we use the so-called \emph{interval arithmetics}, in which instead of working with arbitrary real numbers, we perform computations over intervals which have representable numbers as endpoints. On these objects, an arithmetic is defined in such a way that we are guaranteed that for every $x \in X, y \in Y$
\begin{align*}
x \star y \in X \star Y,
\end{align*}

for any operation $\star$. For example,
\begin{align*}
[\underline{x},\overline{x}] + [\underline{y},\overline{y}] & = [\underline{x} + \underline{y}, \overline{x} + \overline{y}] \\
[\underline{x},\overline{x}] \times [\underline{y},\overline{y}] & = [\min\{\underline{x}\underline{y},\underline{x}\overline{y},\overline{x}\underline{y},\overline{x}\overline{y}\},\max\{\underline{x}\underline{y},\underline{x}\overline{y},\overline{x}\underline{y},\overline{x}\overline{y}\}]
\end{align*}
We can also define the interval version of a function $f(X)$ as an interval $I$ that satisfies that for every $x \in X$ we have $f(x) \in I$.

The article is organized as follows:  in sections 2 and 3 we give some details about how to control the errors $f$, $g$ and the constants that arise in Theorem \ref{stabilitytheorem} by using the computer. Finally, in section 4 we give a complete proof of Theorem \ref{stabilitytheorem}.

\section{Bounds for $f(t)$ and $g(t)$}

\subsection{Representation of the functions and Interpolation}

The first thing one has to decide is how to represent the data and how to pass from the cloud of points in space-time obtained by non-rigorous simulation to a function defined everywhere in $[-\pi,\pi] \times [0,T]$. We need to interpolate in some way. %One of the first things that can come to one's mind is to use the first $N$ Fourier modes. This approach has two disadvantages. The first one is that the linearized water waves equation is not dissipative, hence we will not have a control for the tails uniformly for all time, even in the case where the tails have very small norms at time zero. The second disadvantage is numerical. Suppose $N \sim 10^{3}$, which is a reasonable guess. Since we need to take 5.5 derivatives in the curve, the high order coefficients will be multiplied by roughly a factor of $10^{15}$. If we work with a 64-bit representation, machine epsilon is of the order $10^{-16}$ and we will run into trouble because the computer will not distinguish between zero and non-zero values. Of course, this problem can be solved if we use high precision arithmetics, but in any case we would be multiplying and dividing by very big constants.

In our case, we chose to represent the functions $x$ and $\gamma$ by piecewise polynomials (splines) of high degree (10) in space, and low degree (3) in time. To do so, we first interpolate in space for every node in the time mesh. The interpolation is made via B-Splines. Since the interpolation is reduced to solve a linear (interval) system $Ac = y$, where $A$ is constant in time and space and $y$ depends on the values of the function at time $t$ since the mesh in space is constant, we precondition by multiplying by the non-rigorous inverse of the midpoints of the entries of $A$. We remark that the system is interval-based because we need to produce a curve that is a splash (i.e. there have to be two points $\al_1, \al_2$ such that we can guarantee $x_0(\al_1) = x_0(\al_2)$. Finally, the system is solved using a rigorous Gauss-Seidel iterative method. We also remark that the need for interval-based calculations is only strictly necessary at time $t = 0$ since it is the only point in which we have to guarantee some equality. By working with multiprecision (1024 bits) we can get widths in the coefficients of the order of $10^{-300}$. In order to perform interpolation in time, we fix the values of the function and its time derivative at the mesh points. This gives us lots of systems of 4 equations (the values of the function and its derivative at both endpoints) and 4 unknowns (the 4 coefficients of the degree 3 polynomial) but with an explicit formula for each of them. With this method, our spline will be $C^{1}$ in time but it might not be $C^{2}$.

\subsection{Rigorous bounds for Singular integrals}

In this section we will discuss the computational details of the rigorous calculation of some singular integrals. In particular we will focus on the Hilbert transform, but the methods apply to any integral kernel whose main singularity is homogeneous of degree -1. Parts of the computation  (the $N$ part) are slightly related to the Taylor models with relative remainder presented in M. Jolde\c{s}' thesis \cite{Joldes:rigorous-polynomial-approximations}.

Let us suppose that we have a function $f$ given explicitly by a spline (piecewise polynomial) which is $C^{k-1}$ everywhere and $C^{k}$ except at finitely many  points (the points in which the different pieces of the spline are glued together). We need to calculate rigorously the Hilbert Transform of $f$, that is

\begin{align*}
Hf(x) = \frac{PV}{\pi} \int_{\T} \frac{f(x)-f(y)}{2\tan\left(\frac{x-y}{2}\right)}dy,
\end{align*}

and we want to approximate it by a piecewise polynomial function with less regularity, plus an error that can be bounded in $H^{q}, 0\leq q \leq c < k$ and in $L^{\infty}$. Let us assume that the knots of the spline are $\al_i$, $i = 0, \ldots, N-1$ and that we fix $x \in [\al_i, \al_{i+1}]$ where the indices are taken modulo $N$ and the distance between the indices is taken over $\mathbb{Z}_{N}$. We can split our integral in

\begin{align*}
Hf(x) & = \frac{PV}{\pi} \int_{\T} \frac{f(x)-f(y)}{2\tan\left(\frac{x-y}{2}\right)}dy
= \frac{PV}{\pi} \sum_{j} \int_{\al_{j}}^{\al_{j+1}} \frac{f(x)-f(y)}{2\tan\left(\frac{x-y}{2}\right)}dy \\
& = \frac{PV}{\pi} \sum_{|j-i|>K} \int_{\al_{j}}^{\al_{j+1}} \frac{f(x)-f(y)}{2\tan\left(\frac{x-y}{2}\right)}dy
+ \frac{PV}{\pi} \sum_{|j-i|\leq K} \int_{\al_{j}}^{\al_{j+1}} \frac{f(x)-f(y)}{2\tan\left(\frac{x-y}{2}\right)}dy \\
& \equiv Hf^{F}(x) + Hf^{N}(x).
\end{align*}

Now, if we want to express $Hf^{F}(x)$ as a polynomial, it is easy since the integrand does not have a singularity. Hence
\begin{align*}
Hf^{F}(x) & = \frac{PV}{\pi} \sum_{|j-i|>K} \int_{\al_{j}}^{\al_{j+1}} \frac{f(x)-f(y)}{2\tan\left(\frac{x-y}{2}\right)}dy
= \frac{PV}{\pi} \sum_{|j-i|>K} \int_{\al_{j}}^{\al_{j+1}} F^{j}(x,y) dy \\
& = \sum_{|j-i|>K} \int_{\al_{j}}^{\al_{j+1}} \sum_{n,m} c_{nm} (x-x^{*}(i))^{m}(y-y^{*}(j))^{n} + E(x,y) dy \equiv P(x) + E(x),
\end{align*}

where $E$ accounts for the error and is a polynomial with interval coefficients. Typically, we will use as the points for the Taylor expansions $x^{*}(i) = \al_i$ since we will compare the resulting polynomial with another one of the form $\sum_{j} b_j (x-x_{i})^{j}$ and we will also choose $y^{*}(j) = \frac{\al_{j} + \al_{j+1}}{2}$. This choice is useful for two reasons: first, we will only have to integrate half of the terms since the rest will integrate to zero; and second, the error estimates will be better for this choice of $y^{*}(j)$ in the sense that the coefficients will be smaller. All the computations will be carried out using automatic differentiation. We should remark that we can get estimates for the error $E$ in any of the above mentioned norms without having to recompute it since the relation
\begin{align*}
\partial_{x}^{q} Hf^{F}(x) - \partial_{x}^{q}P(x) = \partial_{x}^{q} E(x)
\end{align*}
holds for every $q < k$.

Now, we move on to the term $Hf^{N}(x)$. In this case, we perform a Taylor expansion in both the denominator
\begin{align*}
2\tan\left(\frac{x-y}{2}\right) = (x-y)+c(x-y)^{3}, \quad c = \text{ small (interval) constant}
\end{align*}
and the numerator
%\begin{align*}
%f(x)-f(y) & =
%(x-\al_i)f'(\al_i) - (y-\al_i)f'(\al_i) +
%\frac{1}{2}(x-\al_i)^{2}f'(\al_i) - \frac12(y-\al_i)^{2}f'(\al_i) \\
%& + \ldots +
%\frac{1}{n!}(x-\al_i)^{n}f^{n}(\eta_1) - \frac{1}{n!}(y-\al_i)^{n}f^{n}(\eta_2)
%\end{align*}
\begin{align*}
f(x) = f(y) + (x-y)f'(y) + \frac{1}{2}(x-y)^2f''(y) + \ldots \frac{1}{n!}(x-y)^{k-1}f^{k-1}(\eta),
\end{align*}
where $\eta$ belongs to an intermediate point between $x$ and $y$, which we can enclose in the convex hull of $[\al_i, \al_{i+1}]$ and $[\al_{j}, \al_{j+1}]$ where the convex hull is understood in the torus. Since typically $K$ will be very small (compared to $N$) there is no ambiguity in the definition. Finally, we can factor out $(x-y)$ and divide both in the numerator and the denominator. Since we know $f(y)$ explicitly, we can perform the explicit integration and get a piecewise polynomial as a result.

\subsection{Estimates of the norm of the Operator $I + T$}

In this subsection we will outline how to compute the norm of the operator $I + T = I + 2 \langle BR(z,\cdot), z_{\al} \rangle$. Since the operator $T$ behaves like a Hilbert Transform plus smoothing terms, we will describe how to calculate rigorously with the help of a computer an estimate for the norm of its inverse. The procedure is more general and can be applied to a bigger family of kernels.
Let $\mathbb{T} = \mathbb{R}/2\pi \mathbb{Z}$, and let $A(x), B(x)$ be real-valued functions on $\mathbb{T}$. Also, let $E(x,y)$ be a real-valued function on $\mathbb{T} \times \mathbb{T}$. We assume $A, B$ and $E$ are given by explicit formulas such as as perhaps piecewise trigonometric polynomials or splines, and $E(x,y)$ is a trigonometric polynomial on each rectangle $I \times J$ of some partition of $\mathbb{T} \times \mathbb{T}$. We suppose $A,B,E$ are smooth enough.

Let $H$ be the Hilbert transform acting on functions on $\mathbb{T}$, i.e.
\begin{align*}
Hf(x) = \frac{PV}{2\pi}\int_{\mathbb{T}} \cot\left(\frac{y}{2}\right)f(x-y)dy.
\end{align*}

Assume that $A$ and $B$ have no common zeros on $\mathbb{T}$.

Let
\begin{align*}
Sf(x) = A(x)f(x) + B(x)Hf(x) + \int_{\T}E(x,y)f(y)dy, \quad f \in L^{2}(\T).
\end{align*}
Thus, $S$ is a singular integral operator.

We hope that $S^{-1}$ exists and has a not-so-big norm on $L^{2}$, but we don't know this yet.

Our goal here is to find approximate solutions $F$ of the equation $SF = f$ for suitable given $f \in L^{2}(\T)$, and to check that $\|SF - f\|_{L^{2}(\T)}  < \delta$ for suitable $\delta$. Our computation of $F$ will be based on heuristic ideas, but the computation of an upper bound for $\|SF - f\|_{L^{2}(\T)}$ will be rigorous. In our case, $A(x) = 1, B(x) = 1$.

To carry this out, let $H_0 \subset H_1 \subset L^{2}(\T)$ be finite-dimensional subspaces, e.g. with $H_i$ consisting of the span of wavelets (from a wavelet bases) having lengthscale $\geq 2^{-N_i}$. Here $N_1 \geq N_0 + 3$ (say). Let $\pi_{i}$ be the orthogonal projection from $L^{2}(\T)$ to $H_i$, and let us solve the equation

\begin{align}
\label{charlieinversionstar}
\pi_1 S \pi_1 F = \pi_0 f.
\end{align}

If $f$ is given explicitly in a wavelet bases, then \eqref{charlieinversionstar} is a linear algebra problem, since $\pi_1 S \pi_1$ is of finite rank, and its matrix (in terms of some given basis for $H_1$) can be computed explicitly.

\begin{itemize}
\item If $\pi_0 f \not \in \text{Range}(\pi_1 S \pi_1)$, then our heuristic procedure fails.
\item If $\pi_0 f \in \text{Range}(\pi_1 S \pi_1)$, then we find $F \in H_1$ such that $\pi_1 S \pi_1 F = \pi_0 f$, i.e. $\pi_1 SF = \pi_0 f$.
\end{itemize}

We then have
\begin{align*}
\|SF - f\|_{L^{2}(\T)} \leq \| (I - \pi_1) SF\|_{L^{2}(\T)} + \|(I - \pi_0) f\|_{L^{2}(\T)},
\end{align*}
and both norms on the right-hand side may be estimated explicitly.

Now, our goal is to make a heuristic computation of an operator of the form
\begin{align*}
\tilde{S}f(x) = \tilde{A}(x)f(f) + \tilde{B}(x)Hf(x) + \int_{\T}\tilde{E}(x,y)f(y)dy
\end{align*}
such that $S\tilde{S} - I$ has small norm on $L^{2}(\T)$.

Here, we will make a heuristic computation of $\tilde{S}$; later we will give a rigorous upper bound for the norm of $S\tilde{S} - I$ on $L^{2}(\T)$. By a heuristic computation of $\tilde{S}$ we mean a heuristic computation of $\tilde{A}, \tilde{B}$ and $\tilde{E}$.

We first find $\tilde{A}$ and $\tilde{B}$ by setting
\begin{align*}
(A + iB)(\tilde{A} + i \tilde{B}) = 1 \Rightarrow
\left\{
\begin{array}{rcl}
A \tilde{A} - B\tilde{B} & = 1& \\
A \tilde{B} + B\tilde{A} & = 0& \\
\end{array}
\right.
\end{align*}

Then, this means that
\begin{align*}
S\tilde{S} = (A \tilde{A} - B\tilde{B}) + (A \tilde{B} + B\tilde{A})H + \text{ Smoothing terms}
= I + \text{ Smoothing terms}
\end{align*}

So, from now on, we suppose that $\tilde{A}$ and $\tilde{B}$ are known. For the operator $I + T$, this means $\tilde{A} = 1/2, \tilde{B} = -1/2$. We want to compute $\tilde{E}$. Now, let $\{\phi_{\nu}\}$ be some orthonormal basis for $L^{2}(\T)$, for example a wavelet basis. By the previous methods, we can try to find functions $\psi_{\nu} \in L^{2}(\T)$ such that $S\psi_{\nu} - \phi_{\nu}$ has small norm. We carry this for $\nu = 1, \ldots, N$ for a large $N$. We now try to make $\tilde{E}$ satisfy

\begin{align}
\tilde{A}(x)\phi_{\nu}(x) + \tilde{B}(x)H\phi_{\nu}(x) + \int_{\T}\tilde{E}(x,y)\phi_{\nu}(y) = \psi_{\nu}(x) \text{ for } \nu = 1, \ldots, N.
\end{align}

Thus, we want
\begin{align}
\label{charlieadmiration}
\int_{\T}\tilde{E}(x,y) \phi_{\nu}(y) dy = \left(\psi_{\nu}(x) - \tilde{A}(x) \phi_{\nu}(x) - \tilde{B}(x)H\phi_{\nu}(x)\right) \equiv \psi_{\nu}^{\#}(x), \quad \nu = 1, \ldots, N.
\end{align}

Note that $\psi_{\nu}^{\#}$ can be computed explicitly.

Since the $\phi_{\nu}$ (all $\nu$) form an orthonormal basis for $L^{2}(\T)$, it is natural to define
\begin{align*}
\tilde{E}(x,y) = \sum_{\nu=1}^{N} \psi_{\nu}^{\#} \phi_{\nu}(y).
\end{align*}
This can be computed explicitly, and it satisfies \eqref{charlieadmiration}. Thus, we can compute
\begin{align}
\label{charlie2admiration}
S\tilde{S} & = (A+BH+E)(\tilde{A} + \tilde{B}H + \tilde{E})\nonumber \\
& = A\tilde{A} + A \tilde{B}H + A \tilde{E} + BH \tilde{A} + BH\tilde{B}H + BH \tilde{E} + E\tilde{A} + E\tilde{B}H + E\tilde{E} \nonumber \\
& = A\tilde{A} + A \tilde{B}H + A \tilde{E} + B\tilde{A}H + B[H,\tilde{A}] - B\tilde{B} + B[H,\tilde{B}]H \nonumber\\
& + BH \tilde{E} + E\tilde{A} + E\tilde{B}H + E\tilde{E} \nonumber \\
& = (A\tilde{A} - B\tilde{B}) + (A \tilde{B} + B\tilde{A})H + \{A \tilde{E} + B[H,\tilde{A}] + B[H,\tilde{B}]H \nonumber\\
& + BH \tilde{E} + E\tilde{A} + E\tilde{B}H + E\tilde{E}\}
\end{align}

We claim that all terms enclosed in curly brackets are integral operators of the form
\begin{align*}
S^{\#}f(x) = \int_{\T}E^{\#}(x,y) f(y) dy,
\end{align*}
for an $E^{\#}$ that we can calculate. Let us go term by term
\begin{itemize}
\item $A \tilde{E}$ has the form $S^{\#}$, with $E^{\#}(x,y) = A(x)\tilde{E}(x,y)$.
\item $B[H,\tilde{A}]$ has the form $S^{\#}$, with $E^{\#}(x,y) = \frac{1}{2\pi}B(x)\cot\left(\frac{x-y}{2}\right)(\tilde{A}(x)-\tilde{A}(y))$.

Note that if $\tilde{A}$ is a piecewise trigonometric polynomial and $C^{k}$, then $E^{\#}$ can easily be computed modulo a small error in $C^{k-1}$.
\item $B[H,\tilde{B}]H$ has the form $S^{\#}$, with
\begin{align*}
E^{\#}(x,y) & = \frac{1}{4\pi^{2}}B(x)PV\int\cot\left(\frac{x-z}{2}\right)(\tilde{B}(x)-\tilde{B}(z))\cot\left(\frac{z-y}{2}\right)dz. \\
& = \frac{1}{4\pi^{2}}B(x)PV\int\left\{\cot\left(\frac{x-z}{2}\right)(\tilde{B}(x)-\tilde{B}(z))-2\tilde{B}'(x)\right\}\cot\left(\frac{z-y}{2}\right)dz.
\end{align*}
\item $BH\tilde{E}$ has the form $S^{\#}$, with
\begin{align*}
E^{\#}(x,y) & = \frac{1}{2\pi}B(x)PV\int\cot\left(\frac{x-z}{2}\right)\tilde{E}(z,y)dz. \\
& = \frac{1}{2\pi}B(x)PV\int\cot\left(\frac{x-z}{2}\right)\left(\tilde{E}(z,y)-\tilde{E}(x,y)\right)dz.
\end{align*}
\item $E\tilde{A}$ has the form $S^{\#}$, with $E^{\#}(x,y) = \tilde{E}(x,y)\tilde{A}(y)$.
\item $E\tilde{B}H$ has the form $S^{\#}$, with
\begin{align*}
E^{\#}(x,y) & = \frac{1}{2\pi}PV\int E(x,z)\tilde{B}(z)\cot\left(\frac{z-y}{2}\right)dz. \\
& = \frac{1}{2\pi}PV\int \left\{E(x,z)\tilde{B}(z)-E(x,y)\tilde{B}(y)\right\}\cot\left(\frac{z-y}{2}\right)dz.
\end{align*}
\item $E\tilde{E}$ has the form $S^{\#}$, with $E^{\#}(x,y) = \int E(x,z)\tilde{E}(z,y)dz$.
\end{itemize}

This proves the claim.

Letting $\mathcal{E}^{\#}f(x) = \int_{\T}E^{\#}(x,y)f(y)dy$ be the operator in curly brackets in \eqref{charlie2admiration}, we see that
\begin{align*}
S\tilde{S} = (A\tilde{A} - B\tilde{B}) + (A\tilde{B} + B\tilde{A})H + \mathcal{E}^{\#},
\end{align*}
and that the function $E^{\#}(x,y)$ can be computed modulo a small error in $C^{0}(\T \times \T)$. Therefore, we obtain an upper bound for the norm of $S\tilde{S} - I$, namely
\begin{align*}
\max| A\tilde{A} - B\tilde{B} - 1| + \max |A\tilde{B} + B\tilde{A}|
+ \max\left\{\max_{x}\int|E^{\#}(x,y)|dy,\max_{y}\int|E^{\#}(x,y)|dx\right\}.
\end{align*}

Defining $S_{err} := S\tilde{S} - I$, we obtain an explicit upper bound $\delta$ for the norm of $S_{err}$ on $L^{2}(\T)$. We hope that $\delta < 1$. If not, then we fail.

Suppose $\delta < 1$. Then
\begin{align*}
S\tilde{S} = I +  S_{err} \Rightarrow S\tilde{S}(I+S_{err})^{-1} = I,
\end{align*}
so we obtain a right inverse for $S$, namely $\tilde{S}(I+S_{err})^{-1}$, which has norm at most
\begin{align}
\label{charlie2star}
\|\tilde{S}\|(1-\delta)^{-1},
\end{align}
where $\|\tilde{S}\|$ denotes the norm of $\tilde{S}$ as an operator on $L^{2}(\T)$. Recall
\begin{align*}
\tilde{S}f(x) = \tilde{A}(x)f(x) + \tilde{B}(x)Hf(x) + \int_{\T}\tilde{E}(x,y)f(y)dy.
\end{align*}
Therefore,
\begin{align*}
\|\tilde{S}\| \leq \max|\tilde{A}(x)| + \max|\tilde{B}(x)| + \max\left\{\max_{x}\int|\tilde{E}(x,y)|dy,\max_{y}\int|\tilde{E}(x,y)|dx\right\}.
\end{align*}

Plugging that bound into \eqref{charlie2star}, we obtain an explicit upper bound for the norm on $L^{2}$ of a right inverse for $S$. Similarly (by looking at $\tilde{S}S$ instead of $S\tilde{S}$), we obtain an upper bound for the norm on $L^{2}$ of a left inverse for $S$.

\begin{remark}
To estimate e.g. $\max_{x}\int_{\T}|E^{\#}(x,y)|dy$ it may be enough just to use the trivial estimate
\begin{align*}
\max_{x}\int_{\T}|E^{\#}(x,y)|dy \leq 2\pi \max_{x,y}|E^{\#}(x,y)|
\end{align*}
\end{remark}
\begin{remark}[Time dependent solutions]
For $t \in [t_0,t_1]$ (a small time interval), let
\begin{align*}
S_t f(x) = A(x,t)f(x) + B(x,t)Hf(x) + \int_{\T}E(x,y,t)f(y)dy,
\end{align*}
where (for each $t$),$A(\cdot,t),B(\cdot,t),E(\cdot,\cdot,t)$ are as assumed above.

If $A,B,E$ depend in a reasonable way on $t$, then one shows easilly that
\begin{align*}
\| S_t - S_{t_0}\| < \eta \text{ for all } t \in [t_0,t_1].
\end{align*}

We can make $\eta$ small by taking $t_1$ close enough to $t_0$. Suppose we prove that $\|S_{t_0}^{-1}\| \leq C_0$ by the previous methods. Then, of course we obtain an upper bound for $\| S_{t}^{-1}\|$ valid for all $t \in [t_0,t_1]$.
\end{remark}

\section{Bounds for $\mathcal{C}(t)$ and $k$}

\subsection{Writing the differential inequality as a differential system of equations}

The calculation of a bound for $\mathcal{C}(t)$ requires more effort than the previous one since one needs to calculate the terms one by one and add all their contributions to $\mathcal{C}(t)$. For example, in order to calculate the evolution of the norm $\|D\|_{H^{k}}(t)$ a systematic approach is to take $k$ derivatives ($k$ ranging from 0 to 4) in the equation for the evolution of $z$ (\ref{CharlieFlat} with $f = g = 0$), take another $k$ derivatives in the equation for $x$ (\ref{CharlieFlat} with arbitrary $f,g$) and subtract them. Let us focus from now on in the term $Q(z)^{2}BR(z,\omega) - Q(x)^2BR(x,\gamma)$ and its derivatives. One notices that in order to write a term in the variables $(z,\omega,\varphi)$ composed of $a$ factors minus its counterpart in the variables $(x,\gamma,\psi)$ in a suitable way (i.e. as a sum of terms that only have factors $x,\gamma,\psi,D,d,\mathcal{D}$) then the number of terms is $2^{a}-1$. The way of writing it is the classical way of adding and subtracting the same term with the purpose of creating differences of terms and eliminate all the occurrences of the variables $(z,\omega,\varphi)$. An example for the Birkhoff-Rott operator (with $Q = 1$) is given next. We should remark that the computation and bounding of the Birkhoff-Rott is the most expensive one, the rest of the terms being easier.
\begin{align*}
&BR(z,\om) - BR(x,\gamma)   = \frac{1}{2\pi}\int \frac{(x(\alpha)-x(\beta))^\perp}{|x(\alpha)-x(\beta)|^2}\left(\om(\beta)-\gamma(\beta)\right)d\beta \\
& + \frac{1}{2\pi}\int \frac{(z(\alpha)-z(\beta))^\perp-(x(\alpha)-x(\beta))^\perp}{|x(\alpha)-x(\beta)|^2}\left(\gamma(\beta)\right)d\beta \\
& + \frac{1}{2\pi}\int \frac{(z(\alpha)-z(\beta))^\perp-(x(\alpha)-x(\beta))^\perp}{|x(\alpha)-x(\beta)|^2}\left(\om(\beta)-\gamma(\beta)\right)d\beta \\
& + \frac{1}{2\pi}\int \left(\frac{1}{|z(\alpha)-z(\beta)|^2} - \frac{1}{|x(\alpha)-x(\beta)|^2}\right) (x(\alpha)-x(\beta))^\perp\gamma(\beta)d\beta \\
& + \frac{1}{2\pi}\int \left(\frac{1}{|z(\alpha)-z(\beta)|^2} - \frac{1}{|x(\alpha)-x(\beta)|^2}\right) (x(\alpha)-x(\beta))^\perp(\om(\beta) - \gamma(\beta))d\beta \\
& + \frac{1}{2\pi}\int \left(\frac{1}{|z(\alpha)-z(\beta)|^2} - \frac{1}{|x(\alpha)-x(\beta)|^2}\right) (z(\alpha)-z(\beta) - (x(\alpha)-x(\beta)))^\perp\gamma(\beta)d\beta \\
& + \frac{1}{2\pi}\int \left(\frac{1}{|z(\alpha)-z(\beta)|^2} - \frac{1}{|x(\alpha)-x(\beta)|^2}\right) (z(\alpha)-z(\beta) - (x(\alpha)-x(\beta)))^\perp(\om(\beta) - \gamma(\beta))d\beta
\end{align*}
After having seen this, it is clear that a tool that can perform symbolic calculations (derivation and basic arithmetic at least) and the correct grouping of the factors is required since the performance at this task by a human is not satisfactory. We developed a tool in 900 lines of C++ code that could do all this and output the collection of terms in Tex. We show an excerpt of the terms concerning the fourth derivative of $BR(z,\omega) - BR(x,\gamma)$. The total number of terms in that case is 2841.

\begin{align*}
& 2\pi \left(\da^{4} BR(x,\gamma) - \da^{4} BR(z,\omega) \right) = \\
&
+\int (\da^{4} x(\al)-\da^{4} x(\al-\beta))^{\perp}d(\al-\beta)\frac{1}{\left|x(\al)-x(\al-\beta)\right|^{2}}d\al\\
&
+\int (\da^{4} x(\al)-\da^{4} x(\al-\beta))^{\perp}d(\al-\beta) \left(\frac{1}{\left|x(\al)-x(\al-\beta)\right|^{2}}-\frac{1}{\left|z(\al)-z(\al-\beta)\right|^{2}}\right)d\al\\
&
+\int (\da^{4} x(\al)-\da^{4} x(\al-\beta))^{\perp}\gamma(\al-\beta)\left(\frac{1}{\left|x(\al)-x(\al-\beta)\right|^{2}} -\frac{1}{\left|z(\al)-z(\al-\beta)\right|^{2}}\right)d\al\\
&
+4\int (\da^{3} x(\al)-\da^{3} x(\al-\beta))^{\perp}\da d(\al-\beta)\frac{1}{\left|x(\al)-x(\al-\beta)\right|^{2}}d\al\\
&
+4\int (\da^{3} x(\al)-\da^{3} x(\al-\beta))^{\perp}\da d(\al-\beta)\left(\frac{1}{\left|x(\al)-x(\al-\beta)\right|^{2}} -\frac{1}{\left|z(\al)-z(\al-\beta)\right|^{2}}\right)d\al\\
&
-8\int (\da^{3} x(\al)-\da^{3} x(\al-\beta))^{\perp}d(\al-\beta)\left(\frac{1}{\left|x(\al)-x(\al-\beta)\right|^{2}}\right)^{2} \\
& \times (\da x(\al)-\da x(\al-\beta))\cdot(D(\al)-D(\al-\beta))d\al\\
&
-8\int (\da^{3} x(\al)-\da^{3} x(\al-\beta))^{\perp}d(\al-\beta)\left(\frac{1}{\left|x(\al)-x(\al-\beta)\right|^{2}}\right)^{2} \\
& \times (\da x(\al)-\da x(\al-\beta))\cdot(x(\al)-x(\al-\beta))d\al\\
&
-8\int (\da^{3} x(\al)-\da^{3} x(\al-\beta))^{\perp}d(\al-\beta)\left(\frac{1}{\left|x(\al)-x(\al-\beta)\right|^{2}}\right)^{2} \\
& \times (\da D(\al)-\da D(\al-\beta))\cdot(D(\al)-D(\al-\beta))d\al\\
&
-8\int (\da^{3} x(\al)-\da^{3} x(\al-\beta))^{\perp}d(\al-\beta)\left(\frac{1}{\left|x(\al)-x(\al-\beta)\right|^{2}}\right)^{2} \\
& \times (x(\al)-x(\al-\beta))\cdot(\da D(\al)-\da D(\al-\beta))d\al\\
&
-8\int (\da^{3} x(\al)-\da^{3} x(\al-\beta))^{\perp}d(\al-\beta)\left(\left(\frac{1}{\left|x(\al)-x(\al-\beta)\right|^{2}}\right)^{2}
- \left(\frac{1}{\left|z(\al)-z(\al-\beta)\right|^{2}}\right)^{2}\right) \\
& \times (\da D(\al)-\da D(\al-\beta))\cdot(D(\al)-D(\al-\beta)) d\al \\
& + 2831 \text{ more terms...}
\end{align*}

However, there is a significant way to reduce the number of terms in the estimates: writing the equation in complex form instead of vector form. Thus, we can write the evolution for $z$ in the following way:

\begin{align*}
\partial_t z^{*}(\al,t) = \frac{1}{2\pi} \int_{\T} \frac{1}{z(\al,t) - z(\beta,t)}\omega(\beta,t)d\beta + c(\al,t)\da z^{*}(\al,t)
\end{align*}

In this formulation,  the fourth derivative accounts for only 140 terms. We present the first 10 below.

\begin{align*}
& 2\pi \left(\da^{4} BR(x,\gamma) - \da^{4} BR(z,\omega) \right) \\
& =
-72\int (\da^{2} x(\al)-\da^{2} x(\al-\beta))(\da x(\al)-\da x(\al-\beta)) \left(\frac{1}{x(\al)-x(\al-\beta)}\right)^{4} \\
& \times (\da D(\al)-\da D(\al-\beta))d(\al-\beta)d\al\\
&
-72\int (\da^{2} x(\al)-\da^{2} x(\al-\beta))(\da x(\al)-\da x(\al-\beta)) \left(\frac{1}{x(\al)-x(\al-\beta)}\right)^{4} \\
& \times (\da D(\al)-\da D(\al-\beta))\gamma(\al-\beta)d\al\\
&
-72\int (\da x(\al)-\da x(\al-\beta))(\da^{2} D(\al)-\da^{2} D(\al-\beta))  \left(\frac{1}{x(\al)-x(\al-\beta)}\right)^{4} \\
& \times (\da D(\al)-\da D(\al-\beta))d(\al-\beta)d\al\\
&
-72\int (\da x(\al)-\da x(\al-\beta))(\da^{2} D(\al)-\da^{2} D(\al-\beta)) \left(\frac{1}{x(\al)-x(\al-\beta)}\right)^{4} \\
& \times (\da D(\al)-\da D(\al-\beta))\gamma(\al-\beta)d\al\\
&
-36\int (\da^{2} D(\al)-\da^{2} D(\al-\beta))\left(\da D(\al)-\da D(\al-\beta)\right)^{2} \left(\frac{1}{x(\al)-x(\al-\beta)}\right)^{4}d(\al-\beta)d\al\\
&
-36\int (\da^{2} D(\al)-\da^{2} D(\al-\beta))\left(\da D(\al)-\da D(\al-\beta)\right)^{2} \left(\frac{1}{x(\al)-x(\al-\beta)}\right)^{4}\gamma(\al-\beta)d\al\\
&
+8\int \left(\frac{1}{x(\al)-x(\al-\beta)}\right)^{3}(\da^{3} D(\al)-\da^{3} D(\al-\beta))(\da D(\al)-\da D(\al-\beta))d(\al-\beta)d\al\\
&
+8\int \left(\frac{1}{x(\al)-x(\al-\beta)}\right)^{3}(\da^{3} D(\al)-\da^{3} D(\al-\beta))(\da D(\al)-\da D(\al-\beta))\gamma(\al-\beta)d\al\\
&
+24\int \left(\frac{1}{x(\al)-x(\al-\beta)}\right)^{3}(\da^{2} D(\al)-\da^{2} D(\al-\beta))(\da D(\al)-\da D(\al-\beta))\da d(\al-\beta)d\al\\
&
+24\int \left(\frac{1}{x(\al)-x(\al-\beta)}\right)^{3}(\da^{2} D(\al)-\da^{2} D(\al-\beta))(\da D(\al)-\da D(\al-\beta))\da \gamma(\al-\beta)d\al\\
& + 130 \text{ more terms...}
\end{align*}

%\subsection{Complex notation}
%\subsection{Automated calculus of the terms}

The final observation is that if we consider $\mathcal{E}(t)$ as a scalar, we might not get suitable estimates. In order to get better estimates, we will modify the energy into a ``vectorized'' version $\mathcal{E}_{v}(t)$, which we will also denote by $\mathcal{E}(t)$ by abuse of notation. This new vectorized energy will be as follows

\begin{align*}
\mathcal{E}(t) =
\left(
\begin{array}{c}
\|D\|_{L^{2}} \\
\|D\|_{\dot{H^{1}}} \\
\|D\|_{\dot{H^{2}}} \\
\|D\|_{\dot{H^{3}}} \\
\|d\|_{L^{2}} \\
\|d\|_{\dot{H^{1}}} \\
\|d\|_{\dot{H^{2}}} \\
\vdots
\end{array}
\right),
\end{align*}
where the homogeneous spaces $\dot{H^{k}}$ have their norm defined by $\|f\|_{\dot{H^{k}}} = \|\da^{k} f\|_{L^{2}}$. With this vectorized system, we avoid both the bounding of any given norm by the full energy and any constant factor arising from interpolation between two Sobolev spaces. Thus, our constant $C(t)$ will roughly be of a size comparable to the largest eigenvalue of the linearized system.

\subsection{Estimates for the linear terms with $Q = 1$}

Since we expect $\mathcal{E}(t)$ to be small, the terms that affect more to the evolution of $\mathcal{E}(t)$ are the linear ones. We now report on the non-rigorous experiments over the linear terms to obtain an approximate bound of the behavior of the full system (i.e. an approximation to the largest eigenvalue of the linearized system). We remark that a multiplication of the estimates by a constant, even a small factor 2 for example, has a big impact on the system, rendering the estimates useless and the estimations not tight enough, because the type of estimates we are going to get are exponential in the product of the time elapsed between the splash and the graph and the constant. Therefore, we should be very careful and fine estimates have to be developed.

First of all, we will work with $Q = 1$ and later move on to the case $Q \neq 1$. We will adopt the following convention to denote the different Kernels (integral operators) that appear:

\begin{align*}
\Theta^{a_1, a_2, a_3, a_4}_{b_1, b_2}(\al,\beta) & = \frac{1}{(x(\al)-x(\beta))^{b_1}}(\da x(\al)-\da x(\beta))^{a_1}(\da^{2}x(\al)-\da^{2}x(\beta))^{a_2} \\
& \times (\da^{3} x(\al)-\da^{3} x(\beta))^{a_3}(\da^{4} x(\al)-\da^{4} x(\beta))^{a_4}\da^{b_2}\gamma(\beta) \\
\Theta^{a_1, a_2, a_3, a_4}_{b_1, -1}(\al,\beta) & = \frac{1}{(x(\al)-x(\beta))^{b_1}}(\da x(\al)-\da x(\beta))^{a_1}(\da^{2}x(\al)-\da^{2}x(\beta))^{a_2} \\
& \times (\da^{3} x(\al)-\da^{3} x(\beta))^{a_3}(\da^{4} x(\al)-\da^{4} x(\beta))^{a_4}. \\
\end{align*}

The operators for which $b_2 \neq -1$ will act on $D$ or its derivatives whereas the operators for which $b_2 = -1$ will act on $d$ or its derivatives. We now describe how to split the Kernels in such a way that they can be computed. For the case where $b_2 \neq -1$ we illustrate this by splitting $\Theta^{0,0,0,0}_{2,0}$, but the technique can be applied to any Kernel.

\begin{align}
\frac{1}{2\pi}\int \Theta^{0,0,0,0}_{2,0}&(D(\al) - D(\beta))d\beta  = \underbrace{\frac{1}{2\pi}D(\al)\int K(\al,\beta)\gamma(\beta)d\beta}_{T_1}
 - \underbrace{\frac{1}{2\pi}\int K(\al,\beta)\gamma(\beta)D(\beta)d\beta}_{T_2} \nonumber \\
& + \underbrace{\frac{1}{2\pi}c_1(\al)\int \frac{D(\al) - D(\beta)}{4\sin^{2}\left(\frac{\al - \beta}{2}\right)}\gamma(\beta)d\beta}_{T_3}
 + \underbrace{\frac{1}{2\pi}c_2(\al)\int \frac{D(\al) - D(\beta)}{2\tan\left(\frac{\al - \beta}{2}\right)}\gamma(\beta)d\beta}_{T_4},
\label{estimacionessplitting}
\end{align}
where
\begin{align*}
K(\al,\beta) & = \frac{1}{(x(\al)-x(\beta))^{2}} - \frac{c_1(\al)}{4\sin^{2}\left(\frac{\al - \beta}{2}\right)} - \frac{c_2(\al)}{2\tan\left(\frac{\al - \beta}{2}\right)}\\
c_1(\al) & = \frac{1}{x_\al^{2}(\al)}\\
c_2(\al) & = \frac{x_{\al \al}(\al)}{x_{\al}^{3}(\al)}.
\end{align*}

We can think of $c_1(\al)$ and $c_2(\al)$ as the Taylor coefficients of $\Theta(\al,\beta)$ around $\beta = \alpha$. We can bound the terms in
\eqref{estimacionessplitting} in the following way:

\begin{align*}
T_4(\al) & = c_2(\al)[H(D\gamma)(\al) - DH(\gamma)(\al)]\\
T_3(\al) & = c_1(\al)[\Lambda(D\gamma)(\al) - D\Lambda(\gamma)(\al)]
\end{align*}

We have then the estimates

\begin{align*}
\|T_4\|_{L^{2}} & \leq \|c_2\|_{L^{\infty}}(\|D\|_{L^{2}}\|\gamma\|_{L^{\infty}} + \|D\|_{L^{2}}\|H\gamma\|_{L^{\infty}}) \\
\|T_3\|_{L^{2}} & \leq \|c_1\|_{L^{\infty}}(\|D\|_{L^{2}}\|\gamma_{\al}\|_{L^{\infty}} + \|D_{\al}\|_{L^{2}}\|\gamma\|_{L^{\infty}} + \|D\|_{L^{2}}\|\Lambda(\gamma)\|_{L^{\infty}}).
\end{align*}

We now move on to $T_1$. We will estimate it in the following way:
\begin{align*}
\int T_1 \overline{D(\al)} d\al = \frac{1}{2\pi}\int |D(\al)|^{2}\int K(\al,\beta) \gamma(\beta)d\beta d\al
 \leq \frac{1}{2\pi}\|D\|_{L^{2}}^{2}\left\|\int K(\cdot,\beta) \gamma(\beta)d\beta\right\|_{L^{\infty}}.
\end{align*}

To estimate the kernel $T_2$ we will use the Generalized Young's inequality \cite{Folland:introduction-pdes}:

\begin{align*}
\|T_2(D)\|_{L^{2}}^{2} = \frac{1}{4\pi^{2}}\int \int \int K(\al,\beta) \gamma(\beta) D(\beta) \overline{K(\al,\sigma)} \overline{\gamma(\sigma)} \overline{D(\sigma)} d\beta d\sigma d\alpha.
\end{align*}

Defining
\begin{align*}
\tilde{K}(\beta,\sigma) = \int K(\al,\beta) \gamma(\beta) \overline{K(\al,\sigma)} \gamma(\sigma) d\alpha,
\end{align*}
we have that
\begin{align*}
\|T_2(D)\|_{L^{2}}^{2} & = \frac{1}{4\pi^{2}}\int \int \tilde{K}(\beta,\sigma)D(\beta)\overline{D(\sigma)} d\beta d\sigma \\
& = \frac{1}{4\pi^{2}}\int D(\beta) \left(\int \tilde{K}(\beta,\sigma)\overline{D(\sigma)} d\sigma\right)d\beta \\
& \leq \frac{1}{4\pi^{2}}\|D\|_{L^{2}} \left\|\int \tilde{K}(\dot,\sigma)d\sigma\right\|_{L^{2}} \\
& \leq \frac{1}{4\pi^{2}}C\|D\|_{L^{2}}^{2}, \quad C = \max\left\{\max_{\beta}\int |\tilde{K}(\beta,\sigma)|d\sigma,\max_{\sigma} \int |\tilde{K}(\beta,\sigma)|d\beta\right\}
\end{align*}

We finally show how to estimate the Kernels with $b_2 = -1$. We will do this by showing how to estimate $\Theta^{0,0,0,0}_{1,-1}$ but the technique can be applied to any Kernel.

\begin{align*}
\frac{1}{2\pi}\int \Theta^{0,0,0,0}_{1,-1}(d(\beta))d\beta & = \underbrace{\frac{1}{2\pi}\int K(\al,\beta)d(\beta)d\beta}_{T_1} \nonumber \\
& + \underbrace{\frac{1}{2\pi}c_1(\al)\int \frac{1}{2\tan\left(\frac{\al - \beta}{2}\right)}d(\beta)d\beta}_{T_2},
\end{align*}
where
\begin{align*}
K(\al,\beta) & = \frac{1}{(x(\al)-x(\beta))} - \frac{c_1(\al)}{2\tan\left(\frac{\al - \beta}{2}\right)}\\
c_1(\al) & = \frac{1}{x_\al(\al)}.
\end{align*}
We can easily estimate these two terms applying to $T_1$ the same estimates (Young's inequality) as for $T_2$ in the previous case and by noting that $T_2$ is $\frac{1}{2}c_1(\al)H(d)$.

\subsection{Estimates for the linear terms with $Q \neq 1$}

To perform the real estimates, where $Q \neq 1$ we will use the estimates from the previous sections. We will explain how to pass from the former ones to the latter ones. We will illustrate this by computing the linear terms of the Birkhoff-Rott operator.

First of all, the total number of terms will increase by a factor 2, since we will have

\begin{align*}
Q^{2}(z)BR(z,\om) - Q^{2}(x)BR(x,\gamma)) & = \underbrace{(Q^{2}(z)-Q^{2}(x))(BR(z,\om) - BR(x,\gamma))}_{\text{nonlinear}} \\
& + \underbrace{Q^{2}(x)(BR(z,\om) - BR(x,\gamma))}_{\text{calculated before}} \\
& + \underbrace{(Q^{2}(z)-Q^{2}(x))BR(x,\gamma)}_{\text{new terms}}
\end{align*}

In order to calculate the old terms with $Q \neq 1$, the only thing we have to do is to incorporate a factor of $\da^{k}Q^{2}(x)(\al)$ in the estimates. The new terms can easily be calculated using that, up to linear order
\begin{align*}
(Q^{2}(z)-Q^{2}(x)) = \frac{1}{8}\left\langle \frac{1+x^4}{x}, \overline{3x^2 - \frac{1}{x^2}}\right\rangle D + O(D^2).
\end{align*}

\section{Proof of Theorem \ref{stabilitytheorem}}
\label{sectionstability}

In this section, we will prove the stability Theorem \ref{stabilitytheorem}.

%%%%%%%%%%%%%%%%%%%%%%%%%%%%%%%%%%%%%%%%%%%%%%%%%%%%%%%%%%%%%%%%%%%%%%%%%%%5
The equations are:

\begin{align*}
\text{SPLASH}&\left\{
\begin{array}{cl}
z_t & = Q_{z}^{2}BR + cz_{\al} \\
c & = \displaystyle \frac{\al + \pi}{2\pi}\int_{-\pi}^{\pi}(Q^2 BR)_{\al}\frac{z_\al}{|z_{\al}|^{2}} d\al - \int_{-\pi}^{\al}(Q^2 BR)_{\beta}\frac{z_{\beta}}{|z_{\beta}|^{2}}d\beta \\
\omega_t & + 2BR_{t} \cdot z_{\al} = - (Q^2)_{\al}|BR|^{2} + 2cBR_{\al} \cdot z_{\al} + (c\varpi)_{\al}\\
 & \displaystyle - \left(\frac{Q^2\varpi^2}{4|z_{\al}|^{2}}\right)_{\al}
- 2(P^{-1}_{2}(z))_{\al}
\end{array}
\right.\\
\text{APPROX}&\left\{
\begin{array}{cl}
x_t & = Q^2(x)BR(x,\gamma)  + bx_{\al} + f\\
b & = \underbrace{\frac{\al + \pi}{2\pi}\int_{-\pi}^{\pi}(Q^2 BR)_{\al}\frac{x_\al}{|x_{\al}|^{2}} d\al - \int_{-\pi}^{\al}(Q^2 BR)_{\beta}\frac{x_{\al}}{|x_{\al}|^{2}}d\beta}_{b_s} \\
&+ \underbrace{\frac{\al + \pi}{2\pi}\int_{-\pi}^{\pi}f_{\al}\frac{x_\al}{|x_{\al}|^{2}} d\al - \int_{-\pi}^{\al}f_{\beta}\frac{x_{\beta}}{|x_{\beta}|^{2}}d\beta}_{b_e}\\
\gamma_t & + 2BR_{t}(x,\gamma) \cdot x_{\al} = - (Q^2(x))_{\al}|BR(x,\gamma)|^{2} + 2bBR_{\al}(x,\gamma) \cdot x_{\al} + (b\gamma)_{\al}  \\
& \qquad \qquad \qquad - \displaystyle \left(\frac{Q^2(x)\gamma^2}{4|x_{\al}|^{2}}\right)_{\al} - 2(P^{-1}_{2}(x))_{\al} + g
\end{array}
\right.
\end{align*}

where
\begin{equation*}
BR(z,\varpi)(\alpha) = \frac{1}{2\pi}PV\int_{-\pi}^{\pi}\frac{(z(\al) - z(\al-\beta))^{\perp}}{|z(\al) - z(\al-\beta)|^{2}}\varpi(\al-\beta)d\beta,
\end{equation*}

$f$ will be the error for $z$ and $g$ will be the error for $\om$.

%%%%%%%%%%%%%%%%%%%%%%%%%%%%%%%%%%%%%%%%%%%%%%%%%%%%%%%%%%%%%%%%%%%%%%%%%%%%%%%%%%%%%%%
\subsection{Computing the difference $z-x$ and $\om - \gamma$}
%%%%%%%%%%%%%%%%%%%%%%%%%%%%%%%%%%%%%%%%%%%%%%%%%%%%%%%%%%%%%%%%%%%%%%%%%%%%%%%%%%%%%%%

We define now:

$$ D \equiv z - x, \quad d \equiv \om - \gamma, \quad \mathcal{D} \equiv \varphi - \psi$$

The energy

$$ E(t) \equiv \frac{1}{2}\left(\|D\|^{2}_{L^{2}} + \int_{-\pi}^{\pi}\frac{Q_{z}^{2}}{|z_{\al}|^{2}}\sigma_{z}|\partial^{4}_{\al}D|^{2} + \|d\|^{2}_{H^{2}} + \|\mathcal{D}\|^{2}_{H^{3+\frac{1}{2}}}\right)$$

and the Rayleigh-Taylor condition

\begin{align*}\sigma_{z} &\equiv \left(BR_{t} + \frac{\varphi}{|z_{\al}|}BR_{\al}\right) \cdot z_{\al}^{\perp} + \frac{\om}{2|z_{\al}|^{2}}\left(z_{\al t} + \frac{\varphi}{|z_{\al}|}z_{\al \al}\right) \cdot z_{\al}^{\perp} \\
&+ Q\left|BR + \frac{\om}{2|z_{\al}|^{2}}z_{\al}\right|^{2}\nabla Q \cdot z_{\al}^{\perp}
 - (\nabla P_{2}^{-1})(z) \cdot z_{\al}^{\perp}\end{align*}

 Note that $\sigma_{z} > 0$. We shall show that

$$\left|\frac{d}{dt}E(t)\right|\leq \mathcal{C}(t)(E(t)+E^{k}(t))+c\delta(t)$$
where $$\mathcal{C}(t)= \mathcal{C}(\|x\|_{H^{5+\frac12}}(t),\|\gamma\|_{H^{3+\frac12}}(t),
\|\psi\|_{H^{4+\frac12}}(t),\|F(x)\|_{L^\infty}(t))$$ and $$\delta(t)=(\|f\|_{H^{5+\frac12}}(t)+\|g\|_{H^{3+\frac12}}(t))^k+(\|f\|_{H^{5+\frac12}}(t)+\|g\|_{H^{3+\frac12}}(t))^2, \text{ $k$ big enough}$$
depend on the norms of $f$ and $g$.

\begin{remark}
From now on, we will denote $E(t) + E(t)^k$ by $P(E(t))$.
\end{remark}

 $ \frac{1}{2}\frac{d}{dt}\|D\|_{L^{2}}^{2} \leq CP(E(t)) + \delta(t)$ is left to the reader. We compute

 $$ \frac{1}{2}\frac{d}{dt}\int_{-\pi}^{\pi}\frac{Q_{z}^{2}}{|z_{\al}|^{2}}\sigma_{z}|\partial^{4}_{\al}D|^{2}
 = \frac{1}{2}\int_{-\pi}^{\pi}\frac{(Q_{z}^{2}\sigma_{z})_{t}}{|z_{\al}|^{2}}|\partial^{4}_{\al}D|^{2}
 + \int_{-\pi}^{\pi}\frac{Q_{z}^{2}}{|z_{\al}|^{2}}\sigma_{z}\partial_{\al}^{4}D \partial_{\al}^{4}D_t$$

 The first integral is easy to bound by $CP(E(t))$, we proceed as in the local existence Theorem I.7 in \cite{Castro-Cordoba-Fefferman-Gancedo-GomezSerrano:finite-time-singularities-free-boundary-euler}. We split
 $$ I = \int_{-\pi}^{\pi}\frac{Q_{z}^{2}}{|z_{\al}|^{2}}\sigma_{z}\partial_{\al}^{4}D \partial_{\al}^{4}D_t = I_1 + I_2 + I_3$$
 where
 \begin{align*}
 I_1 & = \int_{-\pi}^{\pi}\frac{Q_{z}^{2}}{|z_{\al}|^{2}}\sigma_{z}\partial_{\al}^{4}D \partial_{\al}^{4}(Q_{z}^2 BR(z,\om) - Q_{x}^{2}BR(x,\gamma))d\al\\
 I_2 & = \int_{-\pi}^{\pi}\frac{Q_{z}^{2}}{|z_{\al}|^{2}}\sigma_{z}\partial_{\al}^{4}D \partial_{\al}^{4}(cz_{\al} - bx_{\al})d\al\\
 I_3 & = \int_{-\pi}^{\pi}\frac{Q_{z}^{2}}{|z_{\al}|^{2}}\sigma_{z}\partial_{\al}^{4}D \partial_{\al}^{4}fd\al
 \end{align*}
 We have:
 $$ I_3 \leq \frac{1}{2}\int_{-\pi}^{\pi}\frac{Q_{z}^{2}}{|z_{\al}|^{2}}\sigma_{z}|\partial_{\al}^{4}D|^2d\al
 + \frac{1}{2}\int_{-\pi}^{\pi}\frac{Q_{z}^{2}}{|z_{\al}|^{2}}\sigma_{z}|\partial_{\al}^{4}f|^2d\al
 \leq CP(E(t)) + \frac{\|Q_{z}^{2}\sigma_{z}\|_{L^{\infty}}}{2}\delta(t)$$
 Thus, we are done with $I_3$. We now split
 \begin{align*}
 I_1 & = \text{l.o.t} + I_{1,1} + I_{1,2} + I_{1,3} + I_{1,4}\\
 I_{1,1} & = \int_{-\pi}^{\pi}\frac{Q_{z}^{2}}{|z_{\al}|^{2}}\sigma_{z}\partial_{\al}^{4}D (\partial_{\al}^{4}(Q_{z}^2) BR(z,\om) - \partial_{\al}^{4}(Q_{x}^{2})BR(x,\gamma))d\al\\
 I_{1,2} & = \int_{-\pi}^{\pi}\frac{Q_{z}^{2}}{|z_{\al}|^{2}}\sigma_{z}\partial_{\al}^{4}D \left(Q_{z}^2\frac{1}{2\pi}\int_{-\pi}^{\pi}\frac{(\partial_{\al}^{4}z(\al) - \partial_{\al}^{4}z(\al - \beta))^{\perp}}{|z(\al) - z(\al - \beta)|^{2}}\om(\al-\beta)d\beta\right. \\
 & \left.- Q_{x}^2\frac{1}{2\pi}\int_{-\pi}^{\pi}\frac{(\partial_{\al}^{4}x(\al) - \partial_{\al}^{4}x(\al - \beta))^{\perp}}{|x(\al) - x(\al - \beta)|^{2}}\gamma(\al-\beta)d\beta\right)d\al\\
 I_{1,3} & = \int_{-\pi}^{\pi}\frac{Q_{z}^{2}}{|z_{\al}|^{2}}\sigma_{z}\partial_{\al}^{4}D\\ &\times\left(Q_{z}^2\frac{-1}{\pi}\int_{-\pi}^{\pi}\frac{(z(\al) - z(\al - \beta))^{\perp}}{|z(\al) - z(\al - \beta)|^{4}}(z(\al) - z(\al - \beta)) \cdot (\partial_{\al}^{4}z(\al) - \partial_{\al}^{4}z(\al - \beta))\om(\al-\beta)d\beta\right. \\
 & +\left.Q_{x}^2\frac{1}{\pi}\int_{-\pi}^{\pi}\frac{(x(\al) - x(\al - \beta))^{\perp}}{|x(\al) - x(\al - \beta)|^{4}}(x(\al) - x(\al - \beta)) \cdot (\partial_{\al}^{4}x(\al) - \partial_{\al}^{4}x(\al - \beta))\gamma(\al-\beta)d\beta\right) \\
 I_{1,4} & = \int_{-\pi}^{\pi}\frac{Q_{z}^{2}}{|z_{\al}|^{2}}\sigma_{z}\partial_{\al}^{4}D \partial_{\al}^{4}(Q_{z}^2 BR(z,\partial_{\alpha}^{4}\om) - Q_{x}^{2}BR(x,\partial_{\al}^{4}\gamma))d\al
 \end{align*}
 where l.o.t stands for low order terms, nice terms easier to deal with.

 $$ I_{1,1} = \text{l.o.t} + I_{1,1,1} \text{ where }$$

 \begin{align*}
 I_{1,1,1} & = 2\int_{-\pi}^{\pi}\frac{Q_{z}^{2}}{|z_{\al}|^{2}}\sigma_{z}\partial_{\al}^{4}D (\nabla Q(z) \cdot \partial_{\al}^{4} z BR(z,\om) - \nabla Q(x) \cdot \partial_{\al}^{4} x BR(x,\gamma))d\al \\
 & = 2\int_{-\pi}^{\pi}\frac{Q_{z}^{2}}{|z_{\al}|^{2}}\sigma_{z}\partial_{\al}^{4}D \nabla Q(z) \cdot \partial_{\al}^{4} D BR(z,\om)d\al \\
  & + 2\int_{-\pi}^{\pi}\frac{Q_{z}^{2}}{|z_{\al}|^{2}}\sigma_{z}\partial_{\al}^{4}D (\nabla Q(z) \cdot \partial_{\al}^{4} x BR(z,\om) - \nabla Q(x) \cdot \partial_{\al}^{4} x BR(x,\gamma))d\al \\
  & \leq 2\int_{-\pi}^{\pi}\frac{Q_{z}^{2}}{|z_{\al}|^{2}}\sigma_{z}|\partial_{\al}^{4}D|^2 d\al \underbrace{\|\nabla Q(z) BR(z,\om)\|_{L^\infty}}_{\text{bounded as for local existence}} \\
  & + \int_{-\pi}^{\pi}\frac{Q_{z}^{2}}{|z_{\al}|^{2}}\sigma_{z}|\partial_{\al}^{4}D|^2
  + \underbrace{\int_{-\pi}^{\pi}\frac{Q_{z}^{2}}{|z_{\al}|^{2}}\sigma_{z}|\nabla Q(z) \cdot \partial_{\al}^{4} x BR(z,\om) - \nabla Q(x) \cdot \partial_{\al}^{4} x BR(x,\gamma)|^2 d\al}_{\text{l.o.t in $D$ and $d$}} \\
  & \leq CP(E(t))
 \end{align*}

 which means $I_{1,1}$ is done.

 From now on we will denote

 $$ \Delta_{\beta}z(\al) = z(\al) - z(\al - \beta)$$

 \begin{align*}
 I_{1,2} & = I_{1,2,1} + I_{1,2,2} + I_{1,2,3} + I_{1,2,4} \text{ where}\\
 I_{1,2,1} & = \int_{-\pi}^{\pi}\frac{Q_{z}^{2}}{|z_{\al}|^{2}}\sigma_{z}\partial_{\al}^{4}D Q_{z}^2\frac{1}{2\pi}\int_{-\pi}^{\pi}\frac{\Delta_{\beta} \partial_{\al}^{4} D^{\perp}(\al)}{|\Delta_{\beta} z(\al)|^2}\om(\al-\beta)d\beta d\al \\
 I_{1,2,2} & = \int_{-\pi}^{\pi}\frac{Q_{z}^{2}}{|z_{\al}|^{2}}\sigma_{z}\partial_{\al}^{4}D Q_{z}^2\frac{1}{2\pi}\int_{-\pi}^{\pi}\Delta_{\beta}\partial_{\al}^{4}x ^{\perp}(\al) \left(\frac{1}{|\Delta_{\beta}z(\al)|^{2}} - \frac{1}{|\Delta_{\beta}x(\al)|^{2}}\right)\om(\al-\beta)d\beta d\al \\
 I_{1,2,3} & = \int_{-\pi}^{\pi}\frac{Q_{z}^{2}}{|z_{\al}|^{2}}\sigma_{z}\partial_{\al}^{4}D Q_{z}^2\frac{1}{2\pi}\int_{-\pi}^{\pi}\frac{\Delta_{\beta}\partial_{\al}^{4}x ^{\perp}(\al)}{|\Delta_{\beta}x(\al)|^{2}}d(\al-\beta)d\beta d\al \\
 I_{1,2,4} & = \int_{-\pi}^{\pi}\frac{Q_{z}^{2}}{|z_{\al}|^{2}}\sigma_{z}\partial_{\al}^{4}D (Q_{z}^2 - Q_{x}^{2})\frac{1}{2\pi}\int_{-\pi}^{\pi}\frac{\Delta_{\beta}\partial_{\al}^{4}x ^{\perp}(\al)}{|\Delta_{\beta}x(\al)|^{2}}\gamma(\al-\beta)d\beta d\al \\
 \end{align*}

 \begin{align*}
 I_{1,2,1} & = \int_{-\pi}^{\pi}\frac{Q_{z}^{4}}{|z_{\al}|^{2}}\sigma_{z}\partial_{\al}^{4}D(\al)\frac{1}{2\pi}\int_{-\pi}^{\pi}\frac{\Delta_{\al-\beta} \partial_{\al}^{4} D^{\perp}(\al)}{|\Delta_{\al-\beta} z(\al)|^2}\om(\al-\beta)d\beta d\al \\
& = \frac{1}{|z_{\al}|^{2}}\frac{1}{2\pi}\int_{-\pi}^{\pi}\int_{-\pi}^{\pi}\partial_{\al}^{4}D\frac{\Delta_{\al-\beta} \partial_{\al}^{4} D^{\perp}(\al)}{|\Delta_{\al-\beta} z(\al)|^2}\left(\frac{Q_{z}^{4}(\al) \sigma_{z}(\al) \om(\beta) - Q_{z}^{4}(\beta) \sigma_{z}(\beta) \om(\al)}{2}\right. \\
 &\left. + \underbrace{\frac{Q_{z}^{4}(\al) \sigma_{z}(\al) \om(\beta) + Q_{z}^{4}(\beta) \sigma_{z}(\beta) \om(\al)}{2}}_{\text{\tiny this is zero as in local existence $(\partial_{\al}^{4}D \cdot \partial_{\al}^{4}D^{\perp} = 0)$}}\right)d\al d\beta \\
 \Rightarrow I_{1,2,1} & = \frac{1}{2\pi}\int_{-\pi}^{\pi}\frac{\partial_{\al}^{4}D}{|z_{\al}|^{2}}\int_{-\pi}^{\pi}\underbrace{\frac{\Delta_{\al-\beta} \partial_{\al}^{4} D^{\perp}(\al)}{|\Delta_{\al-\beta} z(\al)|^2}}_{\substack{\text{\tiny Hilbert transform}\\ \text{\tiny applied to $\partial_{\al}^{4}D^{\perp}(\al)$}}}\left(\frac{Q_{z}^{4}(\al) \sigma_{z}(\al) \om(\beta) - Q_{z}^{4}(\beta) \sigma_{z}(\beta) \om(\al)}{2} \right) \\
 \Rightarrow I_{1,2,1} & \leq CP(E(t))
 \end{align*}

 For $I_{1,2,2}$ we can make a trick to get less derivatives in $x$.
 \begin{align*}
 I_{1,2,2} & = I_{1,2,2}^{1} + I_{1,2,2}^{2} + I_{1,2,2}^{3} \\
 I_{1,2,2}^{3} & = \frac{1}{2}\int_{-\pi}^{\pi}\frac{Q_{z}^{4}}{|z_{\al}|^{2}}\sigma_{z}\partial_{\al}^{4}D
 \om(\al)\left(\frac{1}{|z_{\al}|^{2}} - \frac{1}{|x_{\al}|^{2}}\right)
\overbrace{\frac{1}{\pi}\int_{-\pi}^{\pi}\frac{\Delta_{\beta}\partial_{\al}^{4}x ^{\perp}(\al)}{\beta^{2}}d\beta}^{\Lambda \partial_{\al}^{4} x} d\al \\
I_{1,2,2}^{2} & = \frac{1}{2\pi}\int_{-\pi}^{\pi}\frac{Q_{z}^{4}}{|z_{\al}|^{2}}\sigma_{z}\partial_{\al}^{4}D \int_{-\pi}^{\pi}\Delta_{\beta}\partial_{\al}^{4}x ^{\perp}(\al) \left(\frac{1}{|\Delta_{\beta}z(\al)|^{2}}
- \frac{1}{|z_{\al}(\al)|^{2}\beta^{2}} + \overbrace{\frac{z_{\al} \cdot z_{\al \al}}{|z_{\al}|^{4} \beta}}^{=0}\right. \\
& \left.-\left(\frac{1}{|\Delta_{\beta}x(\al)|^{2}}
- \frac{1}{|x_{\al}(\al)|^{2}\beta^{2}} + \overbrace{\frac{x_{\al} \cdot x_{\al \al}}{|x_{\al}|^{4} \beta}}^{=0}\right)
\right)\om(\al)d\beta d\al\\
I_{1,2,2}^{1} & = \frac{1}{2\pi}\int_{-\pi}^{\pi}\frac{Q_{z}^{4}}{|z_{\al}|^{2}}\sigma_{z}\partial_{\al}^{4}D \int_{-\pi}^{\pi}\Delta_{\beta}\partial_{\al}^{4}x ^{\perp}(\al) \left(\frac{1}{|\Delta_{\beta}z(\al)|^{2}}
- \frac{1}{|\Delta_{\beta}x(\al)|^{2}}\right)\left(\om(\al - \beta) - \om(\al)\right)d\beta d\al
 \end{align*}

 We use that $\displaystyle \left|\frac{1}{|z_{\al}|^{2}} - \frac{1}{|x_{\al}|^{2}}\right| \leq \frac{|x_{\al}| + |z_{\al}|}{|z_{\al}|^{2}|x_{\al}|^{2}}|D_{\al}|$ to find that
 \begin{align*} I_{1,2,2}^{3} & \leq \frac{1}{4}\int_{-\pi}^{\pi}\frac{Q_{z}^{2}}{|z_{\al}|^{2}}\sigma_{z}|\partial_{\al}^{4}D|^2
 + \|Q_{z}\|_{L^{\infty}}^{6} \|\sigma_{z}\|_{L^{\infty}} \|\om\|_{L^{\infty}}^{2}\left(\frac{|x_{\al}| + |z_{\al}|}{|z_{\al}|^{2}|x_{\al}|^{2}}\right)^{2}\overbrace{\|D_{\al}\|_{L^{\infty}}^{2}}^{ \substack{\text{\tiny Sobolev}\\ \text{\tiny inequalities}}}\overbrace{\|\Lambda \partial_{\al}^{4} x\|_{L^{2}}^{2}}^{\text{Control of $\|x\|_{H^{5}}$}}\\
 & \leq CP(E(t))
 \end{align*}

 We can use that
 $$ \left|\left(\frac{1}{|\Delta_{\beta}z(\al)|^{2}}
- \frac{1}{|z_{\al}(\al)|^{2}\beta^{2}} + \frac{z_{\al} \cdot z_{\al \al}}{|z_{\al}|^{4} \beta}\right)\right| \leq \|z\|_{C^{2}}^{k} \frac{1}{\beta^{1/2}}\|z\|_{C^{2+\frac{1}{2}}}\|F(z)\|_{L^{\infty}}^{k}
$$

and that
\begin{align*}
& \left|\frac{1}{|\Delta_{\beta}z(\al)|^{2}}
- \frac{1}{|z_{\al}(\al)|^{2}\beta^{2}} + \frac{z_{\al} \cdot z_{\al \al}}{|z_{\al}|^{4} \beta}
 - \left(\frac{1}{|\Delta_{\beta}x(\al)|^{2}}
- \frac{1}{|x_{\al}(\al)|^{2}\beta^{2}} + \frac{x_{\al} \cdot x_{\al \al}}{|x_{\al}|^{4} \beta}\right)\right| \\
& \leq \|z\|_{C^{2}}^{k} \|x\|_{C^{2}}^{k}\frac{1}{\beta^{1/2}}\|D\|_{C^{2+\frac{1}{2}}}\|F(z)\|_{L^{\infty}}^{k}\|F(x)\|_{L^{\infty}}^{k}
\end{align*}

to find
\begin{align*}
 I_{1,2,2}^{2} &\leq \frac{1}{8\pi^{2}}\int_{-\pi}^{\pi}\frac{Q_{z}^{2}}{|z_{\al}|^{2}}\sigma_{z}|\partial_{\al}^{4}D|^2
\\&+ C\|Q_{z}\|_{L^{\infty}}^{6} \|\sigma_{z}\|_{L^{\infty}}\|z\|_{C^{2}}^{k} \|x\|_{C^{2}}^{k}\|D\|_{C^{2+\frac{1}{2}}}\|\partial_{\al}^{4} x\|_{L^{2}}^{2}\|F(z)\|_{L^{\infty}}^{k}\|F(x)\|_{L^{\infty}}^{k}\end{align*}

We've used that
$$ \left(\int_{-\pi}^{\pi}d\al\left(\int_{-\pi}^{\pi}\frac{\da^{4}x(\al - \beta)}{|\beta|^{1/2}}d\beta\right)^{2}\right)^{1/2} \leq C\|\partial_{\al}^{4}x\|_{L^{2}}.$$

We split further in $I_{1,2,2}^{1} = I_{1,2,2}^{1,1} + I_{1,2,2}^{1,2}$:

\begin{align*}
I_{1,2,2}^{1,1} & = \frac{1}{2\pi}\int_{-\pi}^{\pi}\frac{Q_{z}^{4}}{|z_{\al}|^{2}}\sigma_{z}\partial_{\al}^{4}D \int_{-\pi}^{\pi}\Delta_{\beta}\partial_{\al}^{4}x ^{\perp}(\al) \left(\frac{1}{|\Delta_{\beta}z(\al)|^{2}}
- \frac{1}{|\Delta_{\beta}x(\al)|^{2}}\right)\\
&\times \left(\om(\al - \beta) - \om(\al) + \om_{\al}(\al) \beta\right)d\beta d\al\\
I_{1,2,2}^{1,2} & = \frac{1}{2\pi}\int_{-\pi}^{\pi}\frac{Q_{z}^{4}}{|z_{\al}|^{2}}\sigma_{z}\partial_{\al}^{4}D  \om_{\al}(\al) \int_{-\pi}^{\pi}\Delta_{\beta}\partial_{\al}^{4}x ^{\perp}(\al)\left(\frac{\beta}{|\Delta_{\beta}z(\al)|^{2}}
- \frac{\beta}{|\Delta_{\beta}x(\al)|^{2}}\right)d\beta d\al\\
\end{align*}

Inside of the $\beta$ integral in $I_{1,2,2}^{1,1}$ there is no principal value, so the appropriate estimate follows:

$$ I_{1,2,2}^{1,1} \leq CP(E(t))$$

For $I_{1,2,2}^{1,2}$ we proceed as for $I_{1,2,2}^{2}$. We decompose adding and subtracting $\displaystyle \frac{1}{|z_{\al}|^{2}\beta} - \frac{1}{|x_{\al}|^{2}\beta}$. Thus, we are done with $I_{1,2,2}$. We decompose $I_{1,2,3} = I_{1,2,3}^{1} + I_{1,2,3}^{2} + I_{1,2,3}^{3}$.

\begin{align*}
I_{1,2,3}^{1} & = \int_{-\pi}^{\pi}\frac{Q_{z}^{2}}{|z_{\al}|^{2}}\sigma_{z}\partial_{\al}^{4}D Q_{z}^2\frac{1}{2\pi}\int_{-\pi}^{\pi}\Delta_{\beta}\partial_{\al}^{4}x ^{\perp}(\al)\\
&\times\left(\frac{1}{|\Delta_{\beta}x(\al)|^{2}}-\frac{1}{|x_\al|^{2}\beta^{2}}+\frac{x_{\al} \cdot x_{\al \al}}{|x_{\al}|^{4}\beta}\right)d(\al-\beta)d\beta d\al \\
I_{1,2,3}^{2} & = - \int_{-\pi}^{\pi}\frac{Q_{z}^{2}}{|z_{\al}|^{2}}\sigma_{z}\partial_{\al}^{4}D Q_{z}^2\frac{\da^{4} x^{\perp}(\al)}{|x_{\al}|^{2}}\frac{1}{2\pi}\int_{-\pi}^{\pi}\frac{\Delta_{\beta}d(\al)}{\beta^{2}}d\beta d\al \\
I_{1,2,3}^{3} & = \int_{-\pi}^{\pi}\frac{Q_{z}^{2}}{|z_{\al}|^{2}}\sigma_{z}\partial_{\al}^{4}D Q_{z}^2\frac{1}{|x_{\al}|^{2}}\frac{1}{2\pi}\int_{-\pi}^{\pi}\frac{\Delta_{\beta}(d\da^{4}x^{\perp})(\al)}{\beta^{2}}d\beta d\al
\end{align*}

It's easy to obtain:

\begin{align*}
I_{1,2,3}^{1} & \leq \frac{1}{4\pi}\int_{-\pi}^{\pi}\frac{Q_{z}^{2}}{|z_{\al}|^{2}}\sigma_{z}|\partial_{\al}^{4}D|^{2}d\al
+ C\|Q_{z}\|_{L^{\infty}}^{6} \|\sigma_{z}\|_{L^{\infty}}\|d\|_{L^{\infty}} \|x\|_{C^{2}}^{k}\|F(x)\|_{L^{\infty}}^{k}\|x\|_{C^{2,\delta}}\|\da^{4} x\|_{L^{2}}^{2}\\& \leq CP(E(t)) \\
 I_{1,2,3}^{2}& \leq CP(E(t)) \text{ analogously since } \|\Lambda d\|_{L^{\infty}} \leq C\|d\|_{H^{2}}\\
 I_{1,2,3}^{3}& \leq CP(E(t)) \text{ using } \|\Lambda (d \da^{4} x^{\perp})\|_{L^{2}} \leq C\|d\|_{H^{2}}\|x\|_{H^{5}}.
\end{align*}
We are done with $I_{1,2,3}$. To deal with $I_{1,2,4}$ se use that

$$ Q_{z}^{2} - Q_{x}^{2} = 2Q((1-t)z + tx)\nabla Q((1-t)z + tx) \cdot D(\al) \text{ for } t \in (0,1).$$

Then it is easy to find

$$ I_{1,2,4} \leq CP(E(t)),$$

and we are done with $I_{1,2}$. We decompose $I_{1,3}$ as

\begin{align*}
I_{1,3} & = I_{1,3,1} + I_{1,3,2} +I_{1,3,3} +I_{1,3,4} +I_{1,3,5} +I_{1,3,6} \\
 I_{1,3,1} & = \int_{-\pi}^{\pi}\frac{Q_{z}^{2}}{|z_{\al}|^{2}}\sigma_{z}\partial_{\al}^{4}D Q_{z}^2\frac{-1}{\pi}\int_{-\pi}^{\pi}\frac{\Delta_{\beta}z^{\perp}(\al)}{|\Delta_{\beta}z(\al)|^4}\Delta_{\beta} z(\al) \cdot \Delta_{\beta} \da^{4}D(\al) \om(\al-\beta) d\beta d\al\\
 I_{1,3,2} & = \int_{-\pi}^{\pi}\frac{Q_{z}^{2}}{|z_{\al}|^{2}}\sigma_{z}\partial_{\al}^{4}D Q_{z}^2\frac{-1}{\pi}\int_{-\pi}^{\pi}\frac{\Delta_{\beta}z^{\perp}(\al)}{|\Delta_{\beta}z(\al)|^4}\Delta_{\beta} z(\al) \cdot \Delta_{\beta} \da^{4}x(\al) d(\al-\beta) d\beta d\al\\
 I_{1,3,3} & = \int_{-\pi}^{\pi}\frac{Q_{z}^{2}}{|z_{\al}|^{2}}\sigma_{z}\partial_{\al}^{4}D Q_{z}^2\frac{-1}{\pi}\int_{-\pi}^{\pi}\frac{\Delta_{\beta}z^{\perp}(\al)}{|\Delta_{\beta}z(\al)|^4}\Delta_{\beta} D \cdot \Delta_{\beta} \da^{4}x(\al) \gamma(\al-\beta) d\beta d\al\\
 I_{1,3,4} & = \int_{-\pi}^{\pi}\frac{Q_{z}^{2}}{|z_{\al}|^{2}}\sigma_{z}\partial_{\al}^{4}D Q_{z}^2\frac{-1}{\pi}\int_{-\pi}^{\pi}\frac{\Delta_{\beta}D^{\perp}(\al)}{|\Delta_{\beta}z(\al)|^4}\Delta_{\beta} x(\al) \cdot \Delta_{\beta} \da^{4}x(\al) \gamma(\al-\beta) d\beta d\al\\
 I_{1,3,5} & = \int_{-\pi}^{\pi}\frac{Q_{z}^{2}}{|z_{\al}|^{2}}\sigma_{z}\partial_{\al}^{4}D Q_{z}^2\frac{-1}{\pi}\int_{-\pi}^{\pi}
 \Delta_{\beta} x^{\perp}(\al)\Delta_{\beta} x(\al) \cdot \Delta_{\beta} \da^{4}x(\al) \gamma(\al-\beta)\\
 &\times \left(
 \frac{1}{|\Delta_{\beta}z(\al)|^4} - \frac{1}{|\Delta_{\beta}x(\al)|^4}\right)d\beta d\al\\
 I_{1,3,6} & = \int_{-\pi}^{\pi}\frac{Q_{z}^{2}}{|z_{\al}|^{2}}\sigma_{z}\partial_{\al}^{4}D (Q_{z}^2-Q_{x}^{2})\frac{-1}{\pi}\int_{-\pi}^{\pi}
 \frac{\Delta_{\beta}x^{\perp}(\al)}{|\Delta_{\beta}x(\al)|^4}
 \Delta_{\beta} x(\al) \cdot \Delta_{\beta} \da^{4}x(\al) \gamma(\al-\beta)d\beta d\al\\
 \end{align*}

 $I_{1,3,j}, \quad j = 2,3,4,5,6$ are easier to deal with (It can be done as before). Therefore we focus on $I_{1,3,1}$.
 \begin{align*}
 I_{1,3,1} & = I_{1,3,1}^{1} + I_{1,3,1}^{2} + I_{1,3,1}^{3} \\
 I_{1,3,1}^{1} & = \int_{-\pi}^{\pi}\frac{Q_{z}^{2}}{|z_{\al}|^{2}}\sigma_{z}\partial_{\al}^{4}D Q_{z}^2\frac{-1}{\pi}\int_{-\pi}^{\pi}\left(\frac{\Delta_{\beta}z^{\perp}(\al)}{|\Delta_{\beta}z(\al)|^4}\Delta_{\beta} z(\al) \om(\al-\beta)
 \right.\\
 &\left.- \frac{\da z^{\perp}(\al)}{|\da z(\al)|^4}\da z(\al - \beta)\om(\al)\frac{1}{\beta^{2}}
 \right) \cdot \Delta_{\beta} \da^{4}D(\al) d\beta d\al\\
 I_{1,3,1}^{2} & = \int_{-\pi}^{\pi}\frac{Q_{z}^{2}}{|z_{\al}|^{2}}\sigma_{z}\partial_{\al}^{4}D Q_{z}^2\frac{-1}{\pi}
 \frac{\da z^{\perp}(\al)}{|\da z(\al)|^4}\om(\al) \da^{4}D(\al) \cdot
 \int_{-\pi}^{\pi} \frac{\da z(\al - \beta) - \da z(\al)}{\beta^{2}}d\beta d\al\\
 I_{1,3,1}^{3} & = \int_{-\pi}^{\pi}\frac{Q_{z}^{2}}{|z_{\al}|^{2}}\sigma_{z}\partial_{\al}^{4}D Q_{z}^2\frac{-1}{\pi}
 \frac{\da z^{\perp}(\al)}{|\da z(\al)|^4}\om(\al)
 \int_{-\pi}^{\pi} \frac{\Delta_{\beta}(\da z \cdot \da^{4} D)(\al)}{\beta^{2}}d\beta d\al\\
 \end{align*}

 In $I_{1,3,1}^{1}$ we find a commutator, which can be handled as before. It is also easy to estimate $I_{1,3,1}^{2}$.

 To deal with $I_{1,3,1}^{3}$ we remember that

\begin{align*}
\da z(\al) \cdot \da^{4}D(\al) & = \da z(\al) \cdot \da^{4} z(\al) - \da x(\al) \cdot \da^{4} x(\al) - \da D(\al) \cdot \da^{4} x(\al)\\
& = -3\da^{2}z(\al) \cdot \da^{3} z(\al) + 3 \da^{2} x(\al) \da^{3} x(\al) - \da D(\al) \cdot \da^{4} x(\al)
\end{align*}

That allows us to decompose further
\begin{align*}
\da z(\al) \cdot \da^{4} D(\al) & = -3\da^{2}z(\al) \cdot \da^{3} D(\al) - 3 \da^{2} D(\al) \da^{3} x(\al) - \da D(\al) \cdot \da^{4} x(\al)
\end{align*}

which yields
\begin{align*}
I_{1,3,1}^{3} & = I_{1,3,1}^{3,1} + I_{1,3,1}^{3,2} + I_{1,3,1}^{3,3} \\
I_{1,3,1}^{3,1} & = \frac{3}{\pi}\int_{-\pi}^{\pi}\frac{Q_{z}^{4}}{|z_{\al}|^{2}}\sigma_{z}\partial_{\al}^{4}D
\cdot
 \frac{\da z^{\perp}(\al)}{|\da z(\al)|^4}\om(\al)
 \int_{-\pi}^{\pi} \frac{\Delta_{\beta}(\da^2 z \cdot \da^{3} D)(\al)}{\beta^{2}}d\beta d\al\\
I_{1,3,1}^{3,2} & = \frac{3}{\pi}\int_{-\pi}^{\pi}\frac{Q_{z}^{4}}{|z_{\al}|^{2}}\sigma_{z}\partial_{\al}^{4}D
\cdot
 \frac{\da z^{\perp}(\al)}{|\da z(\al)|^4}\om(\al)
 \int_{-\pi}^{\pi} \frac{\Delta_{\beta}(\da^2 D \cdot \da^{3} x)(\al)}{\beta^{2}}d\beta d\al\\
I_{1,3,1}^{3,3} & = \frac{3}{\pi}\int_{-\pi}^{\pi}\frac{Q_{z}^{4}}{|z_{\al}|^{2}}\sigma_{z}\partial_{\al}^{4}D
\cdot
 \frac{\da z^{\perp}(\al)}{|\da z(\al)|^4}\om(\al)
 \int_{-\pi}^{\pi} \frac{\Delta_{\beta}(\da D \cdot \da^{4} x)(\al)}{\beta^{2}}d\beta d\al\\
\end{align*}

We use that
\begin{align*}
\left\|\int_{-\pi}^{\pi} \frac{\Delta_{\beta}(\da^2 z \cdot \da^{3} D)(\al)}{\beta^{2}}d\beta\right\|_{L^{2}}^{2} \leq C\left\|\da (\da^{2} z \cdot \da^{3} D)\right\|_{L^{2}}^{2} \leq CP(E(t))
\end{align*}

to control $I_{1,3,1}^{3,1}$. $I_{1,3,1}^{3,2}$ follows similarly. We control $I_{1,3,1}^{3,3}$ using that
\begin{align*}
\left\|\int_{-\pi}^{\pi} \frac{\Delta_{\beta}(\da D \cdot \da^{4} x)(\al)}{\beta^{2}}d\beta\right\|_{L^{2}}^{2} &\leq
\left\|\da (\da D \cdot \da^{4} x)\right\|_{L^{2}}^{2} \\& \leq \|\da D\|_{L^{\infty}}^{2}\|\da^{5} x\|_{L^{2}}^{2} + \|\da^{2} D\|_{L^{\infty}}^{2}\|\da^{4} x\|_{L^{2}}^{2}
\leq CP(E(t))
\end{align*}

This allows us to finish the estimates for $I_{1,3,1}^{3,3}$ and $I_{1,3,1}^{3}$. We are done with $I_{1,3,1}$ and $I_{1,3}$. We now decompose $I_{1,4}$.

\begin{align*}
I_{1,4} &= I_{1,4,1} + I_{1,4,2} + I_{1,4,3} + I_{1,4,4} \\
I_{1,4,1} &= \int_{-\pi}^{\pi}\frac{Q_{z}^{2}}{|z_{\al}|^{2}}\sigma_{z}\partial_{\al}^{4}D \cdot Q_{z}^{2} BR(z,\da^{4} d)d\beta d\al\\
I_{1,4,2} &= \int_{-\pi}^{\pi}\frac{Q_{z}^{2}}{|z_{\al}|^{2}}\sigma_{z}\partial_{\al}^{4}D \cdot Q_{z}^{2} \frac{1}{2\pi}\int_{-\pi}^{\pi}\frac{\Delta_{\beta}D^{\perp}(\al)}{|\Delta_{\beta}z(\al)|^{2}}\da^{4}\gamma(\al - \beta) d\beta d\al\\
I_{1,4,3} &= \int_{-\pi}^{\pi}\frac{Q_{z}^{2}}{|z_{\al}|^{2}}\sigma_{z}\partial_{\al}^{4}D \cdot Q_{z}^{2} \frac{1}{2\pi}\int_{-\pi}^{\pi}\Delta_{\beta}x^{\perp}(\al) \da^{4}\gamma(\al - \beta) \left(\frac{1}{|\Delta_{\beta}z(\al)|^{2}}
- \frac{1}{|\Delta_{\beta}x(\al)|^{2}}\right) d\beta d\al\\
I_{1,4,4} &= \int_{-\pi}^{\pi}\frac{Q_{z}^{2}}{|z_{\al}|^{2}}\sigma_{z}\partial_{\al}^{4}D \cdot (Q_{z}^{2} -Q_{x}^{2})BR(x,\da^{4}\gamma)d\al
\end{align*}

We control $I_{1,4,2}$, $I_{1,4,3}$ and $I_{1,4,4}$ as before. We further split
\begin{align*}
I_{1,4,1} &= I_{1,4,1}^{1} + I_{1,4,1}^{2} + I_{1,4,3} + I_{1,4,4} \\
I_{1,4,1}^{1} & = \int_{-\pi}^{\pi}\frac{Q_{z}^{2}}{|z_{\al}|^{2}}\sigma_{z}\partial_{\al}^{4}D \cdot Q_{z}^{2}\left(BR(z,\da^{4} d) - \frac{1}{2}\frac{\da z^{\perp}(\al)}{|\da z(\al)|^{2}}H(\da^{4} d))\right)d\al\\
I_{1,4,1}^{2} & = \int_{-\pi}^{\pi}\frac{Q_{z}^{2}}{|z_{\al}|^{2}}\sigma_{z}\partial_{\al}^{4}D \cdot \frac{\da z^{\perp}(\al)}{|z_{\al}|}\left(\frac{Q_z^2}{2|z_{\al}|}H(\da^{4} d) - H\left(\frac{Q_{z}^{2}}{2|z_{\al}|}\da^{4} d\right)\right)d\al \\
I_{1,4,1}^{3} & = -\int_{-\pi}^{\pi}\frac{Q_{z}^{2}}{|z_{\al}|^{2}}\sigma_{z}\partial_{\al}^{4}D \cdot \frac{\da z^{\perp}(\al)}{|z_{\al}|}H\left(\left(\frac{Q_{z}^{2}}{2|z_{\al}|} - \frac{Q_{x}^2}{2|x_{\al}|}\right)\da^{4}\gamma\right)d\al \\
I_{1,4,1}^{4} & = \int_{-\pi}^{\pi}\frac{Q_{z}^{2}}{|z_{\al}|^{2}}\sigma_{z}\partial_{\al}^{4}D \cdot \frac{\da z^{\perp}(\al)}{|z_{\al}|}H\left(\frac{Q_{z}^{2}\da^{4}\om}{2|z_{\al}|} - \frac{Q_{x}^2\da^{4}\gamma}{2|x_{\al}|}\right)d\al \\
\end{align*}

There are commutators in $I_{1,4,1}^{1}$ and $I_{1,4,1}^{2}$ so they are easy to estimate. To get the estimate for $I_{1,4,1}^{3}$ we bound

\begin{align*}
\left\|H\left(\left(\frac{Q_{z}^{2}}{2|z_{\al}|} - \frac{Q_{x}^{2}}{2|x_{\al}|}\right)\da^{4}\gamma\right)\right\|_{L^{2}}^{2} \leq \underbrace{\left\|\frac{Q_{z}^{2}}{2|z_{\al}|} - \frac{Q_{x}^{2}}{2|x_{\al}|}\right\|_{L^{\infty}}^{2}}_{\text{at the level of $D(\al)$}}\left\|\da^{4} \gamma\right\|_{L^{2}}^{2} \leq CE^2(t)
\end{align*}

We now remember the following formulas:

\begin{align*}
\varphi & = \frac{Q_{z}^{2}\om}{2|z_{\al}|} - c|z_{\al}| \\
\psi & = \frac{Q_{x}^{2}\gamma}{2|x_{\al}|} - b_{s}|x_{\al}| \\
\end{align*}

These yield

\begin{align*}
I_{1,4,1}^{4} = S + I_{1,4,1}^{4,1} + I_{1,4,1}^{4,2} + \text{ l.o.t},
\end{align*}

where

\begin{align*}
S &= \int_{-\pi}^{\pi}Q_{z}^{2}\sigma_{z} \da^{4}D \cdot \frac{\da z^{\perp}(\al)}{|z_{\al}|^3}H(\da^{4} \Dcal)(\al) d\al\\
I_{1,4,1}^{4,1} &= - \int_{-\pi}^{\pi}\frac{Q_{z}^{2}}{|z_{\al}|^{2}}\sigma_{z}\partial_{\al}^{4}D \cdot \frac{\da z^{\perp}(\al)}{|z_{\al}|}H\left(\frac{Q_{z}\nabla Q(z) \cdot \da^{4} z}{|z_{\al}|}\om - \frac{Q_{x} \nabla Q(x) \cdot \da^{4}x}{|x_{\al}|}\gamma\right)d\al \\
I_{1,4,1}^{4,2} &=  \int_{-\pi}^{\pi}\frac{Q_{z}^{2}}{|z_{\al}|^{2}}\sigma_{z}\partial_{\al}^{4}D \cdot \frac{\da z^{\perp}(\al)}{|z_{\al}|}H\left(\da^{4}\left(c|z_{\al}| - b_{s}|x_{\al}|\right)\right)d\al \\
\end{align*}

$S$ is going to appear later with a negative sign and therefore cancel out. $I_{1,4,1}^{4,1}$ can be bounded as before since it is low order.

We show how to deal with $I_{1,4,1}^{4,2}$. We compute

$$\da^{4}(c|z_{\al}|) = - \da^{3}\left((Q_{z}^{2} BR)_{\al} \cdot \frac{z_{\al}}{|z_{\al}|}\right); \quad \da^{4}(b_s|x_{\al}|) = -\da^{3}\left((Q^{2}_{x}BR)_{\al} \frac{x_{\al}}{|x_{\al}|}\right)$$

Then, in $\da^{4}(c|z_{\al}|) - \da^{4}(b_s|x_{\al}|)$ we consider the most singular terms

\begin{align*}
\da^{4}(c|z_{\al}|) - \da^{4}(b_s|x_{\al}|) &= J_1 + J_2 + J_3 + J_4 + J_5 + \text{ l.o.t.} \\
J_1 & = -2Q_{z} \nabla Q(z) \cdot \da^{4} z BR(z,\om) \cdot \frac{z_{\al}}{|z_{\al}|}
+ 2Q_{x} \nabla Q(x) \cdot \da^{4} x BR(x,\gamma) \cdot \frac{x_{\al}}{|x_{\al}|}  \\
J_2 & = - (Q^{2}_{z}BR)_{\al} \frac{\da^{4} z}{|z_{\al}|} + (Q_{x}^{2} BR)_{\al} \frac{\da^{4} x}{|x_{\al}|} \\
J_3 & = Q^{2}_{z} \frac{1}{\pi}\int_{-\pi}^{\pi}\frac{\Delta_{\beta}z^{\perp}(\al)}{|\Delta_{\beta} z(\al)|^{4}} \frac{z_{\al}(\al)}{|z_{\al}|} \Delta_{\beta} z(\al) \cdot \Delta_{\beta} \da^{4} z(\al) \om(\al - \beta) d\beta \\
& - Q^{2}_{x} \frac{1}{\pi}\int_{-\pi}^{\pi}\frac{\Delta_{\beta}x^{\perp}(\al)}{|\Delta_{\beta} x(\al)|^{4}} \frac{x_{\al}(\al)}{|x_{\al}|} \Delta_{\beta} x(\al) \cdot \Delta_{\beta} \da^{4} x(\al) \gamma(\al - \beta) d\beta \\
J_4 & = -Q_{z}^{2} BR(z,\da^{4}\om) \cdot \frac{z_{\al}}{|z_{\al}|} + Q_{x}^{2} BR(x,\da^{4}\gamma) \cdot \frac{x_{\al}}{|x_{\al}|}
\end{align*}

$J_5$ will be given later. In $J_1$ and $J_2$ we find 4th order terms in derivatives in $z$ and $x$ so they are fine. In $J_3$ we find inside the integrals
\begin{align}
\label{pacostar1}
\Delta_{\beta} z^{\perp}(\al) \cdot z_{\al}(\al) & = (z(\al) - z(\al - \beta) - \beta z_{\al}(\al))^{\perp} \cdot z_{\al}(\al) \\
\label{pacostar2}
\Delta_{\beta} x^{\perp}(\al) \cdot x_{\al}(\al) & = (x(\al) - x(\al - \beta) - \beta x_{\al}(\al))^{\perp} \cdot x_{\al}(\al)
\end{align}

This implies that we find ''Hilbert'' transforms applied to four derivatives of $x$ and $z$. We are done with $J_3$.

In $J_4$ we also find them inside the integrals \eqref{pacostar1} and \eqref{pacostar2} so it is easy to check that we have kernels whose main singularity is homogenous of degree 0 applied to four derivatives of $\da^{4} \om$ and $\da^{4} \gamma$. This implies that we have a Hilbert transform applied to $\da^{3} \om$ and $\da^{3}\gamma$ so we are done with $J_4$. The most dangerous term is $J_5$ which is given by
\begin{align*}
J_5 & = - Q_{z}^{2} \frac{1}{2\pi}\int_{-\pi}^{\pi} \frac{\Delta_{\beta}\da^{4} z^{\perp}(\al)}{|\Delta_{\beta}z(\al)|^{2}}\cdot \frac{z_{\al}(\al)}{|z_{\al}|}\om(\al - \beta) d\beta
+ Q_{x}^{2} \frac{1}{2\pi}\int_{-\pi}^{\pi} \frac{\Delta_{\beta}\da^{4} x^{\perp}(\al)}{|\Delta_{\beta}x(\al)|^{2}}\cdot \frac{x_{\al}(\al)}{|x_{\al}|}\gamma(\al - \beta) d\beta
\end{align*}

We split further
\begin{align*}
J_5 & = J_{5,1} + J_{5,2} \\
J_{5,1} & = - Q_{z}^{2} \frac{1}{2\pi}\frac{z_{\al}(\al)}{|z_{\al}|}\cdot\int_{-\pi}^{\pi} \left(
\frac{\om(\al - \beta)}{|\Delta_{\beta}z(\al)|^{2}} - \frac{\om(\al)}{|z_{\al}|^{2}4\sin^{2}\left(\beta/2\right)}
\right)\Delta_{\beta}\da^{4} z^{\perp}(\al)d\beta\\
& + Q_{x}^{2} \frac{1}{2\pi}\frac{x_{\al}(\al)}{|x_{\al}|}\cdot \int_{-\pi}^{\pi} \left(
\frac{\gamma(\al - \beta)}{|\Delta_{\beta}x(\al)|^{2}} - \frac{\gamma(\al)}{|x_{\al}|^{2}4\sin^{2}\left(\beta/2\right)}
\right)\Delta_{\beta}\da^{4} x^{\perp}(\al)d\beta \\
J_{5,2} & = - Q_{z}^{2}\frac{1}{2} \frac{z_{\al}(\al)}{|z_{\al}|^{3}} \om(\al)\cdot \Lambda(\da^{4}z^{\perp})(\al)
+ Q_{x}^{2}\frac{1}{2} \frac{x_{\al}(\al)}{|x_{\al}|^{3}} \gamma(\al)\cdot \Lambda(\da^{4}x^{\perp})(\al)
\end{align*}

In $J_{5,1}$ we find a Hilbert transform applied to $\da^{4} z^{\perp}$ and $\da^{4} x^{\perp}$ so it is fine. We split further:

\begin{align*}
J_{5,2} & = J_{5,2,1} + J_{5,2,2} + J_{5,2,3} \\
J_{5,2,1} & =  \left(Q_{x}^{2}\frac{1}{2} \frac{x_{\al}(\al)}{|x_{\al}|^{3}} \gamma(\al) - Q_{z}^{2}\frac{1}{2} \frac{z_{\al}(\al)}{|z_{\al}|^{3}} \om(\al)\right)\cdot \Lambda(\da^{4}x^{\perp})(\al)\\
J_{5,2,2} & =  \Lambda\left(Q_{z}^{2}\frac{1}{2} \frac{z_{\al}(\al)}{|z_{\al}|^{3}} \om(\al)\da^{4} D^{\perp}\right)-
Q_{z}^{2}\frac{1}{2} \frac{z_{\al}(\al)}{|z_{\al}|^{3}} \om(\al)\Lambda(\da^{4} D^{\perp})\\
J_{5,2,3} & = - \Lambda\left(Q_{z}^{2}\frac{1}{2} \frac{z_{\al}(\al)}{|z_{\al}|^{3}} \om(\al)\da^{4} D^{\perp}\right)
\end{align*}

$J_{5,2,1}$ can be estimated as before (there are more derivatives: 5 in total, but they are in $x$). In $J_{5,2,2}$ we find a commutator. Finally:

\begin{align*}
I_{1,4,1}^{4,2} \leq CP(E(t)) - \int_{-\pi}^{\pi}Q_{z}^{2}\sigma_{z} \da^{4}D \cdot \frac{z^{\perp}(\al)}{|z_{\al}|}H\left(\Lambda\left(Q_{z}^{2}\frac{1}{2}\frac{z_{\al}}{|z_{\al}|^3}\om \da^{4} D^{\perp}\right)\right)d\al
\end{align*}

We use that $H(\Lambda) = - \da$ and $z_{\al} \cdot \da^{4} D^{\perp} = - z_{\al}^{\perp} \cdot \da^4 D$ to obtain:

\begin{align*}
I_{1,4,1}^{4,2} & \leq CP(E(t))  - \frac{1}{2}\int_{-\pi}^{\pi}Q_{z}^{2}\sigma_{z} \da^{4}D \cdot \frac{z^{\perp}(\al)}{|z_{\al}|}\da\left(\frac{Q_{z}^{2}\om}{|z_{\al}|^2}\da^{4} D\cdot \frac{z_{\al}^{\perp}}{|z_{\al}|}\right)d\al \\
& \leq CP(E(t))  - \underbrace{\frac{1}{2}\int_{-\pi}^{\pi}Q_{z}^{2}\sigma_{z} \da^{4}D \cdot \frac{z^{\perp}(\al)}{|z_{\al}|}
\da^{4}D \cdot \frac{z^{\perp}(\al)}{|z_{\al}|}\da\left(\frac{Q_{z}^{2} \om}{|z_{\al}|}\right)d\al}_{\text{Easy to estimate by $CP(E(t))$}}\\
&  - \frac{1}{2}\int_{-\pi}^{\pi}Q_{z}^{2}\sigma_{z} \frac{Q_{z}^{2}\om}{|z_{\al}|^2}\underbrace{\da^{4}D \cdot \frac{z^{\perp}(\al)}{|z_{\al}|}
\da\left(\da^{4} D\cdot \frac{z_{\al}^{\perp}}{|z_{\al}|}\right)}_{\text{Integration by parts}}d\al
\end{align*}

Then we are done with $I_{1,4,1}^{4,2}$, $I_{1,4,1}^{4}$, $I_{1,4,1}$, $I_{1,4}$ and $I_1$.

To finish with $I$ it remains to control $I_2$. We split it as:

\begin{align*}
I_2 & = I_{2,1} + I_{2,2} + \text{ l.o.t} \\
I_{2,1} & = \int_{-\pi}^{\pi}\frac{Q_{z}^{2}}{|z_{\al}|^{2}}\sigma_{z}\partial_{\al}^{4}D (c \da^{5}z - b \da^{5} x)d\al\\
I_{2,2} & = \int_{-\pi}^{\pi}\frac{Q_{z}^{2}}{|z_{\al}|^{2}}\sigma_{z}\partial_{\al}^{4}D (\da^{4} c z_{\al} - \da^{4} b x_{\al})d\al\\
\end{align*}

The low order terms are easier to deal with. We further split $I_{2,1}$.
\begin{align*}
I_{2,1} &= I_{2,1,1} + I_{2,1,2} + I_{2,1,3} \\
I_{2,1,1} &= \int_{-\pi}^{\pi}\frac{Q_{z}^{2}}{|z_{\al}|^{2}}\sigma_{z}c \underbrace{\partial_{\al}^{4}D \da^{5}D}_{\text{Integration by parts}} d\al\\
I_{2,1,2} &= \int_{-\pi}^{\pi}\frac{Q_{z}^{2}}{|z_{\al}|^{2}}\sigma_{z}\partial_{\al}^{4}D (c - b_s)\underbrace{\da^{5}x}_{\text{5 derivatives, but in $x$}} d\al\\
I_{2,1,3} &= \underbrace{\int_{-\pi}^{\pi}\frac{Q_{z}^{2}}{|z_{\al}|^{2}}\sigma_{z}\partial_{\al}^{4}D b_e\da^{5}x}_{\text{Error term}} d\al\\
\end{align*}

We find $I_{2,1} \leq CP(E(t)) + c\delta(t)$. We decompose $I_{2,2}$.
\begin{align*}
I_{2,2} &= I_{2,2,1} + I_{2,2,2} + I_{2,2,3} \\
I_{2,2,1} & = \int_{-\pi}^{\pi}\frac{Q_{z}^{2}}{|z_{\al}|^{2}}\sigma_{z}\partial_{\al}^{4}D (\da^{4} c - \da^{4} b_s)\cdot z_{\al}d\al\\
I_{2,2,2} & = \int_{-\pi}^{\pi}\frac{Q_{z}^{2}}{|z_{\al}|^{2}}\sigma_{z}\partial_{\al}^{4}D \underbrace{\da^{4} b_s}_{\text{5 derivatives in $x$}} \da Dd\al\\
I_{2,2,3} & = \underbrace{-\int_{-\pi}^{\pi}\frac{Q_{z}^{2}}{|z_{\al}|^{2}}\sigma_{z}\partial_{\al}^{4}D \da^{4}b_e x_{\al}d\al}_{\text{Error term}}\\
\end{align*}

We deal with $I_{2,2,1}$ more carefully. We use that

\begin{align*}
\da^{4} D \cdot z_{\al} & = \da^{4} z \cdot z_{\al} - \da^{4} x \cdot x_{\al} - \da^{4} x \cdot D_{\al} \\
& = - 3 \da^{3} z \cdot \da^{2} z + 3 \da^{3} x \cdot \da^{2} x - \da^{4} x \cdot D_{\al} \\
& = - 3 \da^{3} D \cdot \da^{2} z - 3 \da^{3} x \cdot \da^{2} D - \da^{4} x \cdot D_{\al}
\end{align*}

to obtain
\begin{align*}
I_{2,2,1} & = I_{2,2,1}^{1} + I_{2,2,1}^{2} + I_{2,2,1}^{3} \\
I_{2,2,1}^{1} & = -3\int_{-\pi}^{\pi}\frac{Q_{z}^{2}}{|z_{\al}|^{2}}\sigma_{z}\partial_{\al}^{3}D \cdot \da^{2} z \da^{4}(c - b_s) d\al\\
I_{2,2,1}^{2} & = -3\int_{-\pi}^{\pi}\frac{Q_{z}^{2}}{|z_{\al}|^{2}}\sigma_{z}\partial_{\al}^{3}x \cdot \da^{2} D \da^{4}(c - b_s) d\al\\
I_{2,2,1}^{3} & = -\int_{-\pi}^{\pi}\frac{Q_{z}^{2}}{|z_{\al}|^{2}}\sigma_{z}\partial_{\al}^{4}x \cdot \da D \da^{4}(c - b_s) d\al\\
\end{align*}

We can integrate by parts in all of the above terms to get low order terms. We are finally done with $I$.

\subsection{Computing the difference $\varphi-\psi$}

From the local existence proof we find the equation for $\vp_t$:

\begin{equation}
\left\{
\begin{array}{cl}
\displaystyle \vp_t & = \displaystyle -\vp B_{z}(t) - \frac{Q_z^2}{2|z_{al}|}\da\left(\frac{\vp^{2}}{Q_{z}^{2}}\right) - Q_{z}^{2}\left(BR_t \cdot \frac{z_{\al}}{|z_{\al}|}+\frac{(P^{-1}_{2}(z))_{\al}}{|z_{\al}|}\right) \\
 &  \displaystyle + Q_{z}Q_{t}^{z}\frac{\om}{|z_{\al}|} - 2cBR \cdot \frac{z_{\al}}{|z_{\al}|} Q_{z}Q_{\al}^{z} - c^2|z_{\al}|\frac{Q_{\al}^{z}}{Q_{z}}
 - \frac{Q_{z}^{3}}{|z_{\al}|}|BR|^{2}Q_{\al}^{z} - (c|z_{\al}|)_{t}\\
 \displaystyle B_z(t) & = \displaystyle \frac{1}{2\pi}\int_{-\pi}^{\pi}(Q_{z}^{2}BR)_{\al} \cdot \frac{z_{\al}}{|z_{\al}|}d\al
\end{array}
\right.
\end{equation}

We will show how to find the equation for $\psi_t$. We start from

$$ \psi = \frac{Q_{x}^{2}\gamma}{2|x_{\al}|} - b_s|x_{\al}|$$

and therefore
$$ \frac{\psi^{2}}{Q_{x}^{2}} = \frac{Q_{x}^{2}\gamma^{2}}{4|x_{\al}|} + \frac{b_s^{2}|x_{\al}|^{2}}{Q_{x}^{2}} - \gamma b_s,$$

that yields

$$ -\da\left(\frac{\psi^{2}}{Q_{x}^{2}}\right) = -\da\left(\frac{Q_{x}^{2}\gamma^{2}}{4|x_{\al}|}\right) - \da\left(\frac{b_s^{2}|x_{\al}|^{2}}{Q_{x}^{2}}\right) + \da\left(\gamma b_s\right)$$

The equation for $\gamma_t$ reads:
\begin{align*}
\gamma_t  = & -2BR_t \cdot x_{\al} - (Q_{x}^{2})_{\al}|BR|^{2} + 2b_sBR_{\al} \cdot x_{\al} \\
& -\da\left(\frac{\psi^{2}}{Q_{x}^{2}}\right) + \left(\frac{b_s^{2}|x_{\al}|^{2}}{Q_{x}^{2}}\right)_{\al} - 2(P^{-1}_{2}(z))_{\al}
+ 2b_e BR_{\al} \cdot x_{\al} + (b_e \gamma)_{\al} + g
\end{align*}

Then
\begin{align*}
\om_t  = & Q_x(Q_x)_{t} \frac{\gamma}{|x_{\al}}| - \frac{Q_x^2 \gamma}{2|x_{\al}|^{3}}x_{\al} \cdot x_{\al t}
+ \frac{Q_x^2 \gamma_{t}}{2|x_{\al}|} - (b_s|x_{\al}|)_{t}\\
= & Q_x(Q_x)_{t} \frac{\gamma}{|x_{\al}|} - \frac{Q_x^2 \gamma}{2|x_{\al}|}B_x(t) - \frac{Q_x^2 \gamma}{2|x_{\al}|}\frac{1}{2\pi}\int_{-\pi}^{\pi}f_{\al}\cdot \frac{x_{\al}}{|x_{\al}|^2}d\al \\
& + \frac{Q_x^2}{2|x_{\al}|}\left(-2BR_t \cdot x_{\al} - (Q_{x}^{2})_{\al}|BR|^{2} + 2b_sBR_{\al} \cdot x_{\al}
 -\da\left(\frac{\psi^{2}}{Q_{x}^{2}}\right) + \left(\frac{b_s^{2}|x_{\al}|^{2}}{Q_{x}^{2}}\right)_{\al}\right.\\
  & \left.- 2(P^{-1}_{2}(z))_{\al}
+ 2b_e BR_{\al} \cdot x_{\al} + (b_e \gamma)_{\al} + g\right)
- (b_s|x_{\al}|)_{t}
\end{align*}

We should remark that we have used that

$$ x_{\al} \cdot x_{\al t} = \frac{1}{2\pi}\int_{-\pi}^{\pi}(Q_{x}^{2} BR)_{\al} \cdot x_{\al} d\al + \frac{1}{2\pi}\int_{-\pi}^{\pi}f_{\al} \cdot x_{\al} d\al$$

and

$$ B_x(t) = \frac{1}{2\pi}\int_{-\pi}^{\pi}(Q_x^{2} BR)_{\al} \cdot \frac{x_{\al}}{|x_{\al}|^{2}}d\al$$

Computing, we find that

\begin{align*}
\psi_t = & Q_x(Q_x)_{t} \frac{\gamma}{|x_{\al}|} - \underbrace{\frac{Q_x^2 \gamma}{2|x_{\al}|}B_x(t)}_{(1)} - Q_x^2 BR_t\frac{x_{\al}}{|x_{\al}|} -
\frac{Q_x^3}{|x_{\al}|}|BR|^{2}Q_{\al}^{x} + Q_{x}^{2}b_sBR_{\al} \cdot \frac{x_{\al}}{|x_{\al}|} \\
& - \frac{Q_{x}^{2}}{2|x_{\al}|}\da\left(\frac{\psi^{2}}{Q_{x}^{2}}\right) + \underbrace{\frac{Q_{x}^{2}}{2}\left(\frac{b_s^{2}|x_{\al}|^{2}}{Q_{x}^{2}}\right)_{\al}}_{(1)}
- \frac{Q_{x}^{2}}{|x_{\al}|}(P^{-1}_{2}(z))_{\al} - (b_s|x_{\al}|)_{t} + \mathcal{E}^{1}
\end{align*}

where
$$ \mathcal{E}^{1} =  \frac{Q_{x}^{2}}{|x_{\al}|}BR_{\al}\cdot x_{\al} b_e + \frac{Q_{x}^{2}}{2|x_{\al}|}(b_e \gamma)_{\al} + \frac{Q_{x}^{2}}{|x_{\al}|} g
- \frac{Q_x^2 \gamma}{2|x_{\al}|}\frac{1}{2\pi}\int_{-\pi}^{\pi}f_{\al}\cdot \frac{x_{\al}}{|x_{\al}|^2}d\al$$

are  error terms. We consider
\begin{align*}
(1) = &-\frac{Q_x^2 \gamma}{2|x_{\al}|}B_x(t) +\frac{Q_{x}^{2}}{2}\left(\frac{b_s^{2}|x_{\al}|^{2}}{Q_{x}^{2}}\right)_{\al} \\
= & -\frac{Q_x^2 \gamma}{2|x_{\al}|}B_x(t) +b_s|x_{\al}|(b_s)_{\al} - \frac{Q_{\al}^{x}}{Q_{x}}b_s^{2}|x_{\al}| \\
= & -\frac{Q_x^2 \gamma}{2|x_{\al}|}B_x(t) +b_s|x_{\al}|B_x(t)-b_s(Q_{x}^{2}BR)_{\al} \cdot \frac{x_{\al}}{|x_{\al}|} - \frac{Q_{\al}^{x}}{Q_{x}}b_s^{2}|x_{\al}| \\
= & -B_x(t)\psi-b_s(Q_{x}^{2}BR)_{\al} \cdot \frac{x_{\al}}{|x_{\al}|} - \frac{Q_{\al}^{x}}{Q_{x}}b_s^{2}|x_{\al}|
\end{align*}
It yields

\begin{align*}
\psi_t = & Q_x(Q_x)_{t} \frac{\gamma}{|x_{\al}|} - B_x(t) \psi \underbrace{- b_s(Q_{x}^{2}BR)_{\al} \cdot \frac{x_{\al}}{|x_{\al}|}}_{(2)} - \frac{Q_{\al}^{x}}{Q_{x}}b_s^{2}|x_{\al}| \\
& - Q_x^2 BR_t\frac{x_{\al}}{|x_{\al}|} -
\frac{Q_x^3}{|x_{\al}|}|BR|^{2}Q_{\al}^{x} + \underbrace{Q_{x}^{2}b_sBR_{\al} \cdot \frac{x_{\al}}{|x_{\al}|}}_{(2)} \\
& - \frac{Q_{x}^{2}}{2|x_{\al}|}\da\left(\frac{\psi^{2}}{Q_{x}^{2}}\right)
- \frac{Q_{x}^{2}}{|x_{\al}|}(P^{-1}_{2}(z))_{\al} - (b_s|x_{\al}|)_{t} + \mathcal{E}^{1}
\end{align*}

It is easy to check that
\begin{align*}
(2) = & - b_s(Q_{x}^{2}BR)_{\al} \cdot \frac{x_{\al}}{|x_{\al}|} + Q_{x}^{2}b_sBR_{\al} \cdot \frac{x_{\al}}{|x_{\al}|}
= -2b_s BR \cdot \frac{x_{\al}}{|x_{\al}|} Q_x (Q_x)_{\al},
\end{align*}

then

\begin{align*}
\psi_t = & -B_x(t) \psi - \frac{Q_{x}^{2}}{2|x_{\al}|}\da\left(\frac{\psi^{2}}{Q_{x}^{2}}\right) - Q_x^2\left(BR_t\frac{x_{\al}}{|x_{\al}|}
+ \frac{(P^{-1}_{2}(z))_{\al}}{|x_{\al}|} \right) \\
& + Q_x(Q_x)_{t} \frac{\gamma}{|x_{\al}|}-2b_s BR \cdot \frac{x_{\al}}{|x_{\al}|} Q_x (Q_x)_{\al} - \frac{Q_{\al}^{x}}{Q_{x}}b_s^{2}|x_{\al}|
-\frac{Q_x^3}{|x_{\al}|}|BR|^{2}Q_{\al}^{x} \\
& - (b_s|x_{\al}|)_{t} + \mathcal{E}^{1}
\end{align*}

With this formula it is easy to find that

$$ \frac{1}{2}\frac{d}{dt}\int |\mathcal{D}|^{2}dx \leq CP(E(t)) + c\delta(t)$$

In order to deal with II
$$ II = \int_{-\pi}^{\pi}\Lambda \da^{3} \mathcal{D} \da^{3}\mathcal{D}_{t} d\al$$

we take a derivative in $\al$ in the equation for $\om$ and $\psi$ to reorganize the most dangerous terms. If we find a term of low order, we will denote it by NICE. Since the equations for $\varphi_t$ and $\psi_t$ are analogous except for the $\mathcal{E}^{1}$ term, the NICE terms are going to be easier to estimate in terms of $CP(E(t)) + c\delta(t)$.

\begin{align*}
\psi_{\al t} = & -B_x(t) \psi_{\al} - \da\left(\frac{Q_{x}^{2}}{2|x_{\al}|}\da\left(\frac{\psi^{2}}{Q_{x}^{2}}\right)\right) - \left(Q_x^2\left(\underbrace{BR_t\frac{x_{\al}}{|x_{\al}|}}_{(3)}
+ \frac{(P^{-1}_{2}(z))_{\al}}{|x_{\al}|} \right)\right)_{\al} \\
& + \left(Q_x(Q_x)_{t} \frac{\gamma}{|x_{\al}}|\right)_{\al} - \left(2b_s BR \cdot \frac{x_{\al}}{|x_{\al}|} Q_x (Q_x)_{\al}\right)_{\al} - \left(\frac{Q_{\al}^{x}}{Q_{x}}b_s^{2}|x_{\al}|\right)_{\al}
-\left(\frac{Q_x^3}{|x_{\al}|}|BR|^{2}Q_{\al}^{x}\right)_{\al} \\
& \underbrace{- (b_s|x_{\al}|)_{\al t}}_{(3)} + \mathcal{E}^{1}_{\al}
\end{align*}

Expanding (3):
\begin{align*}
(3) = & -\left(Q_{x}^{2}BR_t\frac{x_{\al}}{|x_{\al}|}\right)_{\al} - (b_s|x_{\al}|)_{\al t} \\
= & -\left(Q_{x}^{2}BR_t\right)_{\al}\frac{x_{\al}}{|x_{\al}|}-Q_{x}^{2}BR_t\left(\frac{x_{\al}}{|x_{\al}|}\right)_{\al}
- \left(|x_{\al}|B_{x}(t) - (Q_{x}^{2} BR)_{\al} \cdot \frac{x_{\al}}{|x_{\al}|}\right)_{t} \\
= & - \left(|x_{\al}|B_{x}(t)\right)_{t} +   (Q_{x}^{2} BR)_{\al} \cdot \left(\frac{x_{\al}}{|x_{\al}|}\right)_{t}
+ 2(Q_{x}(Q_{x})_{t} BR)_{\al} \cdot\frac{x_{\al}}{|x_{\al}|} - Q_{x}^{2}BR_{t} \cdot  \left(\frac{x_{\al}}{|x_{\al}|}\right)_{\al}
\end{align*}
We use that
$$  \left(\frac{x_{\al}}{|x_{\al}|}\right)_{\al} = \frac{x_{\al \al} \cdot x_{\al}^{\perp}}{|x_{\al}|^{2}} \frac{x_{\al}^{\perp}}{|x_{\al}|};
\quad \left(\frac{x_{\al}}{|x_{\al}|}\right)_{t} = \frac{x_{\al t} \cdot x_{\al}^{\perp}}{|x_{\al}|^{2}} \frac{x_{\al}^{\perp}}{|x_{\al}|}$$
to find

\begin{align*}
\psi_{\al t} = & \underbrace{-B_x(t) \psi_{\al}}_{(4)} - \underbrace{\frac{\da^{2}(\psi^{2})}{2|x_{\al}|}}_{(5)}
+ \underbrace{\da\left(\frac{(Q_{x})_{\al}}{|x_{\al}|Q_{x}}\psi^{2}\right)}_{(6)}
- Q_{x}^{2} BR_t \cdot x_{\al}^{\perp} \frac{x_{\al \al} \cdot x_{\al}^{\perp}}{|x_{\al}|^{3}} -
(|x_{\al}|B_x(t))_{t} \\
& + \underbrace{(Q_{x}^{2} BR)_{\al} \cdot x_{\al}^{\perp} \frac{x_{\al t} \cdot x_{\al}^{\perp}}{|x_{\al}|^{3}}}_{(13)}
+ \underbrace{2(Q_x (Q_x)_{t} BR)_{\al} \frac{x_{\al}}{|x_{\al}|}}_{(7)}
- \underbrace{\left(Q_{x}^{2}\frac{(P^{-1}_{2}(z))_{\al}}{|x_{\al}|}\right)_{\al}}_{(8)}\\
&+ \underbrace{\left(Q_x(Q_x)_{t} \frac{\gamma}{|x_{\al}|}\right)_{\al}}_{(9)}
- \underbrace{\left(2b_s BR \cdot \frac{x_{\al}}{|x_{\al}|}Q_x (Q_x)_{\al}\right)_{\al}}_{(10)}
- \underbrace{\left(\frac{(Q_x)_{\al}}{Q_x}b_s^{2}|x_{\al}|\right)_{\al}}_{(11)}\\
&- \underbrace{\left(\frac{Q_x^{3}}{|x_{\al}|}|BR|^{2}(Q_x)_{\al}\right)_{\al}}_{(12)} + \mathcal{E}^{1}_{\al}
\end{align*}

The term $(|x_{\al}|B_x(t))_{t}$ depends only on $t$ so it is not going to appear in computing II.

$$ (4) = -B_x(t) \psi_{\al} \text{ is NICE (at the level of } \psi_{\al})$$

$$ (5) = -\frac{\da^{2}(\psi^{2})}{2|x_{\al}|} \text{ is a transparent term which is NICE (even if we have to deal with } \Lambda^{1/2})$$

$$ (6) = \da\left(\frac{(Q_{x})_{\al}}{|x_{\al}|Q_{x}}\psi^{2}\right) = -\frac{(Q_x)^{2}_{\al}}{|x_{\al}|(Q_x)^{2}}
+ \frac{2(Q_x)_{\al}\psi \psi_{\al}}{|x_{\al}|Q_x} + \frac{\psi^{2}}{Q_x}\left(\frac{(Q_x)_{\al}}{|x_{\al}|}\right)_{\al}$$

The first term is at the level of $\da x$ so it is NICE. The second term is at the level of $\da x$ or $\psi_{\al}$ so it is NICE. We write the last one as
$$\frac{\psi^{2}}{Q_x}\left(\frac{(Q_x)_{\al}}{|x_{\al}|}\right)_{\al}
= \frac{\psi^{2}}{Q_x}x_{\al} \cdot \left(\nabla^{2}Q(x)\cdot \frac{x_{\al}}{|x_{\al}|}\right)
+ \frac{\psi^{2}}{Q_x}x_{\al} \nabla Q(x) \cdot x_{\al}^{\perp} \frac{x_{\al \al} \cdot x_{\al}^{\perp}}{|x_{\al}|^{3}}$$

The first term is at the level of $x_{\al}$ or $\psi$ so it is NICE. For the second term we have used that

$$  \left(\frac{x_{\al}}{|x_{\al}|}\right)_{\al} = \frac{x_{\al \al} \cdot x_{\al}^{\perp}}{|x_{\al}|^{2}} \frac{x_{\al}^{\perp}}{|x_{\al}|}$$
Finally:

$$ (6) = \text{NICE } + \frac{\psi^{2}}{Q_x}x_{\al} \nabla Q(x) \cdot x_{\al}^{\perp} \frac{x_{\al \al} \cdot x_{\al}^{\perp}}{|x_{\al}|^{3}}$$

\begin{align*} (7) = 2(Q_x (Q_x)_{t} BR)_{\al} \frac{x_{\al}}{|x_{\al}|}
&= 2(Q_x)_{\al} (Q_x)_{t} BR \cdot  \frac{x_{\al}}{|x_{\al}|}
+ 2Q_x \left(\frac{(Q_x)_{t}}{|x_{\al}}\right)_{\al} BR \cdot x_{\al}\\
&+ 2Q_x(Q_x)_{t} BR_{\al} \cdot \frac{x_{\al}}{|x_{\al}|}\end{align*}

The first term is at the level of $x_{\al}, x_{t}, BR \sim x_{\al}$ so it is NICE. We use that

$$ \frac{(Q_x)_{t\al}}{|x_{\al}|} = \frac{(Q_x)_{\al t}}{|x_{\al}|}
= \frac{(\nabla Q(x) \cdot x_{\al})_{t}}{|x_{\al}|}
= \left(\nabla Q(x) \cdot \frac{x_{\al}}{|x_{\al}|}\right)_{t} - \nabla Q(x) \cdot x_{\al} \left(\frac{1}{|x_{\al}|}\right)_{t}$$

Using that

$$ \frac{x_{\al} \cdot x_{\al t}}{|x_{\al}|^{2}} = B_{x}(t) + \frac{1}{2\pi}\int_{-\pi}^{\pi}f_{\al} \cdot \frac{x_{\al}}{|x_{\al}|^2}d\al$$
and
$$  \left(\frac{x_{\al}}{|x_{\al}|}\right)_{t} = \frac{x_{\al t} \cdot x_{\al}^{\perp}}{|x_{\al}|^{2}} \cdot\frac{x_{\al}^{\perp}}{|x_{\al}|}$$
we find that
\begin{align}
\label{star}
\frac{(Q_x)_{t\al}}{|x_{\al}|}
&= x_t \cdot \left(\nabla^2 Q(x) \cdot \frac{x_{\al}}{|x_{\al}|}\right) + \nabla Q(x) \cdot x_{\al}^{\perp} \frac{x_{\al t} \cdot x_{\al}^{\perp}}{|x_{\al}|^{3}}\nonumber\\
&+ \nabla Q(x) \cdot\frac{x_{\al}}{|x_{\al}|} B_x(t) + \nabla Q(x) \cdot\frac{x_{\al}}{|x_{\al}|}\frac{1}{2\pi}\int_{-\pi}^{\pi}f_{\al} \cdot \frac{x_{\al}}{|x_{\al}|^2}d\al
\end{align}

That yields

\begin{align*}
(7) & = 2(Q_x (Q_x)_{t} BR)_{\al} \frac{x_{\al}}{|x_{\al}|}
= \text{NICE } + \underbrace{2Q_x BR \cdot x_{\al} x_t \cdot \left(\nabla^2 Q(x) \cdot \frac{x_{\al}}{|x_{\al}|}\right)}_{\text{NICE (at the level of }x_{\al}, x_t, BR)} \\
& + 2Q_x BR \cdot x_{\al} \nabla Q(x) \cdot x_{\al}^{\perp} \frac{x_{\al t} \cdot x_{\al}^{\perp}}{|x_{\al}|^{3}}
+ \underbrace{2Q_x BR \cdot x_{\al} \nabla Q(x) \cdot\frac{x_{\al}}{|x_{\al}|} B_x(t)}_{\text{NICE (at the level of }x_{\al}, x_t, BR)} \\
& + \underbrace{2Q_x BR \cdot x_{\al} \nabla Q(x) \cdot\frac{x_{\al}}{|x_{\al}|}\frac{1}{2\pi}\int_{-\pi}^{\pi}f_{\al} \cdot \frac{x_{\al}}{|x_{\al}|^2}d\al}_{\text{part of error terms}}
+ 2Q_x(Q_x)_{t} BR_{\al} \cdot \frac{x_{\al}}{|x_{\al}|}
\end{align*}

Finally:

\begin{align*}
(7) & = 2(Q_x (Q_x)_{t} BR)_{\al} \frac{x_{\al}}{|x_{\al}|}
= \text{NICE } + 2Q_x BR \cdot x_{\al} \nabla Q(x) \cdot x_{\al}^{\perp} \frac{x_{\al t} \cdot x_{\al}^{\perp}}{|x_{\al}|^{3}}
+ 2Q_x(Q_x)_{t} BR_{\al} \cdot \frac{x_{\al}}{|x_{\al}|}
\end{align*}

\begin{align*}
(8) & = -\left(Q_{x}^{2}\frac{(P^{-1}_{2}(z))_{\al}}{|x_{\al}|}\right)_{\al}
= -\left(Q_{x}^{2}\nabla P^{-1}_{2}(x)\cdot \frac{x_{\al}}{|x_{\al}|}\right)_{\al}
= \underbrace{-2 Q_{x} \nabla Q^{x}_{x} \cdot x_{\al} \nabla P^{-1}_{2}(x) \cdot \frac{x_{\al}}{|x_{\al}|}}_{\text{NICE (at the level of }x_{\al})} \\
& \underbrace{- Q_{x}^{2} x_{\al} \cdot \left(\nabla^{2} P^{-1}_{2}(x) \cdot \frac{x_{\al}}{|x_{\al}|}\right)}_{\text{NICE (at the level of }x_{\al})}
- Q_{x}^{2} \nabla P^{-1}_{2}(x) \cdot x_{\al}^{\perp} \frac{x_{\al \al} \cdot x_{\al}^{\perp}}{|x_{\al}|^{3}}
\end{align*}

which means

\begin{align*}
(8) & = -\left(Q_{x}^{2}\frac{(P^{-1}_{2}(z))_{\al}}{|x_{\al}|}\right)_{\al}
= \text{NICE } - Q_{x}^{2} \nabla P^{-1}_{2}(x) \cdot x_{\al}^{\perp} \frac{x_{\al \al} \cdot x_{\al}^{\perp}}{|x_{\al}|^{3}}
\end{align*}

\begin{align*}
(9) & = \left(Q_x(Q_x)_{t} \frac{\gamma}{|x_{\al}|}\right)_{\al}
= \underbrace{(Q_x)_{\al} (Q_x)_{t}\frac{\gamma}{|x_{\al}|}}_{\text{NICE (at the level of }x_{\al}, x_t)} + Q_x \frac{(Q_x)_{\al t}}{|x_{\al}|}\gamma
+ Q_x(Q_x)_{t} \left(\frac{\gamma}{|x_{\al}|}\right)_{\al}
\end{align*}

We use \eqref{star} to deal with $\frac{(Q_x)_{\al t}}{|x_{\al}|}$. We find that

\begin{align*}
(9) & = \left(Q_x(Q_x)_{t} \frac{\gamma}{|x_{\al}|}\right)_{\al}
= \text{NICE } + Q_x \gamma \nabla Q(x) \cdot x_{\al}^{\perp} \frac{x_{\al t} \cdot x_{\al}^{\perp}}{|x_{\al}|^{3}}
+ Q_x(Q_x)_{t} \left(\frac{\gamma}{|x_{\al}|}\right)_{\al}
\end{align*}

\begin{align*}
(10) & = - \left(2b_s BR \cdot \frac{x_{\al}}{|x_{\al}|}Q_x (Q_x)_{\al}\right)_{\al}
= \underbrace{- 2b_s BR \cdot \frac{x_{\al}}{|x_{\al}|}(Q_x)^{2}_{\al}}_{\text{NICE as before}}
- \left(2b_s BR \cdot \frac{x_{\al}}{|x_{\al}|}\right)_{\al}Q_x (Q_x)_{\al} \\
& - 2b_s BR \cdot x_{\al}Q_x \nabla Q_x(x) \cdot x_{\al}^{\perp} \frac{x_{\al \al} \cdot x_{\al}^{\perp}}{|x_{\al}|^{3}}
 \underbrace{- 2b_s BR \cdot \frac{x_{\al}}{|x_{\al}|}Q_x x_{\al} (\nabla^{2} Q_{x}(x)) \cdot x_{\al}}_{\text{NICE as before}}
\end{align*}
Therefore

\begin{align*}
(10)  = - \left(2b_s BR \cdot \frac{x_{\al}}{|x_{\al}|}Q_x (Q_x)_{\al}\right)_{\al}
& = \text{NICE }
- \left(2b_s BR \cdot \frac{x_{\al}}{|x_{\al}|}\right)_{\al}Q_x (Q_x)_{\al} \\
& - 2b_s BR \cdot x_{\al}Q_x \nabla Q_x(x) \cdot x_{\al}^{\perp} \frac{x_{\al \al} \cdot x_{\al}^{\perp}}{|x_{\al}|^{3}}
\end{align*}

\begin{align*}
(11)  = - \left(\frac{(Q_x)_{\al}}{Q_x}b_s^{2}|x_{\al}|\right)_{\al}
& = - \left(b_s^{2}|x_{\al}|\right)_{\al}\frac{(Q_x)_{\al}}{Q_x}
- \frac{b_s^{2}|x_{\al}|^2}{Q_x}\nabla Q(x) \cdot x_{\al}^{\perp} \frac{x_{\al \al} \cdot x_{\al}^{\perp}}{|x_{\al}|^{3}} \\
& - \frac{x_{\al}(\nabla^{2} Q(x) \cdot x_{\al})}{Q_x}b_s^2|x_{\al}| + \frac{(Q_x)_{\al}^2}{(Q_x)^2}b_s^2|x_{\al}|
\end{align*}

The fact that the last two terms are NICE, allows us to find that

\begin{align*}
(11)  = - \left(\frac{(Q_x)_{\al}}{Q_x}b_s^{2}|x_{\al}|\right)_{\al}
& = \text{NICE } - \left(b_s^{2}|x_{\al}|\right)_{\al}\frac{(Q_x)_{\al}}{Q_x}
- \frac{b_s^{2}|x_{\al}|^2}{Q_x}\nabla Q(x) \cdot x_{\al}^{\perp} \frac{x_{\al \al} \cdot x_{\al}^{\perp}}{|x_{\al}|^{3}}
\end{align*}

Finally:
\begin{align*}
(12) = - \left(\frac{Q_x^{3}}{|x_{\al}|}|BR|^{2}(Q_x)_{\al}\right)_{\al}
& = \underbrace{- 3(Q_x)^2(Q_x)^2_{\al}|BR|^{2}}_{\text{NICE}} - \frac{Q_x^3}{|x_{\al}|}(|BR|^2)_{\al}(Q_x)_{\al} \\
& \underbrace{- \frac{Q_x^3}{|x_{\al}|}|BR|^2x_{\al} \cdot(\nabla^{2}Q(x) \cdot x_{\al})}_{\text{NICE}} - Q_{x}^{3}|BR|^{2}\nabla Q(x) \cdot x_{\al}^{\perp} \frac{x_{\al \al} \cdot x_{\al}^{\perp}}{|x_{\al}|^{3}}
\end{align*}
which implies that

\begin{align*}
(12) = - \left(\frac{Q_x^{3}}{|x_{\al}|}|BR|^{2}(Q_x)_{\al}\right)_{\al}
& = \text{NICE } - \frac{Q_x^3}{|x_{\al}|}(|BR|^2)_{\al}(Q_x)_{\al} - Q_{x}^{3}|BR|^{2}\nabla Q(x) \cdot x_{\al}^{\perp} \frac{x_{\al \al} \cdot x_{\al}^{\perp}}{|x_{\al}|^{3}}
\end{align*}

We gather all the formulas from (4) to (12) absorbing the error terms  by $\tilde{\mathcal{E}}^{1}_{\alpha}$ whenever we encounter them.

It yields:
\begin{align*}
\psi_{\al t} & = \text{NICE } + \underbrace{\frac{\psi^{2}}{Q_x} \nabla Q(x) \cdot x_{\al}^{\perp} \frac{x_{\al \al} \cdot x_{\al}^{\perp}}{|x_{\al}|^{3}}}_{(16)} \underbrace{- Q_{x}^{2} BR_t \cdot x_{\al}^{\perp} \frac{x_{\al \al} \cdot x_{\al}^{\perp}}{|x_{\al}|^{3}}}_{(15)}
\underbrace{- Q_{x}^{2} \nabla P_{2}^{-1}(x) \cdot x_{\al}^{\perp} \frac{x_{\al \al} \cdot x_{\al}^{\perp}}{|x_{\al}|^{3}}}_{(15)} \\
& + \underbrace{Q_{x} \gamma \nabla Q(x) \cdot x_{\al}^{\perp} \frac{x_{\al t} \cdot x_{\al}^{\perp}}{|x_{\al}|^{3}}}_{(18)}
+ \underbrace{Q_x(Q_x)_{t} \left(\frac{\gamma}{|x_{\al}|}\right)_{\al}}_{(14)}
+ \underbrace{2Q_x BR \cdot x_{\al} \nabla Q(x) \cdot x_{\al}^{\perp} \frac{x_{\al t} \cdot x_{\al}^{\perp}}{|x_{\al}|^{3}}}_{(18)} \\
& + \underbrace{Q_x(Q_x)_{t} 2 BR_{\al} \cdot \frac{x_{\al}}{|x_{\al}|}}_{(14)}
\underbrace{- \left(2 b_x BR \cdot \frac{x_{\al}}{|x_{\al}|}\right)_{\al}Q_x(Q_x)_{\al}}_{(17)}
\underbrace{-2b_s BR \cdot x_{\al}Q_x \nabla Q(x) \cdot x_{\al}^{\perp} \frac{x_{\al \al} \cdot x_{\al}^{\perp}}{|x_{\al}|^{3}}}_{(16)}\\
& \underbrace{-(b_s^{2}|x_{\al}|)_{\al}\frac{(Q_x)_{\al}}{Q_x}}_{(17)}
 \underbrace{-\frac{b_s^2|x_{\al}|^{2}}{Q_x}\nabla Q(x) \cdot x_{\al}^{\perp} \frac{x_{\al \al} \cdot x_{\al}^{\perp}}{|x_{\al}|^{3}}}_{(16)}
 \underbrace{-\frac{Q_x^3}{|x_{\al}|}(|BR|^{2})_{\al}(Q_x)_{\al}}_{(17)}
 \\
& \underbrace{-Q_x^3|BR|^{2}\nabla Q(x) \cdot x_{\al}^{\perp} \frac{x_{\al \al} \cdot x_{\al}^{\perp}}{|x_{\al}|^{3}}}_{(16)}+ (Q_x^2 BR)_{\al} \cdot x_{\al}^{\perp} \frac{x_{\al t} \cdot x_{\al}^{\perp}}{|x_{\al}|^{3}} + \tilde{\mathcal{E}}_{\al}^{1}
\end{align*}

We compute
\begin{align*}
(14) & = Q_x(Q_x)_{t} \left(\frac{\gamma}{|x_{\al}|}\right)_{\al}
+ Q_x(Q_x)_{t} 2 BR_{\al} \cdot \frac{x_{\al}}{|x_{\al}|}\\
&= 2 \frac{(Q_x)_{t}}{Q_x}(Q_x)^2\left(\frac{\gamma}{2|x_{\al}|}\right)_{\al}
+ 2\frac{(Q_x)_{t}}{Q_x}(Q_x)^2 BR_{\al} \cdot \frac{x_{\al}}{|x_{\al}|} \\
& = 2\frac{(Q_x)_{t}}{Q_x}\psi_{\al} - 2\frac{(Q_x)_{t}}{Q_x}(Q_x^{2})_{\al}\frac{\gamma}{2|x_{\al}|}
- 2\frac{(Q_x)_{t}}{Q_x}(Q_x^{2})_{\al}BR_{\al} \frac{x_{\al}}{|x_{\al}|}
- 2\frac{(Q_x)_{t}}{Q_x}(|x_{\al}|B_x(t))
\end{align*}

The last formula allows us to conclude that (14)=NICE. We reorganize using (15), (16), (17) and (18).

\begin{align*}
\psi_{\al t} & = \text{NICE } - Q_x^2(BR_t \cdot x_{\al}^{\perp}
+ \nabla P_{2}^{-1}(x) \cdot x_{\al}^{\perp})\frac{x_{\al \al} \cdot x_{\al}^{\perp}}{|x_{\al}|^{3}} \\
& - Q^3\left(|BR|^{2} + \frac{b_s^2|x_{\al}|^{2}}{Q_x^{4}}+2b_s \frac{BR \cdot x_{\al}}{Q_x^{2}} - \frac{\psi^{2}}{Q_x^4}\right)\nabla Q(x) \cdot x_{\al}^{\perp} \frac{x_{\al \al} \cdot x_{\al}^{\perp}}{|x_{\al}|^{3}} \\
& + (Q_x^{2} BR)_{\al} \cdot x_{\al}^{\perp} \frac{x_{\al t} \cdot x_{\al}^{\perp}}{|x_{\al}|^{3}}
+ (Q_x \gamma + 2Q_x BR \cdot x_{\al})\nabla Q(x) \cdot x_{\al}^{\perp} \frac{x_{\al t} \cdot x_{\al}^{\perp}}{|x_{\al}|^{3}} \\
& - \left(\frac{Q_x^{3}(|BR|^{2})_{\al}}{|x_{\al}|} + \frac{(b_s^{2}|x_{\al}|)_{\al}}{Q_x}
+ \left(2b_s BR \cdot \frac{x_{\al}}{|x_{\al}|}\right)_{\al}Q_x\right)(Q_x)_{\al} + \tilde{\mathcal{E}}_{\al}^{1}
\end{align*}

We add and subtract terms in order to find the R-T condition. We remember here that
\begin{align}
\si_z & = \left(BR_t + \frac{\vp}{|z_{\al}|}BR_{\al}\right) \cdot z_{\al}^{\perp} + \frac{\om}{2|z_{\al}|^{2}}\left(z_{\al t} + \frac{\vp}{|z_{\al}|}z_{\al \al}\right)\cdot z_{\al}^{\perp}\nonumber\\
& + Q_z\left|BR + \frac{\om}{2|z_{\al}|^{2}}z_{\al}\right|^{2}\nabla Q(z) \cdot z_{\al}^{\perp} + \nabla P_{2}^{-1}(z) \cdot z_{\al}^{\perp}\nonumber\\
\si_x & = \left(BR_t + \frac{\psi}{|x_{\al}|}BR_{\al}\right) \cdot x_{\al}^{\perp} + \frac{\gamma}{2|x_{\al}|^{2}}\left(x_{\al t} + \frac{\psi}{|x_{\al}|}x_{\al \al}\right)\cdot x_{\al}^{\perp}\nonumber\\
& \label{X} + Q_x\left|BR + \frac{\gamma}{2|x_{\al}|^{2}}x_{\al}\right|^{2}\nabla Q(x) \cdot x_{\al}^{\perp} + \nabla P_{2}^{-1}(x) \cdot x_{\al}^{\perp}
\end{align}

In $\si_{x}$ there are error terms but they are not dangerous. Then, we find

\begin{align*}
\psi_{\al t} & = \text{NICE } \\
&- Q_x^{2}\left(\left(BR_t + \frac{\psi}{|x_{\al}|}BR_{\al}\right) \cdot x_{\al}^{\perp}
+ \frac{\gamma}{2|x_{\al}|^{2}}\left(x_{\al t} + \frac{\psi}{|x_{\al}|}x_{\al \al}\right)\cdot x_{\al}^{\perp}
+ \nabla P_{2}^{-1}(x) \cdot x_{\al}^{\perp}\right)\frac{x_{\al \al} \cdot x_{\al}^{\perp}}{|x_{\al}|^{3}} \\
& \underbrace{+ (Q_x^{2} BR)_{\al} \cdot x_{\al}^{\perp} \frac{x_{\al t} \cdot x_{\al}^{\perp}}{|x_{\al}|^{3}}
+ Q_x^{2}\left(\frac{\psi}{|x_{\al}|}BR_{\al} \cdot x_{\al}^{\perp} + \frac{\gamma}{2|x_{\al}|^{2}}\left(x_{\al t} + \frac{\psi}{|x_{\al}|}x_{\al \al}\right) \cdot x_{\al}^{\perp}\right)\frac{x_{\al \al} \cdot x_{\al}^{\perp}}{|x_{\al}|^{3}}}_{(19)}\\
& - Q^3\left(|BR|^{2} + \frac{b_s^2|x_{\al}|^{2}}{Q_x^{4}}+2b_s \frac{BR \cdot x_{\al}}{Q_x^{2}} - \frac{\psi^{2}}{Q_x^4}\right)\nabla Q(x) \cdot x_{\al}^{\perp} \frac{x_{\al \al} \cdot x_{\al}^{\perp}}{|x_{\al}|^{3}} \\
& + (Q_x \gamma + 2Q_x BR \cdot x_{\al})\nabla Q(x) \cdot x_{\al}^{\perp} \frac{x_{\al t} \cdot x_{\al}^{\perp}}{|x_{\al}|^{3}} \\
& - \left(\frac{Q_x^{3}(|BR|^{2})_{\al}}{|x_{\al}|} + \frac{(b_s^{2}|x_{\al}|)_{\al}}{Q_x}
+ \left(2b_s BR \cdot \frac{x_{\al}}{|x_{\al}|}\right)_{\al}Q_x\right)(Q_x)_{\al} + \tilde{\mathcal{E}}_{\al}^{1}
\end{align*}

Line (19) can be written as

\begin{align*}
(19) & = (Q_x^{2} BR)_{\al} \cdot x_{\al}^{\perp} \frac{x_{\al t} \cdot x_{\al}^{\perp}}{|x_{\al}|^{3}}
+ Q_x^{2}BR_{\al} \cdot x_{\al}^{\perp}\frac{\psi}{|x_{\al}|}\frac{x_{\al \al} \cdot x_{\al}^{\perp}}{|x_{\al}|^{3}}\\
&+ \frac{Q_x^{2}\gamma}{2|x_{\al}|^{2}}\left(x_{\al t} + \frac{\psi}{|x_{\al}|}x_{\al \al}\right)
 \cdot x_{\al}^{\perp}\frac{x_{\al \al} \cdot x_{\al}^{\perp}}{|x_{\al}|^{3}} \\
 & = (Q_x^{2} BR)_{\al} \cdot x_{\al}^{\perp} \frac{x_{\al t} \cdot x_{\al}^{\perp}}{|x_{\al}|^{3}}
+ (Q_x^{2}BR)_{\al} \cdot x_{\al}^{\perp}\frac{\psi}{|x_{\al}|}\frac{x_{\al \al} \cdot x_{\al}^{\perp}}{|x_{\al}|^{3}}\\
&+ \frac{Q_x^{2}\gamma}{2|x_{\al}|^{2}}\left(x_{\al t}\cdot x_{\al}^{\perp} + \frac{\psi}{|x_{\al}|}x_{\al \al}\cdot x_{\al}^{\perp}\right)
\frac{x_{\al \al} \cdot x_{\al}^{\perp}}{|x_{\al}|^{3}}  - 2Q_x(Q_x)_{\al} BR \cdot x_{\al}^{\perp} \frac{\psi}{|x_{\al}|}\frac{x_{\al \al} \cdot x_{\al}^{\perp}}{|x_{\al}|^{3}} \\
& = (Q_x^{2} BR)_{\al} \cdot x_{\al}^{\perp}\frac{1}{|x_{\al}|^{3}}\left( x_{\al t} \cdot x_{\al}^{\perp}
+ \frac{\psi}{|x_{\al}|}x_{\al \al} \cdot x_{\al}^{\perp}\right)\\
&+ \frac{Q_x^{2}\gamma}{2|x_{\al}|^{2}}\frac{1}{|x_{\al}|^{3}}\left(x_{\al t}\cdot x_{\al}^{\perp} + \frac{\psi}{|x_{\al}|}x_{\al \al}\cdot x_{\al}^{\perp}\right)
x_{\al \al} \cdot x_{\al}^{\perp} - 2Q_x(Q_x)_{\al} BR \cdot x_{\al}^{\perp} \frac{\psi}{|x_{\al}|}\frac{x_{\al \al} \cdot x_{\al}^{\perp}}{|x_{\al}|^{3}} \\
& = \frac{1}{|x_{\al}|^{3}}\left( x_{\al t} \cdot x_{\al}^{\perp}
+ \frac{\psi}{|x_{\al}|}x_{\al \al} \cdot x_{\al}^{\perp}\right)\left((Q_x^{2} BR)_{\al} \cdot x_{\al}^{\perp}
+ \frac{Q_x^{2}\gamma}{2|x_{\al}|^{2}}x_{\al \al} \cdot x_{\al}^{\perp}\right)\\
&- 2Q_x(Q_x)_{\al} BR \cdot x_{\al}^{\perp} \frac{\psi}{|x_{\al}|}\frac{x_{\al \al} \cdot x_{\al}^{\perp}}{|x_{\al}|^{3}}
\end{align*}

We expand $x_{\al t}$ to find
\begin{align*}
(19) & = \frac{1}{|x_{\al}|^{3}}\left((Q_x^{2} BR)_{\al} \cdot x_{\al}^{\perp}
+ \frac{Q_x^{2}\gamma}{2|x_{\al}|^{2}}x_{\al \al} \cdot x_{\al}^{\perp}\right)^2\\
&+ \underbrace{\frac{x_{\al \al} \cdot x_{\al}^{\perp}}{|x_{\al}|^{3}}\left((Q_x^{2} BR)_{\al} \cdot x_{\al}^{\perp}
+ \frac{Q_x^{2}\gamma}{2|x_{\al}|^{2}}x_{\al \al} \cdot x_{\al}^{\perp}\right)b_e}_{\text{error term: we incorporate it as }\tilde{\mathcal{E}}_{\al}^{2}}  - 2Q_x(Q_x)_{\al} BR \cdot x_{\al}^{\perp} \frac{\psi}{|x_{\al}|}\frac{x_{\al \al} \cdot x_{\al}^{\perp}}{|x_{\al}|^{3}}
\end{align*}

We denote
\begin{equation}
\label{spiral}
G_{x}(\al) = (Q_x^{2} BR)_{\al} \cdot x_{\al}^{\perp}
+ \frac{Q_x^{2}\gamma}{2|x_{\al}|^{2}}x_{\al \al} \cdot x_{\al}^{\perp}
\end{equation}

We claim that

$$ G_x(\al) = \text{NICE } + |x_{\al}|H(\da \psi)$$

that becomes
$$ (G_x(\al))^2 = \text{NICE}$$

Then
\begin{align*}
(19) = \text{NICE } - 2Q_x(Q_x)_{\al} BR \cdot x_{\al}^{\perp} \frac{\psi}{|x_{\al}|}\frac{x_{\al \al} \cdot x_{\al}^{\perp}}{|x_{\al}|^{3}} + \tilde{\mathcal{E}}_{\al}^{2}
\end{align*}

We write
\begin{align*}
G_x(\al) & = \underbrace{2Q_x(Q_x)_{\al} BR \cdot x_{\al}^{\perp}}_{\text{NICE, at the level of }x_{\al}}
+ \underbrace{Q_x^{2}\frac{1}{2\pi}\int \frac{(x_{\al}(\al) - x_{\al}(\al-\beta)) \cdot x_{\al}(\al)}{|x(\al)-x(\al-\beta)|^{2}}\gamma(\al-\beta)d\beta}_{\text{NICE, we use that }|x_{\al}|^{2} = A_x(t)} \\
& \underbrace{- Q_x^{2}\frac{1}{\pi}\int \frac{(x_{\al}(\al) - x_{\al}(\al-\beta)) \cdot x_{\al}(\al)}{|x(\al)-x(\al-\beta)|^{4}}(x(\al) - x(\al - \beta))(x_{\al}(\al) - x_{\al}(\al - \beta))\gamma(\al-\beta)d\beta}_{\text{NICE, we use that }|x_{\al}|^{2} \text{ only depends on time}} \\
& + \underbrace{Q_x^2 BR(x,\gamma_{\al}) \cdot x_{\al}^{\perp}}_{\text{Hilbert transform applied to }\gamma_{\al}} + \frac{Q_x^{2}\gamma}{2|x_{\al}|^{2}}x_{\al \al} \cdot x_{\al}^{\perp}
\end{align*}

Therefore
\begin{align*}
G_x(\al) & = \text{NICE } + |x_{\al}|Q_x^{2}H\left(\left(\frac{\gamma}{2|x_{\al}|}\right)_{\al}\right) + \frac{Q_x^{2}\gamma}{2|x_{\al}|^{2}}x_{\al \al} \cdot x_{\al}^{\perp} \\
& = \text{NICE } + |x_{\al}|H\left(\left(\frac{Q_x^{2}\gamma}{2|x_{\al}|}\right)_{\al}\right) + \frac{Q_x^{2}\gamma}{2|x_{\al}|^{2}}x_{\al \al} \cdot x_{\al}^{\perp} \\
& = \text{NICE } + |x_{\al}|H(\da \psi) + H\left((b_s |x_{\al}|^{2})_{\al}\right)+ \frac{Q_x^{2}\gamma}{2|x_{\al}|^{2}}x_{\al \al} \cdot x_{\al}^{\perp} \\
& = \text{NICE } + |x_{\al}|H(\psi_{\al}) - H\left((Q_{x}^{2}BR)_{\al} \cdot x_{\al}\right)+ \frac{Q_x^{2}\gamma}{2|x_{\al}|^{2}}x_{\al \al} \cdot x_{\al}^{\perp}
\end{align*}

\begin{align*}
&(Q_{x}^{2}BR)_{\al} \cdot x_{\al}  = \underbrace{2Q_x(Q_x)_{\al} BR \cdot x_{\al}}_{\text{NICE}}
+ Q_x^{2}\frac{1}{2\pi}\int \frac{(x_{\al}(\al) - x_{\al}(\al-\beta))^{\perp} \cdot x_{\al}(\al)}{|x(\al)-x(\al-\beta)|^{2}}\gamma(\al-\beta)d\beta \\
& = \underbrace{- Q_x^{2}\frac{1}{\pi}\int \frac{(x(\al) - x(\al-\beta))^{\perp} \cdot x_{\al}(\al)}{|x(\al)-x(\al-\beta)|^{4}}(x(\al) - x(\al - \beta))(x_{\al}(\al) - x_{\al}(\al - \beta))\gamma(\al-\beta)d\beta}_{\text{NICE, extra cancellation in }(x(\al) - x(\al - \beta))^{\perp} \cdot x_{\al}(\al)} \\
& + \underbrace{Q_x^{2}\frac{1}{2\pi}\int \frac{(x(\al) - x(\al-\beta))^{\perp} \cdot x_{\al}(\al)}{|x(\al)-x(\al-\beta)|^{2}}\gamma(\al-\beta)d\beta}_{\text{NICE, extra cancellation in }(x(\al) - x(\al - \beta))^{\perp} \cdot x_{\al}(\al)}
\end{align*}
This means that
\begin{align*}
(Q_{x}^{2}BR)_{\al} \cdot x_{\al} & = \text{NICE } + \frac{1}{2}H\left(Q_x^2\frac{\da^{2}x^{\perp} \cdot x_{\al}}{|x_{\al}|^{2}}\gamma\right)
\end{align*}
Taking Hilbert transforms:
\begin{align*}
-H\left((Q_{x}^{2}BR)_{\al} \cdot x_{\al}\right) & = \text{NICE } - \frac{1}{2}H^2\left(Q_x^2\frac{\da^{2}x^{\perp} \cdot x_{\al}}{|x_{\al}|^{2}}\gamma\right)
= \text{NICE } + \frac{1}{2}Q_x^2\frac{\da^{2}x^{\perp} \cdot x_{\al}}{|x_{\al}|^{2}}\gamma
\end{align*}

Using that $\da^{2} x^{\perp} \cdot x_{\al} = -\da^{2} x \cdot x_{\al}^{\perp}$ we are done. Thus (19) yields
\begin{align*}
\psi_{\al t} & = \text{NICE } - Q_x^{2}\left(\left(BR_t + \frac{\psi}{|x_{\al}|}BR_{\al}\right) \cdot x_{\al}^{\perp}\right.\\
&+\left. \frac{\gamma}{2|x_{\al}|^{2}}\left(x_{\al t} + \frac{\psi}{|x_{\al}|}x_{\al \al}\right)\cdot x_{\al}^{\perp}
+ \nabla P_{2}^{-1}(x) \cdot x_{\al}^{\perp}\right)\frac{x_{\al \al} \cdot x_{\al}^{\perp}}{|x_{\al}|^{3}} \\
& - Q_{x}^3\left(|BR|^{2} + \frac{b_s^2|x_{\al}|^{2}}{Q_x^{4}}+2b_s \frac{BR \cdot x_{\al}}{Q_x^{2}} - \frac{\psi^{2}}{Q_x^4}\right)\nabla Q(x) \cdot x_{\al}^{\perp} \frac{x_{\al \al} \cdot x_{\al}^{\perp}}{|x_{\al}|^{3}} \\
& + (Q_x \gamma + 2Q_x BR \cdot x_{\al})\nabla Q(x) \cdot x_{\al}^{\perp} \frac{x_{\al t} \cdot x_{\al}^{\perp}}{|x_{\al}|^{3}} \\
& \underbrace{- \left(\frac{Q_x^{3}(|BR|^{2})_{\al}}{|x_{\al}|} + \frac{(b_s^{2}|x_{\al}|)_{\al}}{Q_x}
+ \left(2b_s BR \cdot \frac{x_{\al}}{|x_{\al}|}\right)_{\al}Q_x\right)(Q_x)_{\al}}_{(20)} \\
& \underbrace{- 2Q_x(Q_x)_{\al} BR \cdot x_{\al}^{\perp} \frac{\psi}{|x_{\al}|} \frac{x_{\al \al} \cdot x_{\al}^{\perp}}{|x_{\al}|^{3}}}_{(21)}+ \mathcal{E}_{\al}^{2}, \quad \text{ where } \quad \mathcal{E}_{\al}^{2} = \tilde{\mathcal{E}}_{\al}^{1} + \tilde{\mathcal{E}}_{\al}^{2}
\end{align*}

For (20) we write
\begin{align*}
|x_t|^{2} & = Q_x^{4}|BR|^{2} + b_s^{2}|x_{\al}|^{2} + 2 Q_x^{2} b_s BR \cdot x_{\al} \\
& + \underbrace{b_e^{2}|x_{\al}|^{2} + f^{2} + 2Q_x^{2} BR \cdot x_{\al} b_e + 2b_s b_e |x_{\al}|^{2} + 2Q_x^{2} BR \cdot f
+ 2b_{s} x_{\al} \cdot f + 2b_e x_{\al} \cdot f}_{\text{ error terms } \tilde{\mathcal{E}}_{\al}^{3}} \\
\Rightarrow \frac{|x_t|^{2}}{Q_x|x_{\al}|} & =  \frac{Q_x^{3}|BR|^{2}}{|x_{\al}|} + \frac{b_s^{2}|x_{\al}|}{Q_x} + 2 Q_x b_s BR \cdot \frac{x_{\al}}{|x_{\al}|} + \frac{\tilde{\mathcal{E}}_{\al}^{3}}{Q_x|x_{\al}|} \\
\end{align*}

Now
\begin{align*}
(20) = \text{NICE } - \frac{(|x_t|^{2})_{\al}}{Q_x|x_{\al}|}(Q_x)_{\al} + \frac{\tilde{\mathcal{E}}_{\al}^{3}}{Q_x|x_{\al}|}(Q_x)_{\al}
\end{align*}

which means
\begin{align*}
(20)+(21) & = \text{NICE } - \frac{(|x_t|^{2})_{\al}}{Q_x|x_{\al}|}(Q_x)_{\al}
- 2Q_x(Q_x)_{\al} BR \cdot x_{\al}^{\perp} \frac{\psi}{|x_{\al}|} \frac{x_{\al \al} \cdot x_{\al}^{\perp}}{|x_{\al}|^{3}} + \frac{\tilde{\mathcal{E}}_{\al}^{3}}{Q_x|x_{\al}|}(Q_x)_{\al}
\end{align*}

We write
\begin{align*}
x_{\al t} & = \underbrace{(x_{\al t} \cdot x_{\al}) \frac{x_{\al}}{|x_{\al}|^{2}}}_{\text{only depends on }t}
+ (x_{\al t} \cdot x_{\al}^{\perp}) \frac{x_{\al}^{\perp}}{|x_{\al}|^{2}} \\
& = \left(B_x(t) + \frac{1}{2\pi}\int_{-\pi}^{\pi}f_{\beta} \cdot \frac{x_{\beta}}{|x_{\beta}|^{2}}d\beta\right)x_{\al} + \left((Q_{x}^{2} BR)_{\al} \cdot x_{\al}^{\perp} + b x_{\al \al} \cdot x_{\al}^{\perp} + f_{\al} \cdot x_{\al}^{\perp}\right)\frac{x_{\al}^{\perp}}{|x_{\al}|^{2}} \\
& = \left(B_x(t) + \frac{1}{2\pi}\int_{-\pi}^{\pi}f_{\beta} \cdot \frac{x_{\beta}}{|x_{\beta}|^{2}}d\beta\right)x_{\al} + \left((Q_{x}^{2} BR)_{\al} \cdot x_{\al}^{\perp} + b_s x_{\al \al} \cdot x_{\al}^{\perp}\right)\frac{x_{\al}^{\perp}}{|x_{\al}|^{2}}\\
&+ \left(b_e x_{\al \al} \cdot x_{\al}^{\perp} + f_{\al} \cdot x_{\al}^{\perp}\right)\frac{x_{\al}^{\perp}}{|x_{\al}|^{2}} \\
& = \left(B_x(t) + \frac{1}{2\pi}\int_{-\pi}^{\pi}f_{\beta} \cdot \frac{x_{\beta}}{|x_{\beta}|^{2}}d\beta\right)x_{\al} +
\underbrace{\left((Q_{x}^{2} BR)_{\al} \cdot x_{\al}^{\perp} + \frac{Q_x^{2}\gamma}{2|x_{\al}|^{2}} x_{\al \al} \cdot x_{\al}^{\perp}\right)}_{G_x(\al) \text{ as in }\eqref{spiral}}\frac{x_{\al}^{\perp}}{|x_{\al}|^{2}}\\
& - \frac{\psi}{|x_{\al}|}x_{\al \al} \cdot x_{\al}^{\perp} \frac{x_{\al}^{\perp}}{|x_{\al}|^{2}} \left(b_e x_{\al \al} \cdot x_{\al}^{\perp} + f_{\al} \cdot x_{\al}^{\perp}\right)\frac{x_{\al}^{\perp}}{|x_{\al}|^{2}}
\end{align*}

Writing $x_t = (Q_x^{2}BR) + b_s x_{\al} + b_e x_{\al} + f_{\al}$ we compute

\begin{align*}
x_{\al t} \cdot x_{\al} & = \underbrace{Q_{x}^{2} BR \cdot x_{\al}}_{\text{NICE}}\left(B_x(t) + \underbrace{\frac{1}{2\pi}\int_{-\pi}^{\pi}f_{\beta} \cdot \frac{x_{\beta}}{|x_{\beta}|^{2}}d\beta}_{\text{error}}\right) + \underbrace{G_x(\al) Q_{x}^{2} BR \cdot \frac{x_{\al}^{\perp}}{|x_{\al}|^{2}}}_{\text{NICE because }G_x \text{ is nice}} \\
& - \frac{\psi}{|x_{\al}|}x_{\al \al} \cdot x_{\al}^{\perp} Q_x^{2} BR \cdot \frac{x_{\al}^{\perp}}{|x_{\al}|^{2}}
+ Q_x^{2} BR \cdot \frac{x_{\al}^{\perp}}{|x_{\al}|^{2}}\left(b_e x_{\al \al} \cdot x_{\al}^{\perp} + \underbrace{f_{\al} \cdot x_{\al}^{\perp}}_{\text{error}}\right)\frac{x_{\al}^{\perp}}{|x_{\al}|^{2}} \\
& + \underbrace{b_s\left(B_x(t) + \frac{1}{2\pi}\int_{-\pi}^{\pi}f_{\beta} \cdot \frac{x_{\beta}}{|x_{\beta}|^{2}}d\beta\right)|x_{\al}|^{2}}_{\text{NICE}}
+ b_e\left(\underbrace{B_x(t)}_{\text{error}} + \frac{1}{2\pi}\int_{-\pi}^{\pi}f_{\beta} \cdot \frac{x_{\beta}}{|x_{\beta}|^{2}}d\beta\right)|x_{\al}|^{2} + \hat{\mathcal{E}}
\end{align*}

where $\hat{\mathcal{E}}$ is an error term. To simplify we write

$$ x_{\al t} \cdot x_{\al} = \text{NICE }- \frac{\psi}{|x_{\al}|}x_{\al \al} \cdot x_{\al}^{\perp} Q_x^{2} BR \cdot \frac{x_{\al}^{\perp}}{|x_{\al}|^{2}} + \text{ errors}$$

Setting the above formula in the expression of (20)+(21) allows us to find
$$ (20) + (21) = \text{NICE }+ \text{ errors }$$

This yields
\begin{align*}
\psi_{\al t} & = \text{NICE } - Q_x^{2}\left(\left(BR_t + \frac{\psi}{|x_{\al}|}BR_{\al}\right)
+ \frac{\gamma}{2|x_{\al}|^{2}}\left(x_{\al t} + \frac{\psi}{|x_{\al}|}x_{\al \al}\right)
+ \nabla P_{2}^{-1}(x)\right)\cdot x_{\al}^{\perp}\frac{x_{\al \al} \cdot x_{\al}^{\perp}}{|x_{\al}|^{3}} \\
& - Q_{x}^3\left(|BR|^{2} + \frac{b_s^2|x_{\al}|^{2}}{Q_x^{4}}+2b_s \frac{BR \cdot x_{\al}}{Q_x^{2}} - \frac{\psi^{2}}{Q_x^4}\right)\nabla Q(x) \cdot x_{\al}^{\perp} \frac{x_{\al \al} \cdot x_{\al}^{\perp}}{|x_{\al}|^{3}} \\
& + (Q_x \gamma + 2Q_x BR \cdot x_{\al})\nabla Q(x) \cdot x_{\al}^{\perp} \frac{x_{\al t} \cdot x_{\al}^{\perp}}{|x_{\al}|^{3}} + \mathcal{E}_{\al}^{3}
\end{align*}

being $\mathcal{E}_{\al}^{3}$ a new error term. We now complete the formula for $\si_x$ in \eqref{X} to find
\begin{align*}
\psi_{\al t} & = \text{NICE } - Q_x^{2}\si_{x} \frac{x_{\al \al} \cdot x_{\al}^{\perp}}{|x_{\al}|^{3}} \\
& + \underbrace{Q_x \left|BR + \frac{\gamma}{2|x_{\al}|^{2}}x_{\al}\right|^{2}
\nabla Q(x) \cdot x_{\al}^{\perp} \frac{x_{\al \al} \cdot x_{\al}^{\perp}}{|x_{\al}|^{3}}}_{(22)} \\
& + \underbrace{Q_{x}^3\left(-|BR|^{2} - \frac{b_s^2|x_{\al}|^{2}}{Q_x^{4}}-2b_s \frac{BR \cdot x_{\al}}{Q_x^{2}} + \frac{\psi^{2}}{Q_x^4}\right)\nabla Q(x) \cdot x_{\al}^{\perp} \frac{x_{\al \al} \cdot x_{\al}^{\perp}}{|x_{\al}|^{3}}}_{(23)} \\
& + \underbrace{(Q_x \gamma + 2Q_x BR \cdot x_{\al})\nabla Q(x) \cdot x_{\al}^{\perp} \frac{x_{\al t} \cdot x_{\al}^{\perp}}{|x_{\al}|^{3}}}_{(24)} + \mathcal{E}_{\al}^{3}
\end{align*}

Expanding
$$ \frac{\psi^{2}}{Q_x^{4}} = \frac{\gamma^{2}}{4|x_{\al}|^{2}} + \frac{b_s^{2}|x_{\al}|^{2}}{Q_x^{4}} - \frac{\gamma b_s}{Q_x^{2}}$$

we find

\begin{align*}
(22) + (23) = Q_{x}^3\left(\frac{\gamma^{2}}{2|x_{\al}|^{2}} + BR \cdot x_{\al} \frac{\gamma}{|x_{\al}|^{2}}-2b_s \frac{BR \cdot x_{\al}}{Q_x^{2}}
- \frac{\gamma b_s}{Q_x^{2}}\right)\nabla Q(x) \cdot x_{\al}^{\perp} \frac{x_{\al \al} \cdot x_{\al}^{\perp}}{|x_{\al}|^{3}}
\end{align*}

Writing
\begin{align*}
x_{\al t} \cdot x_{\al}^{\perp} = (Q_x^{2} BR)_{\al} x_{\al}^{\perp} + b_{s} x_{\al \al} \cdot x_{\al}^{\perp} + \text{ errors}
\end{align*}

we obtain that
\begin{align*}
(24)&  = \left(Q_x \gamma + 2Q_x BR \cdot x_{\al}\right)\nabla Q(x) \cdot x_{\al}^{\perp} \frac{(Q_x^{2} BR)_{\al} \cdot x_{\al}^{\perp}}{|x_{\al}|^{3}} \\
& + \left(Q_x \gamma + 2Q_x BR \cdot x_{\al}\right)\nabla Q(x) \cdot x_{\al}^{\perp} b_s\frac{x_{\al \al} \cdot x_{\al}^{\perp}}{|x_{\al}|^{3}}
+ \text{ errors}
\end{align*}

Thus
\begin{align*}
(22)+(23)&+(24)  = Q_{x}^3\left(\frac{\gamma^{2}}{2|x_{\al}|^{2}} + BR \cdot x_{\al} \frac{\gamma}{|x_{\al}|^{2}}
\right)\nabla Q(x) \cdot x_{\al}^{\perp} \frac{x_{\al \al} \cdot x_{\al}^{\perp}}{|x_{\al}|^{3}} \\
& + \left(Q_x \gamma + 2Q_x BR \cdot x_{\al}\right)\nabla Q(x) \cdot x_{\al}^{\perp} \frac{(Q_x^{2} BR)_{\al} \cdot x_{\al}^{\perp}}{|x_{\al}|^{3}}
+ \text{ errors} \\
& = Q_x \nabla Q(x) \cdot x_{\al}^{\perp} \left(\gamma + 2BR \cdot x_{\al}\right)\left(\frac{Q_x^{2}\gamma}{2|x_{\al}|^{2}}\frac{x_{\al \al} \cdot x_{\al}^{\perp}}{|x_{\al}|^{3}} + \frac{(Q_x^{2}BR)_{\al} \cdot x_{\al}^{\perp}}{|x_{\al}|^{3}}\right) + \text{ errors} \\
& = Q_x \nabla Q(x) \cdot x_{\al}^{\perp} \left(\gamma + 2BR \cdot x_{\al}\right)\frac{1}{|x_{\al}|^{3}}D_{x}(\al) + \text{ errors}\\
& = \text{NICE } + \text{ errors}
\end{align*}

Finally, we obtain

$$ \psi_{\al t} = \text{NICE}(x,\gamma,\psi) - Q_{x}^{2}\si_{x} \frac{x_{\al \al} \cdot x_{\al}^{\perp}}{|x_{\al}|^{3}}
+ \mathcal{E}_{\al}^{4}$$

For $\vp_{\al t}$ we find
$$ \vp_{\al t} = \text{NICE}(z,\om,\vp) - Q_{z}^{2}\si_{z} \frac{z_{\al \al} \cdot z_{\al}^{\perp}}{|z_{\al}|^{3}},$$

since we can apply the same methods as before to the equations with $f = g = 0$, which are satisfied by $(z,\om,\vp)$. Then:

\begin{align*}
II & = \int_{-\pi}^{\pi} \Lambda \da^{3} \mathcal{D} \cdot \da^{3} \mathcal{D}_{t} =
\int_{-\pi}^{\pi} \Lambda \da^{3} \mathcal{D}\left(\text{NICE}(z,\om,\vp) - \text{NICE}(x,\gamma,\psi)\right)d\al \\
& - \int_{-\pi}^{\pi}\Lambda \da^{3} \mathcal{D}\left(\da^{2}\left(Q_{z}^{2} \si_{z} \frac{z_{\al \al} \cdot z_{\al}^{\perp}}{|z_{\al}|^{3}} - Q_x^{2} \si_{x} \frac{x_{\al \al} \cdot x_{\al}^{\perp}}{|x_{\al}|^{3}}\right)\right)
- \int_{-\pi}^{\pi} \Lambda \da^{3} \mathcal{D} \mathcal{E}_{\al}^{4} d\al
\equiv II_{1} + II_{2} + II_{3}
\end{align*}

$$ II_1 \leq CP(E(t)) \quad \text{ because we are dealing with the NICE term}$$
$$ II_{3} \leq CP(E(t)) + c\delta(t) \quad \text{ because of the errors}$$

It remains to estimate $II_{2}$. We consider the most singular terms
\begin{align*}
II_2 & = II_{2,1} + II_{2,2} + II_{2,3} + \text{ l.o.t} \\
II_{2,1} & = - \int_{-\pi}^{\pi}\Lambda (\da^{3} \mathcal{D})\left((Q_{z}^{2})_{\al \al} \si_{z} \frac{z_{\al \al} \cdot z_{\al}^{\perp}}{|z_{\al}|^{3}} - (Q_x^{2})_{\al \al} \si_{x} \frac{x_{\al \al} \cdot x_{\al}^{\perp}}{|x_{\al}|^{3}}\right)d\al \\
II_{2,2} & = - \int_{-\pi}^{\pi}\Lambda (\da^{3} \mathcal{D})\left(Q_{z}^{2} \si_{z} \frac{\da^{4}z \cdot z_{\al}^{\perp}}{|z_{\al}|^{3}} - Q_x^{2} \si_{x} \frac{\da^{4}x \cdot x_{\al}^{\perp}}{|x_{\al}|^{3}}\right)d\al \\
II_{2,3} & = - \int_{-\pi}^{\pi}\Lambda (\da^{3} \mathcal{D})\left(Q_{z}^{2} \da^{2} \si_{z} \frac{z_{\al \al} \cdot z_{\al}^{\perp}}{|z_{\al}|^{3}} - Q_x^{2} \da^{2} \si_{x} \frac{x_{\al \al} \cdot x_{\al}^{\perp}}{|x_{\al}|^{3}}\right)d\al
\end{align*}
\begin{align*}
II_{2,1} & = - \int_{-\pi}^{\pi}H (\da^{3} \mathcal{D})\left((Q_{z}^{2})_{\al \al} \si_{z} \frac{z_{\al \al} \cdot z_{\al}^{\perp}}{|z_{\al}|^{3}} - (Q_x^{2})_{\al \al} \si_{x} \frac{x_{\al \al} \cdot x_{\al}^{\perp}}{|x_{\al}|^{3}}\right)_{\al}d\al\\
&\leq CP(E(t)) + c\delta(t) \quad \text{ as before}
\end{align*}
For $II_{2,2}$ we decompose further
\begin{align*}
II_{2,2} & = - S + \widetilde{II}_{2,2}, \quad \text{ where} \\
\widetilde{II}_{2,2} & = - \int_{-\pi}^{\pi}\Lambda (\da^{3} \mathcal{D})\left(Q_{z}^{2} \si_{z} \frac{\da^{4}x \cdot z_{\al}^{\perp}}{|z_{\al}|^{3}} - Q_x^{2} \si_{x} \frac{\da^{4}x \cdot x_{\al}^{\perp}}{|x_{\al}|^{3}}\right)d\al \\
S &= \int_{-\pi}^{\pi}Q_{z}^{2}\sigma_{z} \da^{4}D \cdot \frac{\da z^{\perp}(\al)}{|z_{\al}|^3}H(\da^{4} \Dcal)(\al) d\al
\end{align*}

We find that
\begin{align*}
\widetilde{II}_{2,2} \leq CP(E(t)) + c\delta(t)
\end{align*}

and $-S$ cancels out with $S$. We are done with $II_{2,2}$. We write

\begin{align*}
II_{2,3} & = \int_{-\pi}^{\pi}H (\da^{3} \mathcal{D})\left(Q_{z}^{2} \da^{3} \si_{z} \frac{z_{\al \al} \cdot z_{\al}^{\perp}}{|z_{\al}|^{3}} - Q_x^{2} \da^{3} \si_{x} \frac{x_{\al \al} \cdot x_{\al}^{\perp}}{|x_{\al}|^{3}}\right)d\al + \quad \text{ l.o.t}
\end{align*}

We claim that
\begin{align}
\label{casitapaco}
Q_x^2 \da^3 \si_{x} = |x_{\al}|H(\da^{3} \psi_{t}) - b_s |x_{\al}| H(\da^{4} \psi) + \text{ errors } + \text{NICE}(x,\gamma,\psi)
\end{align}

In the local existence we get

\begin{align*}
Q_z^2 \da^3 \si_{z} = |z_{\al}|H(\da^{3} \vp_{t}) - c |z_{\al}| H(\da^{4} \vp) + \text{NICE}(z,\om,\vp)
\end{align*}

This implies
\begin{align*}
II_{2,3} & = II_{2,3,1} + II_{2,3,2} + II_{2,3,3} + II_{2,3,4} \\
II_{2,3,1} & = \int_{-\pi}^{\pi}H (\da^{3} \mathcal{D})\left(|z_{\al}|H(\da^{3} \vp_t) \frac{z_{\al \al} \cdot z_{\al}^{\perp}}{|z_{\al}|^{3}} - |x_{\al}| H(\da^{3} \psi_t)\frac{x_{\al \al} \cdot x_{\al}^{\perp}}{|x_{\al}|^{3}}\right)d\al \\
II_{2,3,2} & = -\int_{-\pi}^{\pi}H (\da^{3} \mathcal{D})\left(c|z_{\al}|H(\da^{4} \vp) \frac{z_{\al \al} \cdot z_{\al}^{\perp}}{|z_{\al}|^{3}} - b_s|x_{\al}| H(\da^{4} \psi)\frac{x_{\al \al} \cdot x_{\al}^{\perp}}{|x_{\al}|^{3}}\right)d\al \\
II_{2,3,3} & = \int_{-\pi}^{\pi}H (\da^{3} \mathcal{D})\left(\text{NICE}(z,\om,\vp) \frac{z_{\al \al} \cdot z_{\al}^{\perp}}{|z_{\al}|^{3}} - \text{NICE}(x,\gamma,\psi)\frac{x_{\al \al} \cdot x_{\al}^{\perp}}{|x_{\al}|^{3}}\right)d\al \\
II_{2,3,4} & = -\int_{-\pi}^{\pi}H (\da^{3} \mathcal{D})\frac{x_{\al \al} \cdot x_{\al}^{\perp}}{|x_{\al}|^{3}}d\al + \text{ errors}
\end{align*}

It is easy to find
$$ II_{2,3,4} \leq CP(E(t)) + c\delta(t), \quad \text{ error terms}$$
$$ II_{2,3,3} \leq CP(E(t)), \quad \text{ l.o.t}$$

In $II_{2,3,2}$ we split further:

\begin{align*}
II_{2,3,2} & = II_{2,3,2}^{1} + II_{2,3,2}^{2} \\
II_{2,3,2}^{1} & = -\int_{-\pi}^{\pi}H (\da^{3} \mathcal{D})|z_{\al}|H(\da^{4} \Dcal) \frac{z_{\al \al} \cdot z_{\al}^{\perp}}{|z_{\al}|^{3}}d\al \\
II_{2,3,2}^{2} & = -\int_{-\pi}^{\pi}H (\da^{3} \mathcal{D})H(\da^{4} \psi)\left(c|z_{\al}| \frac{z_{\al \al} \cdot z_{\al}^{\perp}}{|z_{\al}|^{3}} - b_s|x_{\al}| \frac{x_{\al \al} \cdot x_{\al}^{\perp}}{|x_{\al}|^{3}}\right)d\al
\end{align*}

Then
$$ II_{2,3,2}^{1} = \frac{1}{2}\int_{-\pi}^{\pi}|H (\da^{3} \mathcal{D})|^{2}\left(|z_{\al}| \frac{z_{\al \al} \cdot z_{\al}^{\perp}}{|z_{\al}|^{3}}\right)_{\al}d\al \leq CP(E(t))$$

For $II_{2,3,2}^{2}$:

$$ II_{2,3,2}^{2} = -\int_{-\pi}^{\pi}\Lambda^{1/2}(\da^{3}\psi)\Lambda^{1/2}\left(H (\da^{3} \mathcal{D})\left(c|z_{\al}| \frac{z_{\al \al} \cdot z_{\al}^{\perp}}{|z_{\al}|^{3}} - b_s|x_{\al}| \frac{x_{\al \al} \cdot x_{\al}^{\perp}}{|x_{\al}|^{3}}\right)\right)d\al \leq CP(E(t))$$

It remains
\begin{align*}
II_{2,3,1} & = II_{2,3,1}^{1} + II_{2,3,1}^{2} \\
II_{2,3,1}^{1} &  = \int_{-\pi}^{\pi}H (\da^{3} \mathcal{D})H(\da^{3} \Dcal_t) \frac{z_{\al \al} \cdot z_{\al}^{\perp}}{|z_{\al}|^{2}} d\al \\
II_{2,3,1}^{2} & = \int_{-\pi}^{\pi}H (\da^{3} \mathcal{D}) \underbrace{H(\da^{3}\psi_t)}_{\text{approx. sol.}}\left(\frac{z_{\al \al} \cdot z_{\al}^{\perp}}{|z_{\al}|^{2}} - \frac{x_{\al \al} \cdot x_{\al}^{\perp}}{|x_{\al}|^{2}}\right)d\al \\
\end{align*}

Then
$$ II_{2,3,1}^{1} \leq CP(E(t)) + c\delta(t)$$

At this point we remember that we had to deal with
$$ II = \int_{-\pi}^{\pi} \Lambda (\da^{3} \mathcal{D})\da^{3} \mathcal{D}_{t}d\al$$

so in $II_{2,3,1}^{1}$ we find one derivative less (or $1/2$ derivatives less) and this shows that we can bound

$$ II_{2,3,1}^{1} \leq CP(E(t)) + c\delta(t)$$

by brute force. It remains to show claim \eqref{casitapaco}. We remember

\begin{align*}
Q_{x}^{2}\si_{x} & = Q_x^{2}\left(BR_t + \frac{\psi}{|x_{\al}|}BR_{\al}\right) \cdot x_{\al}^{\perp}
+ \frac{Q_x^2\gamma}{2|x_{\al}|^{2}}\left(x_{\al t} + \frac{\psi}{|x_{\al}|}x_{\al \al}\right) \cdot x_{\al}^{\perp} \\
& + \underbrace{Q_x^{3}\left|BR + \frac{\gamma}{2|x_{\al}|^{2}}x_{\al}\right|^{2}\nabla Q(x) \cdot x_{\al}^{\perp}}_{\text{this term is in }H^3 \text{ so it is NICE}} + \underbrace{Q_x^{2} \nabla P_{2}^{-1}(x) \cdot x_{\al}^{\perp}}_{\text{this term is also in }H^3}
\end{align*}

We write
\begin{align*}
\frac{Q_x^2\gamma}{2|x_{\al}|^{2}}&\left(x_{\al t} + \frac{\psi}{|x_{\al}|}x_{\al \al}\right) \cdot x_{\al}^{\perp}\\
& = \frac{Q_x^2\gamma}{2|x_{\al}|^{2}}\left((Q_x BR)_{\al} \cdot x_{\al}^{\perp} + b_s x_{\al \al} \cdot x_{\al}^{\perp}
+ b_e x_{\al \al} \cdot x_{\al}^{\perp} + f_{\al} \cdot x_{\al}^{\perp} + \frac{\psi}{|x_{\al}|}x_{\al \al} \cdot x_{\al}^{\perp}\right) \\
& = \frac{Q_x^2\gamma}{2|x_{\al}|^{2}}\left((Q_x BR)_{\al} \cdot x_{\al}^{\perp} + \left(b_s + \frac{\psi}{|x_{\al}|}\right) x_{\al \al} \cdot x_{\al}^{\perp}\right)\\
& + \frac{Q_x^2\gamma}{2|x_{\al}|^{2}}\left(b_e x_{\al \al} \cdot x_{\al}^{\perp} + f_{\al} \cdot x_{\al}^{\perp}\right) \\
& = \frac{Q_x^2\gamma}{2|x_{\al}|^{2}}\left((Q_x BR)_{\al} \cdot x_{\al}^{\perp} + \frac{Q_x^{2}\gamma}{2|x_{\al}|^{2}} x_{\al \al} \cdot x_{\al}^{\perp}\right) + \text{ errors} \\
& = \frac{Q_x^2\gamma}{2|x_{\al}|^{2}}G_x(\al) + \text{ errors} = \text{NICE } + \text{ errors}
\end{align*}

Finally, the most singular terms in $Q_x^{2} \si_{x}$ are

$$ L = Q_x^{2} BR_t \cdot x_{\al}^{\perp} + \frac{Q_x^{2}\psi}{|x_{\al}|}BR_{\al} \cdot x_{\al}^{\perp}$$

We take 3 derivatives and consider the most dangerous characters:
\begin{align*}
L & = M_1 + M_2 + M_3 + \text{ l.o.t} \\
M_1 & = Q_x^{2} BR(x,\da^{3} \gamma_t) \cdot x_{\al}^{\perp} + \frac{Q_x^{2}\psi}{|x_{\al}|}BR(x,\da^{4}\gamma) \cdot x_{\al}^{\perp} \\
M_2 & = Q_x^{2} \frac{1}{2\pi}\int_{-\pi}^{\pi}\frac{(\da^{3} x_{t}(\al) - \da^{3}x_{t}(\al-\beta))\cdot x_{\al}(\al)}{|x(\al) - x(\al-\beta)|^{2}}
\gamma(\al-\beta)d\beta \\
& + \frac{Q_x^{2}\psi}{|x_{\al}|} \frac{1}{2\pi}\int_{-\pi}^{\pi}\frac{(\da^{4} x(\al) - \da^{4}x(\al-\beta))\cdot x_{\al}(\al)}{|x(\al) - x(\al-\beta)|^{2}}
\gamma(\al-\beta)d\beta \\
M_3 & = -\frac{Q_x^2}{\pi}\int_{-\pi}^{\pi}\frac{\Delta_{\beta} x(\al) \cdot x_{\al}(\al)}{|\Delta_{\beta}x(\al)|^{4}}\Delta_{\beta} x(\al) \cdot
\Delta_{\beta} \da^{3}x_t(\al) \gamma(\al-\beta)d\beta \\
&  -\frac{\psi Q_x^2}{|x_{\al}|\pi}\int_{-\pi}^{\pi}\frac{\Delta_{\beta} x(\al) \cdot x_{\al}(\al)}{|\Delta_{\beta}x(\al)|^{4}}\Delta_{\beta} x(\al) \cdot
\Delta_{\beta} \da^{4}x(\al) \gamma(\al-\beta)d\beta
\end{align*}

In $M_2$ we find
$$ M_2 = \frac{Q_x^2\gamma}{2|x_{\al}|^{2}}\Lambda(\da^{3} x_t \cdot x_{\al})
+ \frac{Q_x^{2}\psi\gamma}{|x_{\al}|^{3}}\Lambda(\da^{4}x \cdot x_{\al}) + \text{ l.o.t}$$

For the second term we use the usual trick
$$ \da^{4} x \cdot x_{\al} = -3\da^{3} x \cdot x_{\al \al}$$

For the first term we remember that
\begin{align*}
& |x_{\al}|^{2} = A(t) \Rightarrow x_{\al} \cdot x_{\al t} = \frac{1}{2}A'(t) \Rightarrow (x_{\al} \cdot x_{\al t})_{\al} = 0 \\
& \Rightarrow x_{\al \al} \cdot x_{\al t} + x_{\al} \cdot x_{\al \al t} = 0 \Rightarrow
x_{\al \al \al} \cdot x_{\al t} + 2x_{\al \al} \cdot x_{\al \al t} + x_{\al} \cdot x_{\al \al \al t} = 0\\
& \Rightarrow x_{\al} \cdot x_{\al \al \al t} = - 2 x_{\al \al} \cdot x_{\al \al t} - x_{\al \al \al} \cdot x_{\al t}
\end{align*}

This allows us to control $M_2$. For $M_3$ we find

$$ M_3 = -\frac{Q_x^2\gamma}{|x_{\al}|^{2}}\Lambda(x_{\al} \cdot \da^{3} x_t)
- \frac{Q_x^{2}\psi\gamma}{|x_{\al}|^{3}}\Lambda(x_{\al} \cdot \da^{4}x) + \text{ l.o.t}$$

so it can be estimated in the same way as $M_2$. There remains $M_1$.

$$ M_1  = Q_x^{2} BR(x,\da^{3} \gamma_t) \cdot x_{\al}^{\perp} + \frac{Q_x^{2}\psi}{|x_{\al}|}BR(x,\da^{4}\gamma) \cdot x_{\al}^{\perp}$$

Using that $\Delta_{\beta} x^{\perp}(\al) \cdot x_{\al}^{\perp}(\al) = \Delta_{\beta} x(\al) \cdot x_{\al}(\al)$ we find

\begin{align}
\label{ouroboros}
M_1 = \frac{Q_x^{2}}{2}H(\da^{3}\gamma_{t}) + \frac{Q_x^{2} \psi}{2|x_{\al}|}H(\da^{4}\gamma) + \text{ l.o.t}
\end{align}

We compute
\begin{align}
\frac{Q_x^{2}}{2}H(\da^{3}\gamma_{t}) & = H\left(\da^{3}\left(\frac{Q_x^{2}\gamma}{2}\right)_{t}\right) + \text{ NICE} \nonumber \\
& = H(\da^{3}(|x_{\al}|\psi)_{t}) + H(\da^{3}(|x_{\al}|b_s)_{t}) + \text{ NICE} \nonumber \\
& = |x_{\al}|H(\da^{3}\psi_{t}) + H(\da^{2}\partial_t(-(Q_x^{2} BR)_{\al} \cdot x_{\al})) + \text{ NICE} \label{anchor}
\end{align}

We compute the most singular term in
\begin{align*}
\da^{2}\partial_t(-(Q_x^{2} BR)_{\al} \cdot x_{\al}) & = -\frac{Q_x^{2}}{2\pi}\int_{-\pi}^{\pi}\frac{(\da^{3}x_t(\al) - \da^{3}x_t(\al-\beta))^{\perp} \cdot x_{\al}(\al)}{|\Delta_{\beta} x(\al)|^{2}}\gamma(\al-\beta)d\beta \\
& + \underbrace{\frac{Q_x^{2}}{\pi}\int_{-\pi}^{\pi}\frac{(\Delta_{\beta} x(\al))^{\perp} \cdot x_{\al}}{|\Delta_{\beta}x(\al)|^{4}}\Delta_{\beta} x(\al) \Delta_{\beta} \da^{3} x_{t}(\al) \gamma(\al-\beta)d\beta}_{\text{extra cancellation}} \\
& - \underbrace{\frac{Q_x^{2}}{2\pi}\int_{-\pi}^{\pi}\frac{(\Delta_{\beta} x(\al))^{\perp} \cdot x_{\al}}{|\Delta_{\beta}x(\al)|^{2}}\da^{3} \gamma_{t}(\al-\beta)d\beta}_{\text{extra cancellation}} + \text{ l.o.t. } + \text{ NICE}
\end{align*}

This shows that

\begin{align*}
\da^{2}\partial_t(-(Q_x^{2} BR)_{\al} \cdot x_{\al}) = -\frac{Q_x^{2}\gamma}{2|x_{\al}|^{2}}\Lambda(\da^{3} x_{t}^{\perp} \cdot x_{\al})  + \text{ l.o.t. } + \text{ NICE}
\end{align*}

That gives
\begin{align*}
\da^{2}\partial_t(-(Q_x^{2} BR)_{\al} \cdot x_{\al}) = -\Lambda\left(\frac{Q_x^{2}\gamma}{2|x_{\al}|^{2}}\da^{3} x_{t}^{\perp} \cdot x_{\al}\right)  + \text{ l.o.t. } + \text{ NICE}
\end{align*}

which implies
\begin{align*}
H(\da^{2}\partial_t(-(Q_x^{2} BR)_{\al} \cdot x_{\al})) &= \da\left(\frac{Q_x^{2}\gamma}{2|x_{\al}|^{2}}\da^{3} x_{t}^{\perp} \cdot x_{\al}\right)  + \text{ l.o.t. } + \text{ NICE} \\
&= -\frac{Q_x^{2}\gamma}{2|x_{\al}|^{2}}\da\left(\da^{3} x_{t} \cdot x_{\al}^{\perp}\right) + \text{ NICE}
\end{align*}

Plugging the above formula in \eqref{anchor} we find that
\begin{align*}
\frac{Q_x^{2}}{2}H(\da^{3}\gamma_{t}) & = |x_{\al}|H(\da^{3}\psi_{t}) - \frac{Q_x^{2}\gamma}{2|x_{\al}|^{2}}\da\left(\da^{3} x_{t} \cdot x_{\al}^{\perp}\right) + \text{ NICE} \\
& = |x_{\al}|H(\da^{3}\psi_{t}) - \frac{Q_x^{2}\gamma}{2|x_{\al}|^{2}}\da\left(\da^{3}(Q_x^{2} BR) \cdot x_{\al}^{\perp}\right)
- \frac{Q_x^{2}\gamma}{2|x_{\al}|^{2}}\da\left(b_s \da^{4} x \cdot x_{\al}^{\perp}\right) + \text{ l.o.t}\\
& + \text{ NICE} + \text{ errors}\\
\end{align*}

As we did before, in $\da(\da^{3}(Q_x^{2}BR) \cdot x_{\al}^{\perp})$, the most dangerous term is given by $Q_{x}^{2}\frac{1}{2}H(\da^{4}\gamma)$, the tangential terms appear, which implies
\begin{align*}
\da(\da^{3}(Q_x^{2}BR) \cdot x_{\al}^{\perp}) & = Q_x^{2}\frac{1}{2}H(\da^{4}\gamma) + \text{ NICE}
\end{align*}
and therefore
\begin{align*}
\frac{Q_x^{2}}{2}H(\da^{3}\gamma_{t}) = |x_{\al}|H(\da^{3}\psi_{t}) - \frac{Q_x^{2}\gamma}{2|x_{\al}|^{2}}\frac{Q_{x}^{2}}{2}H(\da^{4}\gamma)
- \frac{Q_x^{2}\gamma}{2|x_{\al}|^{2}}b_s\da\left(\da^{4} x \cdot x_{\al}^{\perp}\right) + \text{ NICE} + \text{ errors}\\
\end{align*}

We use \eqref{ouroboros} to find
\begin{align*}
 \frac{Q_x^{2}}{2}&H(\da^{3}\gamma_{t}) + \frac{Q_x^{2}\psi}{2|x_{\al}|}H(\da^{4}\gamma)
 \\
 & = |x_{\al}|H(\da^{3}\psi_{t}) - \frac{Q_x^{2}}{2}b_sH(\da^{4}\gamma)
- \frac{Q_x^{2}\gamma}{2|x_{\al}|^{2}}b_s\da\left(\da^{4} x \cdot x_{\al}^{\perp}\right) + \text{ NICE} + \text{ errors} \\
& = |x_{\al}|H(\da^{3}\psi_{t}) - b_s|x_{\al}|H\left(\da^{4}\left(\frac{Q_x^{2}\gamma}{2|x_{\al}|}\right)\right)
- \frac{Q_x^{2}\gamma}{2|x_{\al}|^{2}}b_s\da\left(\da^{4} x \cdot x_{\al}^{\perp}\right) + \text{ NICE} + \text{ errors} \\
& = |x_{\al}|H(\da^{3}\psi_{t}) - b_s|x_{\al}|H\left(\da^{4}\psi\right) - b_s|x_{\al}|H(\da^{4}(b_s|x_{\al}|))
- \frac{Q_x^{2}\gamma}{2|x_{\al}|^{2}}b_s\da\left(\da^{4} x \cdot x_{\al}^{\perp}\right) \\
&+ \text{ NICE} + \text{ errors} \\
\end{align*}

We will show that

$$ - b_s|x_{\al}|H(\da^{4}(b_s|x_{\al}|))
- \frac{Q_x^{2}\gamma}{2|x_{\al}|^{2}}b_s\da\left(\da^{4} x \cdot x_{\al}^{\perp}\right)$$

is NICE and then we are done.

\begin{align*}
- b_s|x_{\al}|H(\da^{4}(b_s|x_{\al}|))&- \frac{Q_x^{2}\gamma}{2|x_{\al}|^{2}}b_s\da\left(\da^{4} x \cdot x_{\al}^{\perp}\right)
 = - b_sH(\da^{4}(b_s|x_{\al}|^{2}))- \frac{Q_x^{2}\gamma}{2|x_{\al}|^{2}}b_s\da\left(\da^{4} x \cdot x_{\al}^{\perp}\right) \\
& = b_sH(\da^{3}((Q_x^{2} BR)_{\al} \cdot x_{\al}))- \frac{Q_x^{2}\gamma}{2|x_{\al}|^{2}}b_s\da\left(\da^{4} x \cdot x_{\al}^{\perp}\right)
\end{align*}

We repeat the calculation for dealing with the most dangerous terms in
\begin{align*}
\da^{3}((Q_x^{2} BR)_{\al} \cdot x_{\al}) = \Lambda\left(\da^{4} x^{\perp} \cdot x_{\al} \frac{\gamma Q_{x}^{2}}{2|x_{\al}|^{2}}\right) + \text{ l.o.t}
\end{align*}

In the l.o.t we use that $\Delta_{\beta} x^{\perp}(\al) \cdot x(\al) $ gives an extra cancellation. We find that
\begin{align*}
b_sH(\da^{3}&((Q_x^{2} BR)_{\al} \cdot x_{\al}))  - \frac{Q_x^{2}\gamma}{2|x_{\al}|^{2}}b_s\da\left(\da^{4} x \cdot x_{\al}^{\perp}\right) \\
& = b_sH(\Lambda\left(\da^{4} x^{\perp} \cdot x_{\al} \frac{\gamma Q_{x}^{2}}{2|x_{\al}|^{2}}\right))- \frac{Q_x^{2}\gamma}{2|x_{\al}|^{2}}b_s\da\left(\da^{4} x \cdot x_{\al}^{\perp}\right) + \text{ NICE} \\
& = -b_s\da\left(\da^{4} x^{\perp} \cdot x_{\al} \frac{\gamma Q_{x}^{2}}{2|x_{\al}|^{2}}\right)- \frac{Q_x^{2}\gamma}{2|x_{\al}|^{2}}b_s\da\left(\da^{4} x \cdot x_{\al}^{\perp}\right) + \text{ l.o.t} + \text{ NICE}
\end{align*}

Using that $\da^{4} x^{\perp} \cdot x_{\al} = - \da^{4} x \cdot x_{\al}^{\perp} $ we are done.

\section*{Acknowledgements}

All the authors were partially supported by the grant MTM2011-26696 (Spain)
 and ICMAT Severo Ochoa project SEV-2011-0087. AC was partially supported by the ERC grant 307179-GFTIPFD.
 CF was supported by NSF grant DMS-09-0104.

\bibliographystyle{abbrv}
\bibliography{references}

\begin{tabular}{ll}
\textbf{Angel Castro} &  \\
{\small Departamento de Matem\'aticas} &\\
{\small Universidad Aut\'onoma de Madrid } & \\
{\small Instituto de Ciencias Matem\'aticas-CSIC} &\\
{\small Campus de Cantoblanco} & \\
{\small Email: angel\_castro@icmat.es} & \\
   & \\
\textbf{Diego C\'ordoba} &  \textbf{Charles Fefferman}\\
{\small Instituto de Ciencias Matem\'aticas} & {\small Department of Mathematics}\\
{\small Consejo Superior de Investigaciones Cient\'ificas} & {\small Princeton University}\\
{\small C/ Nicol\'{a}s Cabrera, 13-15} & {\small 1102 Fine Hall, Washington Rd, }\\
{\small Campus Cantoblanco UAM, 28049 Madrid} & {\small Princeton, NJ 08544, USA}\\
{\small Email: dcg@icmat.es} & {\small Email: cf@math.princeton.edu}\\
 & \\
\textbf{Francisco Gancedo} &  \textbf{Javier G\'omez-Serrano}\\
{\small Departamento de An\'alisis Matem\'atico} & {\small Department of Mathematics}\\
{\small Universidad de Sevilla} & {\small Princeton University}\\
{\small C/ Tarfia, s/n } & {\small 1102 Fine Hall, Washington Rd,} \\
{\small Campus Reina Mercedes, 41012 Sevilla}  & {\small Princeton, NJ 08544, USA} \\
{\small Email: fgancedo@us.es} & {\small Email: jg27@math.princeton.edu}\\
\end{tabular}

\end{document}